\documentclass{elsarticle}
\usepackage[left=1in,top=1in,right=1in,bottom=1in]{geometry}
\usepackage{graphicx}
\usepackage{placeins}
\usepackage{subfig}
\usepackage{amsmath}
\usepackage{amssymb}
\usepackage{algorithm}%
\usepackage{algorithmicx}%
\usepackage{algpseudocode}%
\usepackage[english]{babel}
\usepackage{amsthm}
\usepackage{url}
\newtheorem*{remark}{Remark}

\usepackage{xcolor}

\begin{document}

\begin{frontmatter}

\title{Pressure-Stabilized Fixed-Stress Iterative Solutions of Compositional Poromechanics}

\author[1]{Ryan M. Aronson\corref{cor1}}
\ead{rmaronso@stanford.edu}

\author[2]{Nicola Castelletto}

\author[3]{François P. Hamon}

\author[2]{Joshua A. White}

\author[1]{Hamdi A. Tchelepi}

\cortext[cor1]{Corresponding author}

\affiliation[1]{organization={Stanford University},
country={United States}}

\affiliation[2]{organization={Atmospheric, Earth, and Energy Division, Lawrence Livermore National Laboratory},
country={United States}}

\affiliation[3]{organization={TotalEnergies},
country={Denmark}}

\journal{arXiv}



\begin{abstract}
    We consider the numerical behavior of the fixed-stress splitting method for coupled poromechanics as undrained regimes are approached. We explain that pressure stability is related to the splitting error of the scheme, not the fact that the discrete saddle point matrix never appears in the fixed-stress approach. This observation reconciles previous results regarding the pressure stability of the splitting method. Using examples of compositional poromechanics with application to geological CO$_2$ sequestration, we see that solutions obtained using the fixed-stress scheme with a low order finite element-finite volume discretization which is not inherently inf-sup stable can exhibit the same pressure oscillations obtained with the corresponding fully implicit scheme. Moreover, pressure jump stabilization can effectively remove these spurious oscillations in the fixed-stress setting, while also improving the efficiency of the scheme in terms of the number of iterations required at every time step to reach convergence. 
\end{abstract}


\end{frontmatter}

\section{Introduction}

The focus of this paper is the development of numerical schemes for solving poromechanical problems in which undrained, or nearly undrained, regions exist. These are regions where fluid cannot easily flow within the pore space of the solid skeleton, so that pore pressure changes cannot rapidly dissipate. As a practical example of growing interest, undrained behavior frequently appears in simulations of CO$_2$ sequestration, where storage reservoirs are usually bounded above and below by low-permeability sealing formations. These seals prevent the buoyant CO$_2$ from travelling beyond the reservoir in a process known as structural trapping \cite{ipcc}. Accurate numerical modeling of flow and mechanics in both the drained reservoir and undrained seals is therefore very important to ensure the safety of the injection process.  One must ensure that the seals retain integrity under loading, and that stress perturbations do not activate large faults in the vicinity \cite{rutqvist2012geomechanics}. 


The complexity of undrained regions is that, when combined with nearly incompressible materials, the governing equations fundamentally change in character. In particular, it is well known that the governing equations reduce to a saddle point structure \cite{zienkiewicz1990static}, very similar to the incompressible Stokes equations of fluid mechanics. Thus for these types of problems we must consider the inf-sup stability of the discretization spaces used to represent the discrete mechanics and flow degrees of freedom. Since this reduction only occurs when considering coupled flow and mechanics with nearly incompressible materials in undrained regimes -- which to this point have not been commonly considered in standard reservoir simulation applications -- it is common in practice to ignore the inf-sup condition when selecting a discretization for simulations of large reservoirs. A common choice is to use a piecewise linear, nodal finite element scheme for the mechanics equations and a cell-centered, piecewise constant finite volume scheme for the flow equations. This is a convenient choice as it is the lowest-order conforming discretization pair, meaning it is computationally efficient for the large-scale simulations needed to model basin-scale CO$_2$ storage. Moreover, the choice of using a lowest-order finite volume scheme for the flow problem results in a scheme which conserves mass locally on each element and is monotone even in the presence of discontinuous material properties. However, these spaces are not inf-sup stable, meaning that in undrained regions spurious pressure oscillations will occur, just as incompressible flow. Within the poromechanical community, many methods have been developed to smooth these spurious pressure modes for various discretization types \cite{truty2006stabilized, aguilar2008numerical, preisig2011stabilization, wan2003stabilized, white2008stabilized, berger2015stabilized, camargo2021macroelement,frigo2021efficient}.

When considering pressure stabilization, most of the existing works in the literature consider a fully implicit, or monolithic, approach to solving the coupled equations of flow and mechanics. In these methods the flow and mechanics unknowns are solved for simultaneously at every time step and these schemes are attractive as they can be made unconditionally stable \cite{garipov2018unified}. However, implementation of a fully implicit simulation framework requires the development of highly sophisticated multiphysics solvers and preconditioners for the resulting linear systems, like those in \cite{bui2020scalable, bui2021multigrid}. 

Given these challenges, a large number of alternative techniques for coupling the mechanics and flow equations together have been proposed in the literature \cite{dean2006comparison}. These techniques are usually iterative, in the sense that they solve the flow and mechanics equations separately and pass data or solutions back and forth. The main advantage of this approach is that existing, highly optimized solvers for single physics applications can be used, thus allowing for fast development and improved flexibility. Within iteratively coupled techniques, methods can be further grouped by how accurately the coupling is captured. For example, methods can be loosely coupled \cite{park1983stabilization, settari2001advances}, where the subproblems are solved separately and there is no attempt to ensure the coupled problem is being accurately approximated at every time step. Another class includes sequential implicit methods. There exist numerous variations, including some which can be shown to be convergent to the fully implicit, or monolithic, solution \cite{mikelic2013convergence}. In this work we focus on the fixed-stress splitting method, which has become very popular in coupled reservoir simulation applications \cite{settari1998coupled, kim2011stability, kim2011stability_spe, mikelic2013convergence, castelletto2015accuracy, garipov2018unified}, though we note that there are other methods resulting from different choices of solver order and solution passing schemes \cite{kim2011stability_drained}. 

There are very few studies in the existing literature, however, which consider the performance of these sequential coupling strategies on problems in which undrained regions occur. The change in character of the governing equations means that usual analyses of splitting or fractional step methods does not apply directly \cite{janenko1971method}. The work by Yoon and Kim \cite{yoon2018spatial} claims that the fixed-stress scheme is not subject to the inf-sup condition, as the discrete saddle point system is never formed nor inverted in the scheme, only matrices which are full rank by construction. In a later work by Storvik \textit{et al.} \cite{storvik2019optimization}, the authors noted that an inf-sup stable discretization is needed in order to develop their results on optimizing the convergence rate of the fixed-stress split, and they noted that because of this assumption, theirs was the first proof of convergence of the fixed-stress method for undrained problems. This seems to contradict the previous claims of \cite{yoon2018spatial}, and reconciling these results by providing a more in depth study of the relation between splitting schemes and pressure stability is the main goal of the current work. 

This manuscript is ordered as follows. First we review the governing equations of compositional multiphase flow coupled with geomechanics as well as the spatial and time discretization. We then discuss the fixed-stress method, paying special attention to choices regarding sequential iterations and convergence when defining the method which are important in developing understanding of the possible pressure stability of the method. The following section discusses the relation between pressure stability and splitting schemes in detail and appeals to methods used in incompressible flow as examples for clarity. It is here that we form the main point of this work, namely that it is the error which comes from the splitting method which can possibly provide an implicit pressure regularizing effect, not merely the lack of any (nearly) singular saddle point systems being solved. Through this argument it is much easier to see that, if the fixed-stress scheme includes iterations to convergence within each time step, the inf-sup condition must still be satisfied by the underlying spatial discretization in order to avoid spurious pressure oscillations, even though all of the matrices being inverted in the scheme are full rank. In the case of a spatial discretization which is not inherently inf-sup stable, we can use the same stabilizations developed for fully implicit methods to smooth spurious oscillations, and in this work we utilize pressure jump stabilization \cite{berger2015stabilized, camargo2021macroelement, aronson2023pressure}. Our numerical results on problems including undrained burden regions confirm that, without stabilization, the pressure field obtained with the fixed-stress method contains the same oscillations present in fully implicit results. We believe that this is the first time this has been explicitly shown in the field. Including pressure stabilization not only smooths these oscillations, but also decreases the number of sequential (outer) iterations required at any time step to reach convergence, matching similar results in \cite{storvik2018optimization} obtained with inherently stable elements.

\section{Governing Equations}

Following the theory first introduced by Biot \cite{biot1941general} and since extended to more complex scenarios \cite{coussy2004poromechanics, wang2000theory}, the coupled poromechanical behavior of a solid skeleton saturated by fluid components is fully described by mechanical force and mass conservation equations. The mechanical equilibrium of the solid skeleton can be written as

\begin{equation}
    \nabla \cdot \boldsymbol{\sigma}' - b \nabla p + \rho \mathbf{g} = \mathbf{0}, 
\end{equation}

\noindent where $\boldsymbol{\sigma}'$ is the effective stress tensor, $b$ is the Biot coefficient, $p$ is the fluid pressure (note we ignore capillarity effects in this work), $\rho = (1-\phi)\rho_s + \phi \rho_f$ is the homogenized density, $\rho_s$ is the density of the solid, $\rho_f$ is the fluid density, $\phi$ is the porosity, and $\mathbf{g}$ is the gravitational acceleration.
The effective stress and fluid pressure combine to define the total stress via

\begin{equation}
    \boldsymbol{\sigma} = \boldsymbol{\sigma}' - bp\mathbf{I},
\end{equation}

\noindent where $\mathbf{I}$ is the identity tensor.

In this work we assume small strains and linear elastic solids, thus the effective stress is determined by the mechanical strain $\epsilon$ and the elasticity modulus tensor $\mathbb{C}$:

\begin{equation}
    \boldsymbol{\sigma}' = \mathbb{C} : \boldsymbol{\epsilon},
\end{equation}

\noindent and the strain is simply the symmetric gradient of the displacement $\mathbf{u}$:

\begin{equation}
    \boldsymbol{\epsilon} = \nabla^s \mathbf{u} = \frac{1}{2}(\nabla \mathbf{u} + \nabla^T \mathbf{u}).
\end{equation}

We assume this solid skeleton is fully saturated with $n_c$ components flowing in $n_p$ miscible phases. Thus the mechanical equilibrium equation for the solid skeleton is augmented with conservation equations for each component, which are of the form

\begin{equation}
    \left( \dot{\overline{ \phi m_c }} \right) + \nabla \cdot \left( \sum_{\ell = 1}^{n_p} x_{c\ell} \rho_\ell \mathbf{v}_\ell \right) +  \sum_{\ell = 1}^{n_p} x_{c\ell} \rho_\ell q_\ell = 0 \qquad c = 1, \ldots, n_c.
    \label{eq:comp_mass_balance}
\end{equation}

\noindent Here we use $m_c = \sum_{\ell = 1}^{n_p} x_{c\ell} \rho_\ell s_\ell$ to represent the component mass of component $c$, and for phase $\ell$, $x_{c\ell}$ is the mass fraction of component $c$ in this phase, $\rho_{\ell}$ is the density, $s_{\ell}$ is the saturation, $\mathbf{v}_{\ell}$ is the Darcy velocity, and $q_{\ell}$ is the source term used to represent injection and production. The dependence of the porosity $\phi$ on pressure and stress is given by

\begin{equation}
    \phi = \phi_0 + \frac{(b-\phi_0)(1-b)}{K_{dr}}(p - p_0) + b\epsilon_v,
\end{equation}

\noindent where $\phi_0$ and $p_0$ are the reference porosity and pressure, respectively, $K_{dr}$ is the drained bulk modulus of the skeleton, and $\epsilon_v$ is the volumetric strain $\epsilon_v = \boldsymbol{\epsilon}_{ii}$, where we have adopted the notation that repeated indices represent summations. 

To fully close the system, we start with the volume constraint requiring that the saturations of each phase sum to one everywhere:

\begin{equation}
    \sum_{\ell=1}^{n_p} s_{\ell} = 1. 
\end{equation}

\noindent We also have the mass constraint, i.e. the sum of the component fractions in a given phase must be equal to one:

\begin{equation}
    \sum_{c=1}^{n_c} x_{c\ell} = 1  \qquad \ell = 1, \ldots, n_p.
\end{equation}

\noindent We assume thermodynamic equilibrium, which is expressed in terms of the fugacity of component $c$ in phase $\ell$: $f_{c,\ell}(p, T, x_{c,\ell})$ as 

\begin{equation}
    f_{c\ell}(p, T, x_{c\ell}) - f_{ck}(p, T, x_{ck}) = 0 \qquad \forall \ell \neq k, \quad c = 1, ... , n_c,
\end{equation}

\noindent where $T$ is the temperature (though we assume isothermal conditions in this work). Finally the velocity of each phase is determined using the relative permeability of phase $\ell$, $k_{r\ell}$, the viscosity of phase $\ell$, $\mu_\ell$, and the permeability tensor $\mathbf{k}$ via the multiphase extension of Darcy's law:

\begin{equation}
    \mathbf{v}_\ell = -\frac{k_{r\ell}}{\mu_\ell} \mathbf{k} \cdot \left( \nabla p - \rho_\ell \mathbf{g} \right).
    \label{eq:darcy}
\end{equation}

In certain portions of the following sections, we will also refer to the governing equations of single-phase flow for simplicity. In this case we set $n_p = n_c = 1$, and then we can rewrite the mass balance equation (Equation \eqref{eq:comp_mass_balance}) as 

\begin{equation}
    \dot{ m} + \nabla \cdot \left( \rho \mathbf{v} \right) = \rho q,
\end{equation}

\noindent where $m$ is the fluid content increment. Defining $M$ to be Biot's modulus and using the constitutive relation 

\begin{equation}
    m = b \nabla \cdot \mathbf{u} + \frac{1}{M}p,
\end{equation}

\noindent as well as the relation between effective and total stress, we can rewrite the mass conservation equation as

\begin{equation}
    \left( \frac{1}{M} + \frac{b^2}{K_{dr}} \right) \dot{p} + \frac{b}{K_{dr}}
    \dot{\sigma_v} + \nabla \cdot \mathbf{v} =  q. 
    \label{eq:single_phase_mass}
\end{equation}

\noindent Note that now the only flow unknown is the fluid pressure, and if the skeleton is linear elastic and has constant permeability, the coupled problem also becomes linear. This is the formulation used in \cite{kim2011stability} when defining the fixed-stress splitting. 

\subsection{Undrained Conditions}

In this work we are particularly interested in the behavior of numerical schemes used to approximate the solutions of the differential equations above as undrained conditions are approached. A problem is considered undrained if the product $||\mathbf{k}|| \delta t$ is small or zero, which occurs when either small time steps $\delta t$ are taken relative to the loading rate, or the skeleton is essentially impermeable. In this work we focus on cases with small $||\mathbf{k}||$, typical in burden or seal formations surrounding permeable aquifers or reservoirs. Note that many studies focus on the saddle-point structure that appears at the initial time of the simulation or which occurs with small time steps \cite{murad1992improved, murad1994stability}. In these cases, the errors made by employing unstable element pairs usually decay in time. This is not necessarily the case when the permeability is near zero, however, and the pressure oscillations can actually in grow in time in this case \cite{wan2003stabilized}. 

For much of our analysis, understanding the numerical behavior of discretization schemes in the case of single-phase flow is sufficient. If the solid and fluid materials are nearly incompressible, the undrained mass equation reduces to an incompressibility constraint, which prevents the solid skeleton from deforming volumetrically:

\begin{equation}
    \nabla \cdot \dot{\mathbf{u}} = 0.
\end{equation}

This reduction was also shown for multiphase poromechanical systems in \cite{camargo2021macroelement, aronson2023pressure}. Finally, we emphasize that, while the exact saddle point system structure is only recovered if the permeability magnitude is exactly zero and the materials are exactly incompressible, if undrained conditions are merely approached the results obtained using discretizations that are not inf-sup stable may be polluted by spurious pressure modes. In particular, we will also show that the brine phase present in simulations of CO$_2$ sequestration is close enough to incompressible for the instability to arise.

\section{Discretization}

There are, of course, many different possible discretization schemes for coupled poromechanics including Galerkin finite elements for both flow and mechanics \cite{murad1994stability}, mixed finite elements (which can provide local mass conservation) \cite{ferronato2010fully}, finite volume methods for both flow and mechanics \cite{nordbotten2014cell}, as well as discontinuous and enriched Galerkin methods \cite{liu2004discontinuous,choo2018enriched}. 

In this work we focus on a mixed discretization scheme where the mechanics unknowns are discretized using nodal, piecewise linear basis functions, while the flow unknowns are approximated using a cell centered, piecewise constant basis. Note that, for compositional simulations, one must decide which set of variables are to be the primary unknowns and which are determined implicitly from the primary unknowns. In this work we elect to use the overall composition formulation, in which the primary unknowns are selected to be the pressure and the overall component densities \cite{voskov2012comparison}.

The spatial domain $\Omega$ is approximated using a set of disjoint elements $K$ so that $\Omega \approx \bigcup K$, and in this work we mainly focus on hexahedral and tetrahedral elements. We denote the mechanics test space as 

\begin{equation}
    \boldsymbol{\mathcal{V}}^h = \{ \boldsymbol{\psi}^h : \boldsymbol{\psi}^h \in [C^0(\Omega)]^3, \textup{ } \boldsymbol{\psi}^h|_{K} \in [\mathbb{X}_1(K)], \textup{ } \boldsymbol{\psi}^h = \mathbf{0} \textup{ on } \Gamma_D^u\},
\end{equation}

\noindent where the function space $\mathbb{X}_1(K)$ is the space of trilinear nodal basis functions $\mathbb{Q}_1(K)$ for hexahedra and the space of linear nodal basis functions $\mathbb{P}_1(K)$ for tetrahedral elements. Then the discrete, weak form of the mechanics equation can be stated at time $t_{n+1}$ as, $\forall \boldsymbol{\psi}^h \in \boldsymbol{\mathcal{V}}^h$

\begin{equation}
    \int_\Omega \nabla^s \boldsymbol{\psi}^h : (\boldsymbol{\sigma}')^h_{n+1} d\Omega - \int_\Omega \nabla \cdot \boldsymbol{\psi}^h b p^h_{n+1} d \Omega - \int_\Omega \boldsymbol{\psi}^h \cdot \rho^h_{n+1} \mathbf{g} d \Omega - \int_{\Gamma_N^u} \boldsymbol{\psi}^h \cdot \mathbf{T}_{n+1} d \Gamma = 0, 
    \label{eq:discrete_mech}
\end{equation}

\noindent where $\mathbf{T}_{n+1}$ is the applied boundary traction at time $t_{n+1}$.

Now we move on to consider the flow equations, which are discretized in space using a cell-centered, piecewise constant finite volume scheme as mentioned above. Note that the flow equations are time dependent, thus we must also select a time discretization scheme. In this work we utilize backward Euler integration. The discretized mass conservation equations can be written as, $\forall K \in \Omega$

\begin{equation}
    V^K\frac{\left( \phi m_c \right)_{n+1}^K -  \left( \phi m_c \right)_{n}^K } {t_{n+1} - t_{n}}  + \sum_{f\in \partial K} F_{c, n+1}^f = V^K \left( \sum_{\ell = 1}^{n_p} x_{c\ell} \rho_\ell q_\ell \right)^K_{n+1} \qquad c = 1, \ldots, n_c.
    \label{eq:discrete_mass}
\end{equation}

\noindent Here we have used $V$ to denote the volume of a cell, $F_{c, n+1}^f$ to denote the numerical flux of component $c$ across face $f$ at time $t_{n+1}$, and the superscript $K$ to denote the cell-centered values in cell $K$. In this work, we model source terms by specifying an injected CO$_2$ mass rate and do not rely on standard Peaceman-type well models.

To compute the numerical fluxes $F_{c, n+1}^f$ we use a linear two-point flux approximation (TPFA), common in the reservoir simulation community \cite{aziz1979petroleum}. The flux across face $f$ between cells $K$ and $L$ is computed using a single point upwind mass term and a potential difference:

\begin{equation}
    F_{c, n+1}^f = \left( \sum_{\ell = 1}^{n_p} x_{c\ell} \rho_\ell \frac{k_{r,\ell}}{\mu_\ell} \right)^{upw}_{n+1} \Upsilon^f \left( p^{K}_{n+1} - p^{L}_{n+1} + g^f_{c, n+1} \right), 
\end{equation}

\noindent where $g^f_{c, n+1}$ represents the gravitational potential difference across the face $f$ and the upwind direction is determined by the sign of the full potential difference term. The factor $\Upsilon^f$ is the transmissibility coefficient computed from the mesh geometry and the permeability \cite{aziz1979petroleum, white2019two}. Note that, despite its popularity, the TPFA scheme can be non-convergent upon spatial refinement on non K-orthogonal grids, and one must resort to more complex schemes to achieve provable convergence \cite{terekhov2017cell}.

Fully implicit approaches solve the coupled nonlinear problem resulting from Equations \eqref{eq:discrete_mech} and \eqref{eq:discrete_mass} monolithically at every time step, using, for example, a Newton-Raphson approach. This involves the formation and inversion of the Jacobian matrix at every nonlinear iteration. Considering an incompressible, single-phase case for simplicity, as undrained conditions are approached, this matrix reduces to the block form

\begin{equation}
    J = 
    \begin{bmatrix}
    A &  -B^T \\
    B& C 
    \end{bmatrix}
    \rightarrow
    \begin{bmatrix}
    A &  -B^T \\
    B& 0 
    \end{bmatrix},
    \label{eq:discrete_saddle}
\end{equation}

\noindent as the $C$-block contains the accumulation and flux contributions, which are zero for incompressible and impermeable materials. This reduced matrix is a standard form of a discrete saddle-point system, where the $A$-block represents the discrete stiffness operator while $B$ is defined by the discrete divergence.

When stabilized methods have been developed in the past for undrained poromechanical problems, the invertibility of this matrix is of key importance. One interpretation of the discrete inf-sup condition is that of a requirement that this Jacobian matrix is full rank. Indeed, stabilized methods can remedy this by inserting a block matrix into the lower-right block which makes the coupled Jacobian matrix full rank. However, as we shall see in the following sections, the inf-sup condition may still be required when employing sequential schemes like the fixed-stress split. 

\section{The Fixed-Stress Split}

Rather than solve the system monolithically, in this work we consider the fixed-stress sequential approach in which flow and mechanics are solved separately. The key idea of the fixed-stress split is the assumption that the rate of the mean, total volumetric stress, given by 

\begin{equation}
    \sigma_v = \frac{1}{3}\boldsymbol{\sigma}_{ii},
\end{equation}

\noindent is assumed constant between successive iterations between the flow and mechanics solvers. The fixed-stress method was proven to be stable and convergent in \cite{kim2011stability, mikelic2013convergence}, converges to the monolithic solution if iterations between the mechanics and flow solves are repeated within individual time steps, and can be interpreted as a preconditioned Richardson iteration or a preconditioner for the monolithic method \cite{castelletto2015accuracy}. 

The choice of inclusion of sequential, or outer, iterations within individual time steps is important to the analysis of the spatial stability of the fixed-stress split in the undrained limit, and so we clearly distinguish the sequentially iterative and non-iterative methods. In the case where sequential iterations are not included, the fixed-stress assumption is \cite{kim2011stability}

\begin{equation}
     \sigma^{n+1}_{v} - \sigma^{n}_{v} = \sigma^{n}_{v} - \sigma^{n-1}_{v}.
     \label{eq:FS_split}
\end{equation}

\noindent This means that the flow problem is solved with the total stress determined by an extrapolation using an explicit time derivative approximation. The structure of a time step of the coupled scheme is then described by Algorithm \ref{alg:FS_noiter}. 

\begin{algorithm}
\caption{Fixed-Stress (non-iterative)}
\label{alg:FS_noiter}
\begin{algorithmic}[1]
        \State Begin new time step (advancing from $t_{n}$ to $t_{n+1}$)
        \State Solve flow problem assuming $\sigma^{n+1}_{v} = 2 \sigma^{n}_{v} - \sigma^{n-1}_{v}$
        \State Solve mechanics problem with updated flow variables
        \State Advance to next time step
 \end{algorithmic}
\end{algorithm}

For single-phase flow problems coupled with mechanics (which we recall is a linear problem in our case), we can describe the fixed-stress splitting using a block linear algebra formulation. In particular, by using the assumption in Equation \eqref{eq:FS_split} in conjunction with Equation \eqref{eq:single_phase_mass}, the fixed-stress method replaces the fully implicit discrete system by \cite{kim2011stability}

\begin{equation}
    \begin{bmatrix}
    A & -B^T \\
    0 & C+R 
    \end{bmatrix}
    \begin{bmatrix}
    u^{n+1}-u^n  \\
    p^{n+1}-p^n  
    \end{bmatrix}
    =
    \begin{bmatrix}
    Q_u \\
    Q_p
    \end{bmatrix}
    +
    \begin{bmatrix}
    0 & 0 \\
    -B & R 
    \end{bmatrix}
    \begin{bmatrix}
    u^{n}-u^{n-1}  \\
    p^{n}-p^{n-1}  
    \end{bmatrix}.
    \label{eq:block_FS}
\end{equation}

\noindent In these expressions we note that the matrix $R$ is the diagonal matrix with entries equal to $b^2 / K_{dr}$, while the remaining blocks are simply those that result from the discretization of the differential operators, with $Q_u$ and $Q_p$ representing forcing terms. 

In the case where sequential iterations are performed within each time step, the fixed-stress approximation is given by

\begin{equation}
    (\sigma^{n+1}_{v} - \sigma^{n}_{v})_{k+1} = (\sigma^{n+1}_{v} - \sigma^{n}_{v})_k \implies \sigma^{n+1}_{v, k+1} = \sigma^{n+1}_{v, k},
\end{equation}

\noindent where $k$ represents the sequential iteration index. In this case it is common to assume the initial state at a time step is equal to the converged state at the previous time step, and no explicit extrapolation is needed in contrast to the non-iterated scheme \cite{garipov2018unified, storvik2019optimization}. The form of a time step of this scheme is summarized in Algorithm \ref{alg:FS_iter}.

\begin{algorithm}
\caption{Fixed-Stress (iterative)}
\label{alg:FS_iter}
\begin{algorithmic}[1]
        \State Begin new time step (advancing from $t_{n}$ to $t_{n+1}$)
        \State Set unknowns at time $t_{n+1}$ to converged result from time $t_{n}$
        \State Set $k = 0$
        \While{Not converged}
        \State Solve flow problem assuming $\sigma^{n+1}_{v, k+1} = \sigma^{n+1}_{v, k}$ 
        \State Solve mechanics problem with updated flow variables
        \State Update $k = k+1$
        \EndWhile
        \State Advance to next time step
 \end{algorithmic}
\end{algorithm}

\noindent Note that the solutions of the flow and mechanics problems listed in steps 5 and 6 of Algorithm \ref{alg:FS_iter} may also require iterative solution procedures if the corresponding problems are nonlinear. We also remark that the sequentially iterated scheme can also be cast in a block linear algebraic way in the case of single phase flow, in which case the result is the same as in Equation \eqref{eq:block_FS}, but with the unknowns being the $k+1$ and $k$-th estimates of the solution change at time $t_{n+1}$ \cite{castelletto2015accuracy}. 

In the development of the above methods, the key assumption was the conservation of the physical total stress between sequential iterations or time steps. However, one can also choose to conserve a fictitious total stress-like quantity when designing iterative schemes, or equivalently use a different bulk modulus within the definition of the physical total stress, and these schemes may achieve faster convergence \cite{mikelic2013convergence, castelletto2015accuracy}. For instance, one can choose to hold the quantity $K_{dr} \nabla \cdot \mathbf{u} - \alpha b p$ constant, where $\alpha$ is a free parameter. This results in entries of $R$ in the block system (Equation \eqref{eq:block_FS}) which are multiplied by $\alpha$, and convergence of the fixed-stress method is guaranteed for 
$\alpha$ values greater than $1/2$. Indeed, it has been found experimentally and analytically that a value in the range between $1/2$ and $1$ typically produces the fastest sequential convergence, in terms of number of iterations within a time step \cite{storvik2019optimization, both2019numerical}. We will typically use $\alpha = 1.0$ in our studies unless otherwise stated, corresponding to conservation of the standard physical total volumetric stress.  

\section{Pressure Stability of Splitting Schemes}

While the spatial stability of the monolithic method in the undrained limit has been well studied at this point, the behavior of the fixed-stress scheme in this same regime is less known. As stated previously \cite{yoon2018spatial} it is claimed that the fixed-stress split provides pressure stability in the undrained limit regardless of the underlying discretization pair for mechanics and flow degrees of freedom. We note that it seems that this work employs the fixed-stress method without sequential iterations within time steps. Meanwhile, \cite{storvik2019optimization} explains that the assumption of an inf-sup stable element pair is necessary to prove the convergence of the fixed-stress method in the undrained, incompressible regime. In the preprint of this same manuscript, these authors also report slower convergence of the sequential iterations within time steps when the splitting method with equal-order, piecewise linear interpolations (which are not inf-sup stable) is employed for undrained problems \cite{storvik2018optimization}.

To reconcile these results, we start by again considering \cite{yoon2018spatial}, which claims that the fixed-stress scheme is always inf-sup stable. The crux of the authors' argument is that, in the discrete setting, the saddle point matrix never arises and all of the coefficient matrices being inverted all have full rank. Thus the inf-sup criterion need never be satisfied. However, the same matrices are inverted in the iterative scheme, which would thus suggest that an inf-sup stable discretization is not needed in the iterative scheme either, contradicting the results in \cite{storvik2019optimization, storvik2018optimization}. This rank-based argument was also presented within the incompressible flow community when studying pressure projection schemes (which also perform an operator split much like the fixed-stress methodology), but this understanding was shown to be subtly wrong \cite{guermond1998stability, guermond2006overview}. Instead, these works show that it is the error introduced by the splitting scheme which can have a regularizing effect on the pressure, smoothing the spurious, oscillatory pressure modes associated with a lack of inf-sup stability. However, there are certainly splitting schemes which do not provide sufficient pressure stabilization and will admit spurious pressure oscillations, even though the discrete matrices are nonsingular at every step of the algorithm. For concreteness, we provide examples by considering the stability analysis of two simple splitting schemes presented in \cite{guermond1998stability, guermond2006overview} for incompressible flows, to elucidate the relation between operator splits and inf-sup stability before returning to the specific case of the fixed-stress split. 

\subsection{Splitting Schemes for Incompressible Flows}

In the case of incompressible Stokes flow, the Chorin splitting (or pressure-projection scheme) \cite{chorin1968numerical, chorin1969convergence} is performed by first generating an intermediate velocity field considering the unsteady and viscous diffusion terms, then the pressure which makes this velocity field incompressible is computed via a Poisson problem, and finally the velocity at the end of the time step is computed using the gradient of the new pressure. It is now known \cite{guermond1998stability, guermond2006overview} that the splitting scheme results in a modified saddle point structure where the mass equation is of the form

\begin{equation}
    \nabla \cdot \tilde{\mathbf{u}}^{n+1} = \delta t \Delta p^{n+1},
\end{equation}

\noindent where $\tilde{\mathbf{u}}$ is the intermediate velocity field. In words, the splitting error has introduced a diffusive pressure effect, and it is this mechanism which dampens spurious pressure oscillations if an inf-sup unstable element pair is used for spatial discretization. Indeed, the term $\delta t \Delta p$ essentially resembles terms that are introduced in stabilized methods, such as the fluid pressure Laplacian technique for coupled poromechanics \cite{truty2006stabilized, aguilar2008numerical, preisig2011stabilization}. However, as detailed in \cite{guermond1998stability, guermond2006overview}, the presence of the $\delta t$ factor means that if the time step size is selected to be too small, pressure oscillations will reappear. Indeed, convergence of this projection method can only be proved for sufficiently large time step sizes if the underlying discretization spaces are not inf-sup stable \cite{badia2007convergence}. 

A slightly more advanced splitting is the incremental pressure projection scheme \cite{goda1979multistep}. The idea of this scheme is simply the inclusion of the explicit pressure gradient when computing the intermediate velocity field. The projection step is then used to determine the increment in pressure between the old and new time steps, rather than simply the new pressure as in the Chorin splitting. In this case the modified mass equation can be written as \cite{guermond1998stability}

\begin{equation}
    \nabla \cdot \tilde{\mathbf{u}}^{n+1} = \delta t \Delta (p^{n+1} - p^n).
\end{equation}

\noindent The incremental scheme can be shown to provide possible second order accuracy in time for velocity, in contrast to the first order accuracy of the non-incremental scheme \cite{guermond2006overview}. However, it has also been shown that the splitting error in this case is not sufficient to smooth spurious pressure oscillations for any choice of time step size \cite{guermond1998stability}. Note that the discrete matrices are still full rank in this scheme and the discrete saddle point matrix is not inverted at any point, and thus this argument is not sufficient to prove pressure stability of a split scheme. 

\subsection{Splitting Errors of the Fixed-Stress Splitting}

With this interpretation of the stability of splitting schemes, we now return to consideration of the fixed-stress splitting in undrained conditions. We start with the non-sequentially-iterated scheme, and focus on the case of nearly incompressible, single phase flow for simplicity, though the conclusions should be general in the case of multiphase extensions. As undrained conditions are approached, the scheme shown in Equation \eqref{eq:block_FS} reduces to

\begin{equation}
    \begin{bmatrix}
    A & -B^T \\
    0 & R 
    \end{bmatrix}
    \begin{bmatrix}
    u^{n+1}-u^n  \\
    p^{n+1}-p^n  
    \end{bmatrix}
    =
    \begin{bmatrix}
    Q_u \\
    Q_p  
    \end{bmatrix}
    +
    \begin{bmatrix}
    0 & 0 \\
    -B & R 
    \end{bmatrix}
    \begin{bmatrix}
    u^{n}-u^{n-1}  \\
    p^{n}-p^{n-1}  
    \end{bmatrix},
    \label{eq:block_FS_undrained}
\end{equation}

\noindent  where the $C$ contribution has vanished, as in the fully implicit case (Equation \eqref{eq:discrete_saddle}). The saddle point matrix corresponding to the discretization of the fully implicit system has been replaced by a matrix which is clearly full rank, as $A$ and $R$ are positive definite. However, we now know that this is not sufficient to conclude that the scheme will be pressure stable. Using simple algebraic manipulations and ignoring any forcing terms for simplicity, we can rewrite the fixed-stress scheme as

\begin{equation}
    \begin{bmatrix}
    A & -B^T \\
    -B & 0 
    \end{bmatrix}
    \begin{bmatrix}
    u^{n+1}-u^n  \\
    p^{n+1}-p^n  
    \end{bmatrix}
    =
    \begin{bmatrix}
    0  \\
    -B(u^{n+1} - 2u^n + u^{n-1}) + R(p^{n+1} - 2p^n + p^{n-1})  
    \end{bmatrix}.
    \label{eq:block_FS_saddle}
\end{equation}

\noindent Much like splitting schemes for incompressible flow, Equation \eqref{eq:block_FS_saddle} shows that the fixed-stress scheme can be written in the form of a saddle-point problem with a modified right hand side, representing the splitting error. One must show that the error term $-B(u^{n+1} - 2u^n + u^{n-1}) + R(p^{n+1} - 2p^n + p^{n-1})$ produces a diffusive effect on the pressure if one wishes to claim that the fixed-stress method is pressure stable, as in \cite{yoon2018spatial}. It is likely, however, that any stabilizing effect will depend on the problem parameters, similar to the Chorin splitting for incompressible flow. 

Performing a similar study on the sequentially iterated scheme yields a clear understanding of its lack of pressure stability. In this case the scheme can be written as

\begin{equation}
\begin{split}
    &\begin{bmatrix}
    A & -B^T \\
    -B & 0 
    \end{bmatrix}
    \begin{bmatrix}
    (u^{n+1}-u^n)_{k+1}  \\
    (p^{n+1}-p^n)_{k+1}  
    \end{bmatrix}
    \\ &=
    \begin{bmatrix}
    0  \\
    -B[(u^{n+1}-u^n)_{k+1}-(u^{n+1}-u^n)_{k}] + R[(p^{n+1}-p^n)_{k+1}-(p^{n+1}-p^n)_k]  
    \end{bmatrix}.
\end{split}
\label{eq:block_FS_saddle_iter}
\end{equation}

\noindent As the goal of the sequentially iterative scheme is to iterate until convergence, we see that there is no splitting error on the right hand side of the mass equation which could provide a pressure stabilization effect. This is consistent with the notion that the iterated scheme converges to the solution of the fully implicit method, and once again we emphasize that all of the matrices inverted during the iterated scheme are still nonsingular. 

\subsection{A Simple Numerical Example}

To confirm our understandings of the pressure stability of the fixed-stress scheme, we now consider the results from a simple numerical test, namely an undrained modification of the single-phase Barry-Mercer problem \cite{barry1999exact, camargo2021macroelement}. We consider the unit square domain in 2D, discretized with a 10-by-10 grid. We set the skeleton Young's modulus to 10 kPa and the Poisson ratio to $0.2$. We consider an incompressible fluid with unit viscosity. The permeability is assumed isotropic with value $1.0 \times 10^{-12}$ m$^2$ and we simulate 10 time steps of size 10 seconds, making the problem essentially undrained. Finally, we apply a forcing in the flow equation of the form $\sin(\pi t /100)$ to the cell centered at $(0.35, 0.15)$. 

\begin{figure}
\centering
\subfloat[Fully implicit]{\label{sfig:undrained_fim}\includegraphics[width=0.5\textwidth]{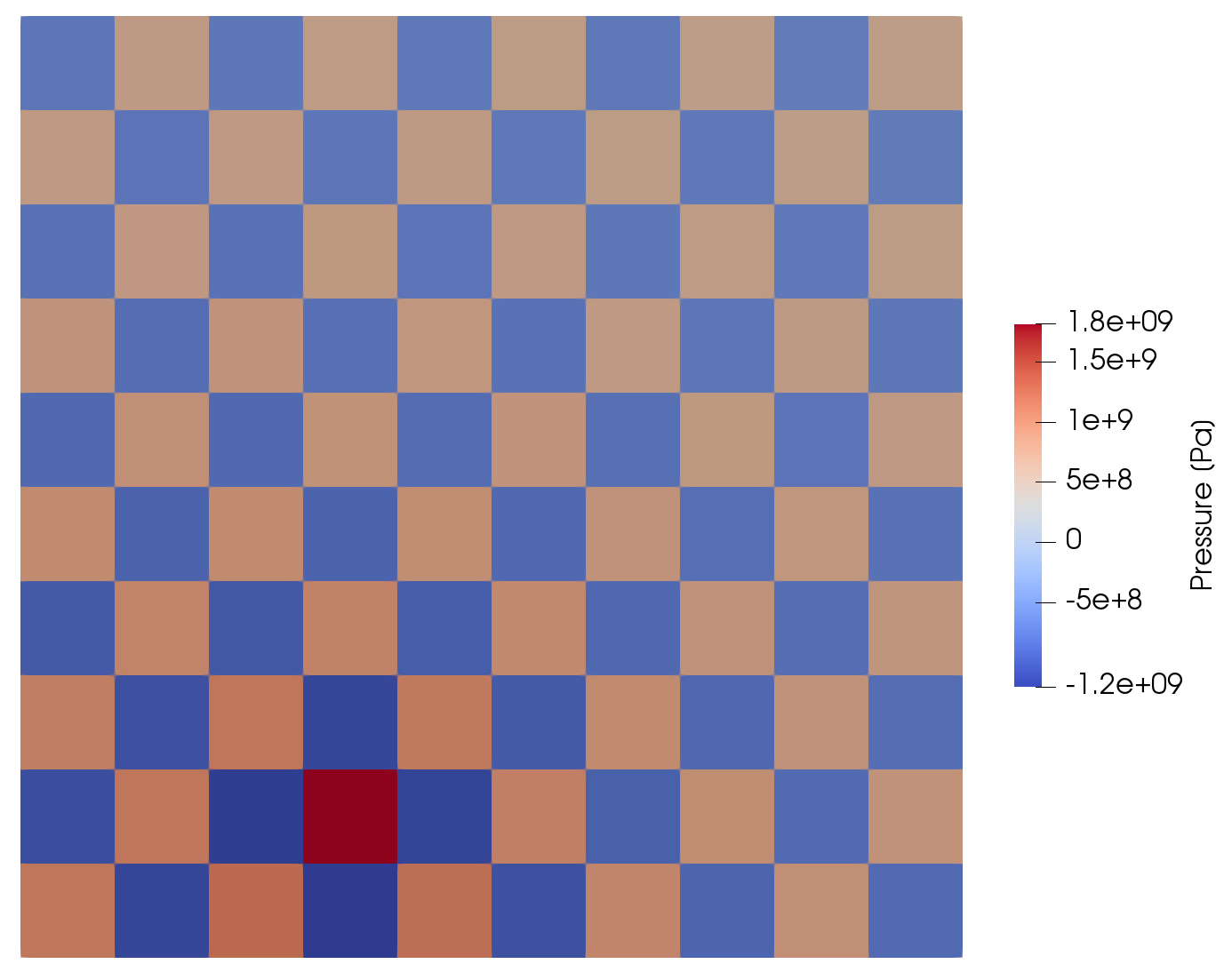}} \hfill
\subfloat[Fixed-stress without iteration]{\label{sfig:undrained_FS_noiter}\includegraphics[width=0.5\textwidth]{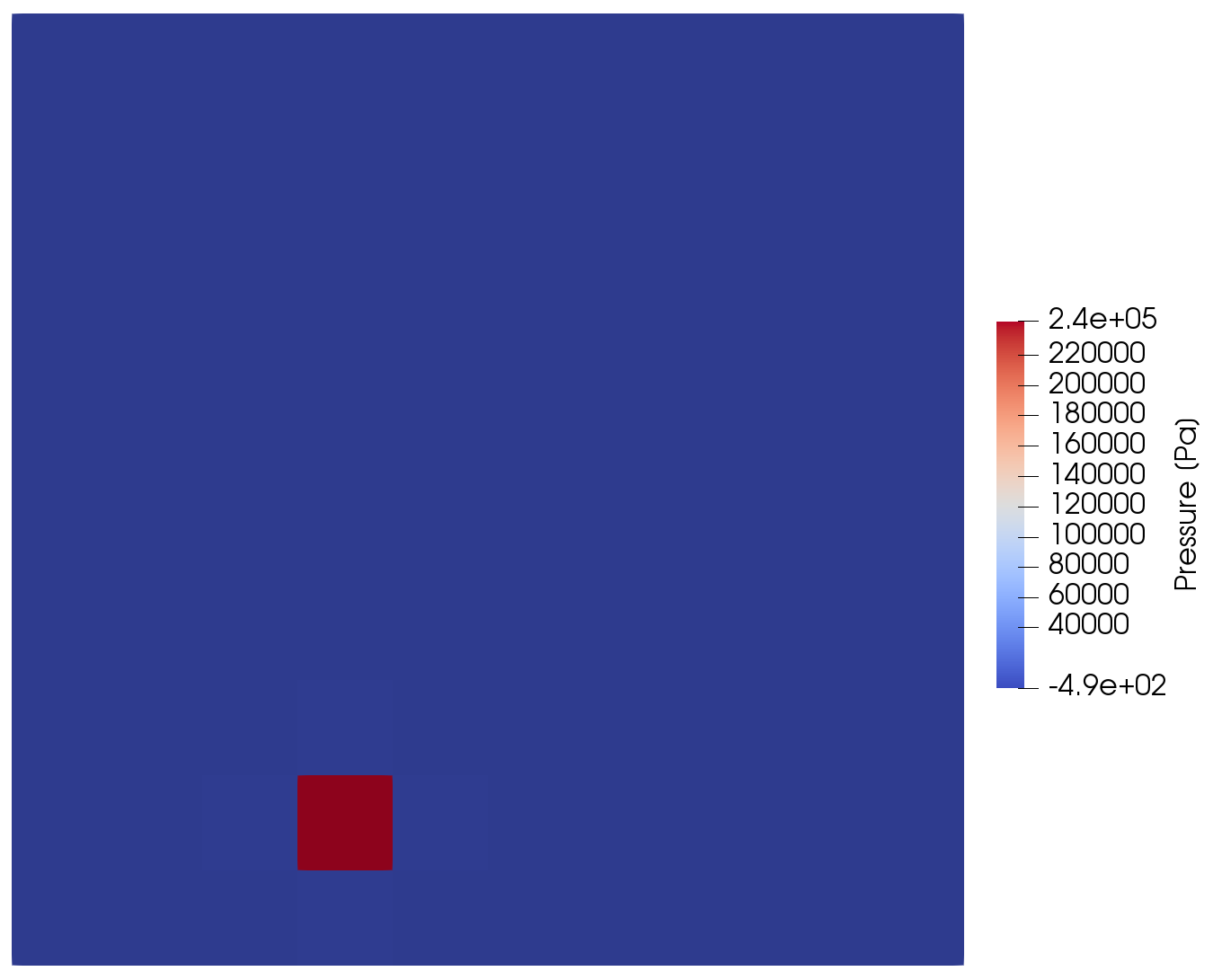}} \\
\caption{Undrained Barry-Mercer pressure solution after 10 time steps}
\label{fig:undrained}
\end{figure}

Figure \ref{sfig:undrained_fim} shows the resulting pressure field using the fully implicit method, for reference. Clearly the solution is polluted by the checkerboard mode in this undrained scenario, as expected. Figure \ref{sfig:undrained_FS_noiter}, meanwhile, shows the solution obtained with the fixed-stress method without sequential iterations within time steps. Immediately, this solution seems to be pressure stable as there is no obvious checkerboard mode, and this matches the conclusions found in \cite{yoon2018spatial}. 

\begin{figure}
\centering
\subfloat[10 iterations per time step]{\label{sfig:undrained_FS_10iter}\includegraphics[width=0.5\textwidth]{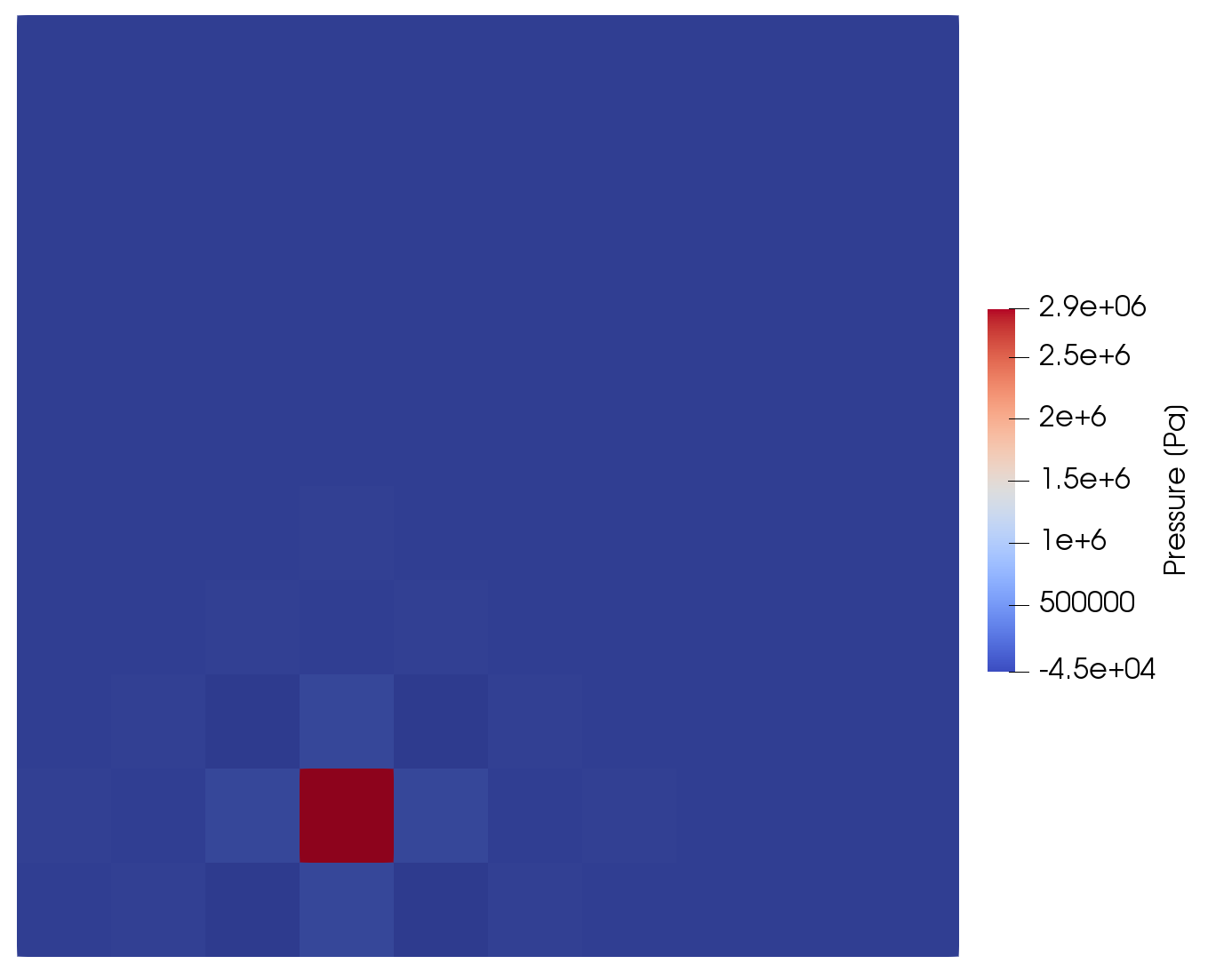}} \hfill
\subfloat[50 iterations per time step]{\label{sfig:undrained_FS_50iter}\includegraphics[width=0.5\textwidth]{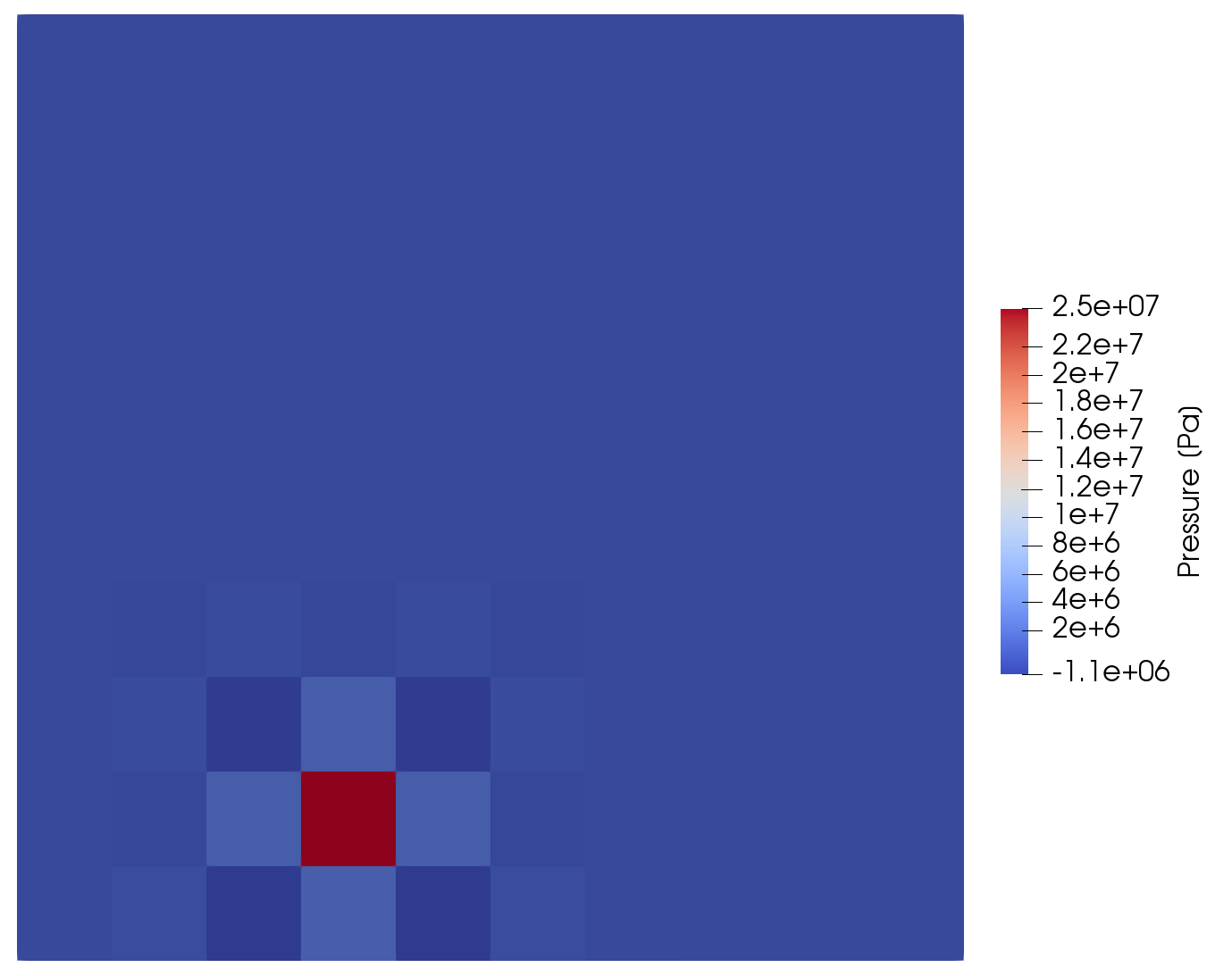}} \\
\subfloat[100 iterations per time step]{\label{sfig:undrained_FS_100iter}\includegraphics[width=0.5\textwidth]{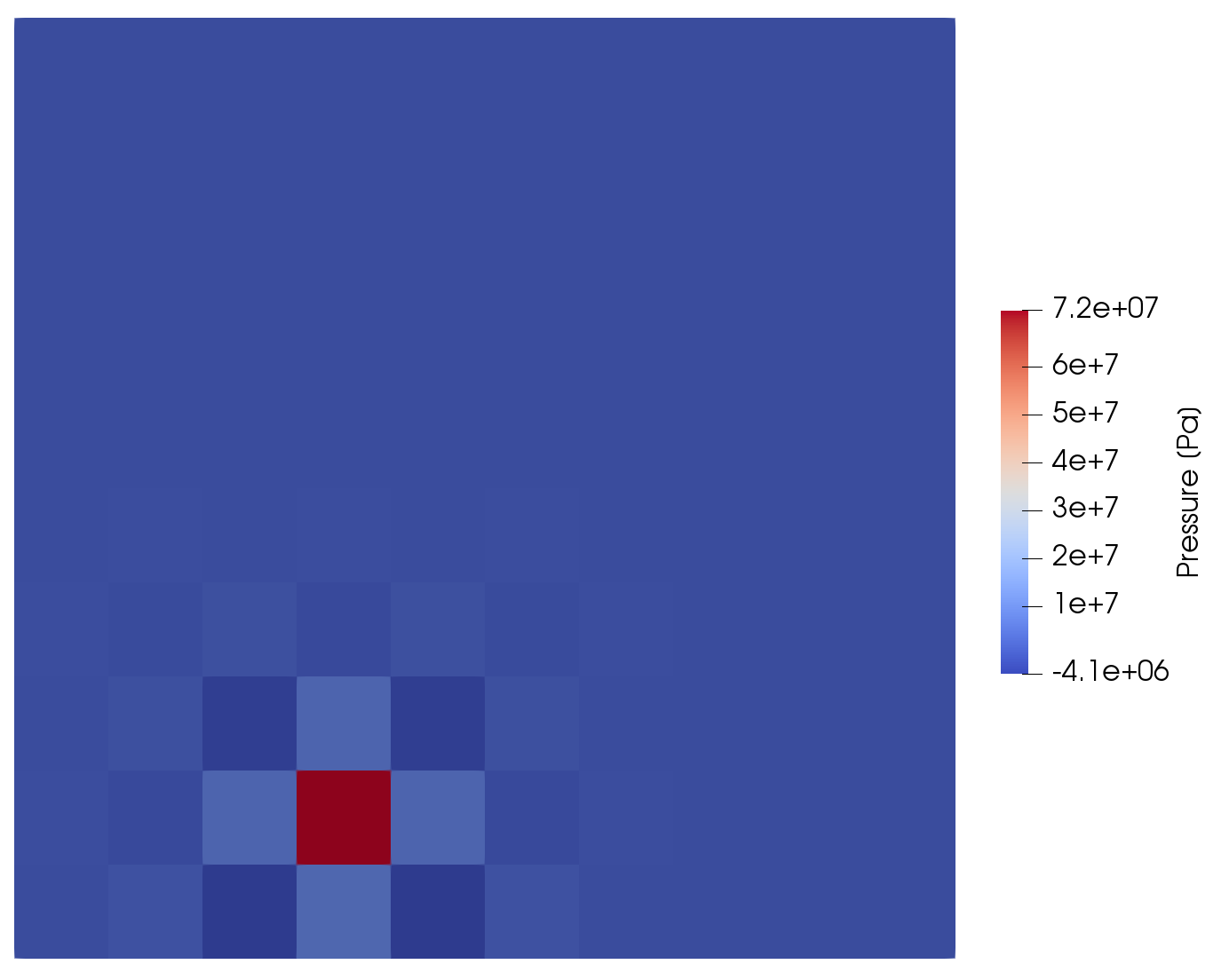}} \hfill
\subfloat[500 iterations per time step]{\label{sfig:undrained_FS_500iter}\includegraphics[width=0.5\textwidth]{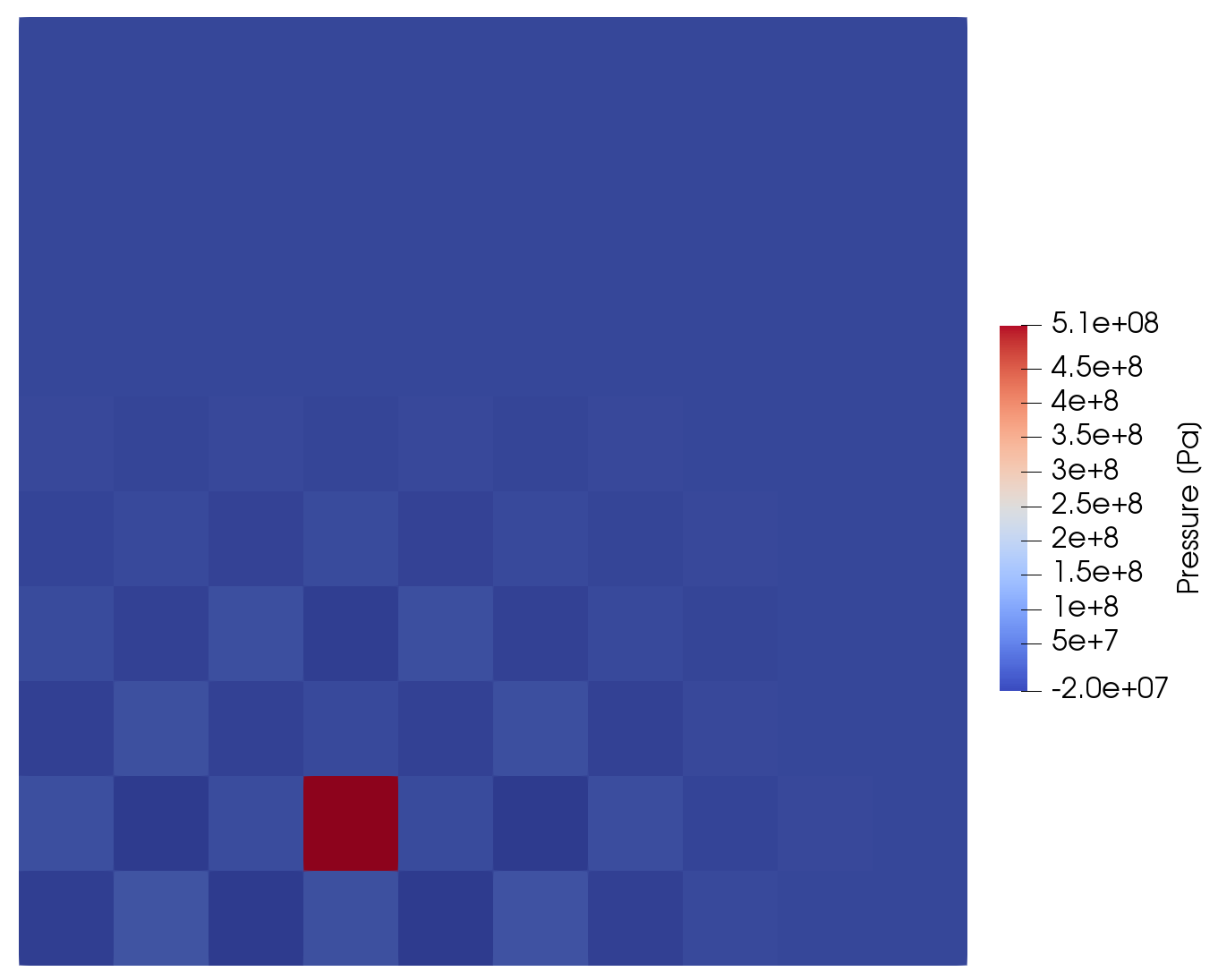}} \\
\caption{Undrained Barry-Mercer pressure solution after 10 time steps with iterative fixed-stress scheme}
\label{fig:undrained_iter}
\end{figure}

Now we consider the behavior of the fixed-stress method when sequential iterations are taken within time steps. Figure \ref{fig:undrained_iter} shows the pressure fields after 10 time steps when increasing numbers of iterations are taken in each time step. In this case the pressure solution is converging to the fully implicit solution as the number of iterations is increased. This confirms that the fixed-stress method certainly does not contain any pressure stabilizing mechanism when iterated to convergence. Moreover, we see that the convergence is very slow, which is consistent with the results in \cite{storvik2018optimization}, and makes sense when the fixed-stress method is interpreted as a Richardson iteration \cite{castelletto2015accuracy}. In particular, in the case of undrained conditions and a non-inf-sup stable discretization, the spectral radius of the iteration matrix cannot be bounded away from one and depends on the eigenvalues associated with any spurious pressure modes in the fully implicit system.

In the rest of this work we focus on the stabilization of the sequentially iterated scheme using the same stabilization schemes as in the fully implicit setting. We also apply the scheme to compositional multiphase flows and realistic examples using, for example, realistic compressibilities of the CO$_2$-rich and brine phases instead of a fully incompressible fluid, and with nearly undrained conditions only appearing in the burden regions. These will demonstrate that the effects discussed above appear not only in this academic test case, but also in plausible simulations of CO$_2$ injection and storage. 

\begin{remark}

    While the focus in this work is on clarifying that full rank matrices are not sufficient to generate a pressure stable splitting scheme, we also note that more study is warranted of the pressure stability properties of the non-iterated fixed-stress scheme, given the possibility of the stabilizing term depending on the problem setup. Even in the simple test considered here, we believe that interesting results are generated by considering the evolution of the solution over time, while also zooming in the colormap. These results are shown in Figure \ref{fig:undrained_FS_dts}. After one time step the pressure solution does seem to be monotone, and we note that many of the results in \cite{yoon2018spatial} only consider the stability after one time step. However, we see that as the solution advances in time, the spurious pressure mode begins to appear and increase in strength. This would suggest that the non-iterative fixed-stress method may not be able to provide sufficient pressure stability in the general setting where an undrained condition persists throughout the simulation, which may include many time steps. We note that some of the examples in \cite{yoon2018spatial} involving longer time integration became drained as time advanced, meaning the spurious pressure errors may have decayed by the end of the simulation.

    \begin{figure}
    \centering
    \subfloat[Pressure after 1 time step]{\label{sfig:undrained_FS_1dt}\includegraphics[width=0.5\textwidth]{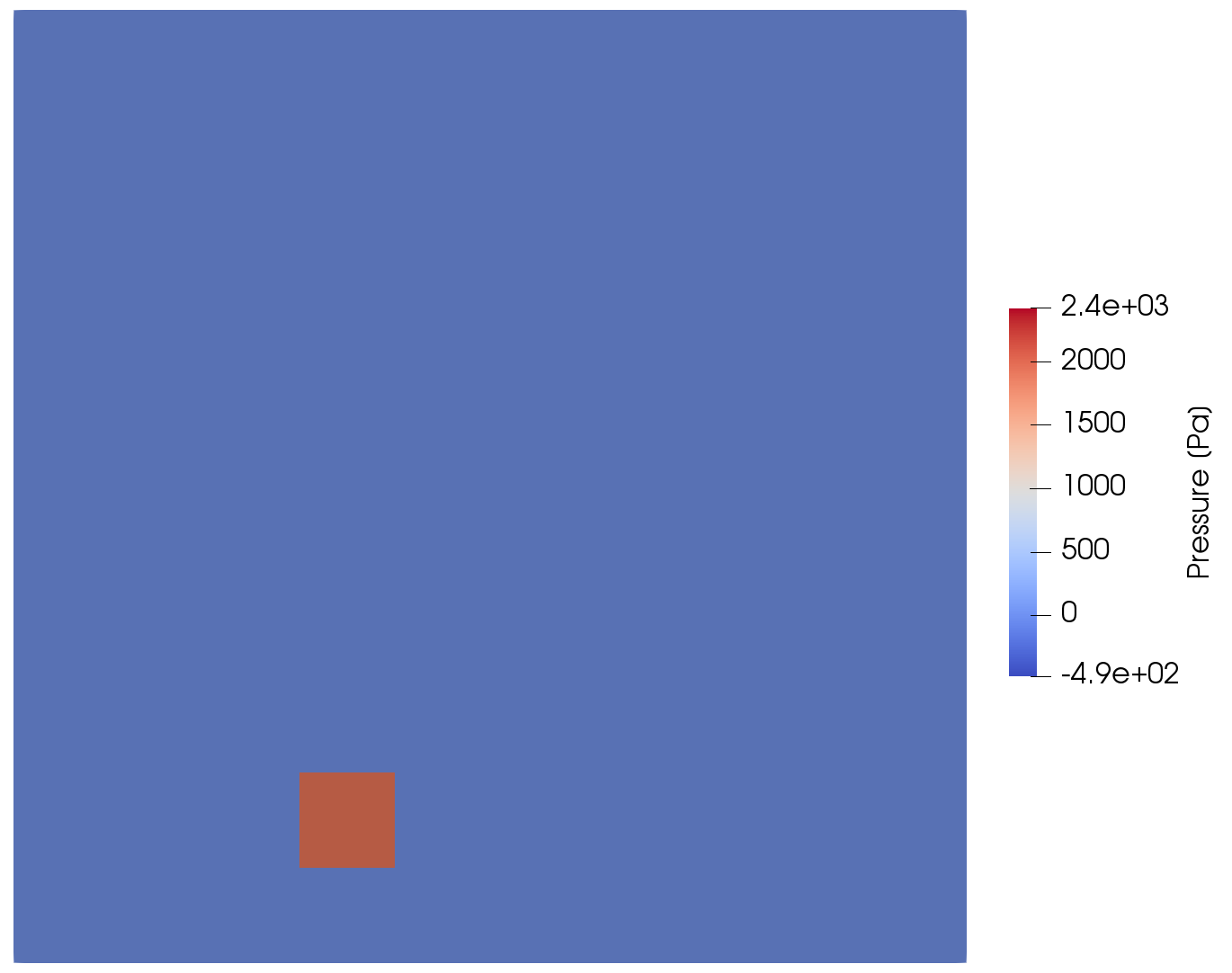}} \hfill
    \subfloat[Pressure after 4 time steps]{\label{sfig:undrained_FS_4dt}\includegraphics[width=0.5\textwidth]{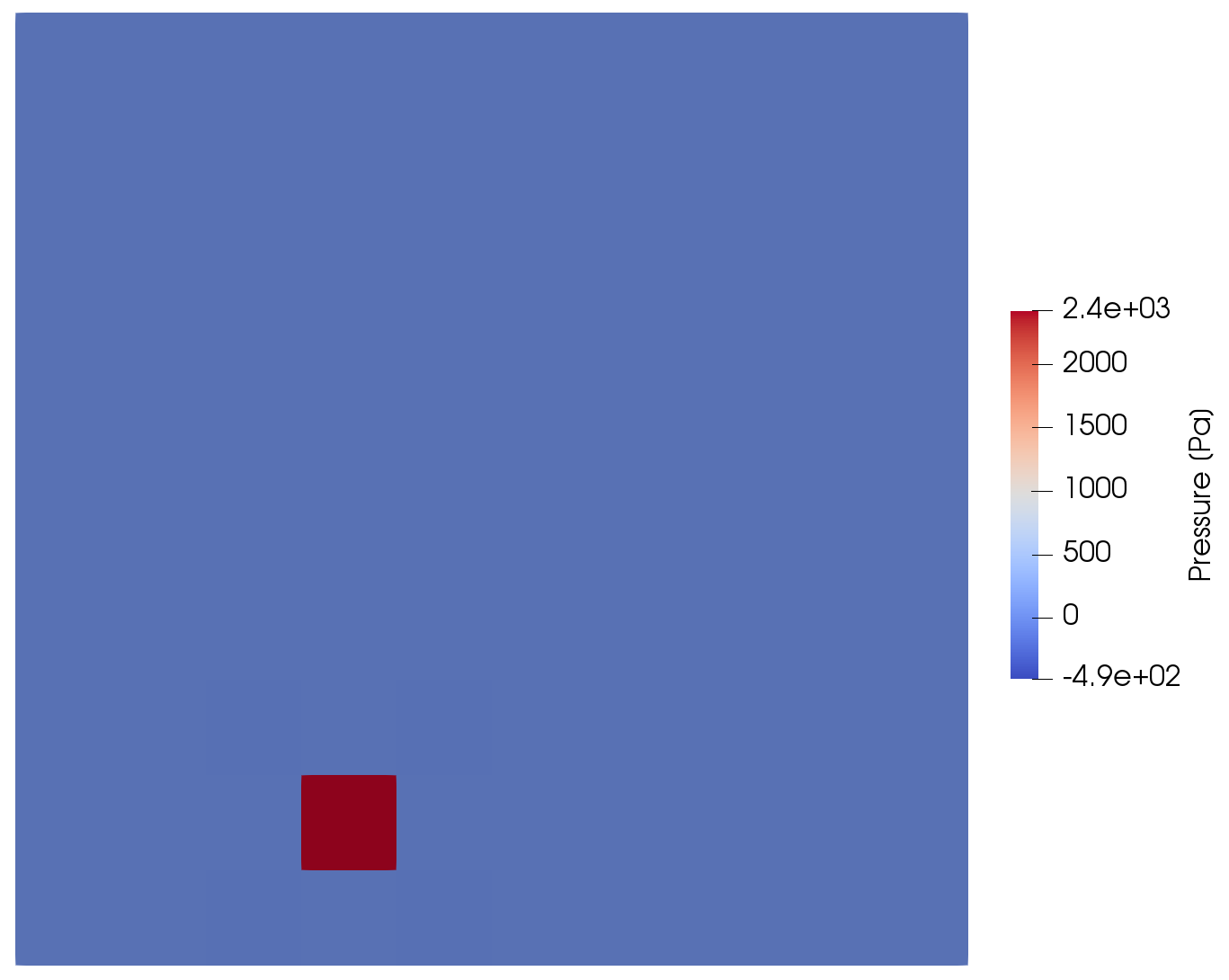}} \\
    \subfloat[Pressure after 7 time steps]{\label{sfig:undrained_FS_7dt}\includegraphics[width=0.5\textwidth]{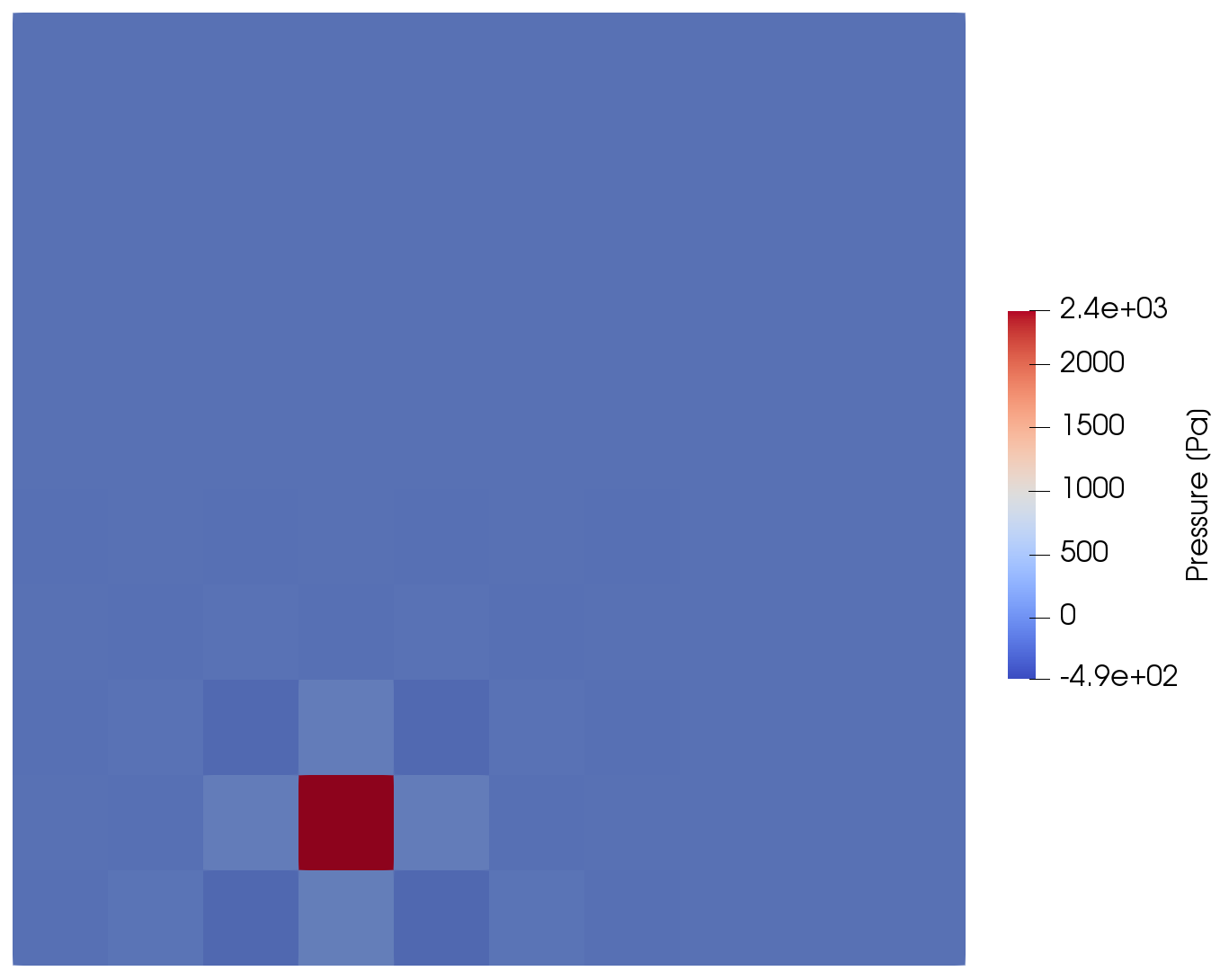}} \hfill
    \subfloat[Pressure after 10 time steps]{\label{sfig:undrained_FS_10dt}\includegraphics[width=0.5\textwidth]{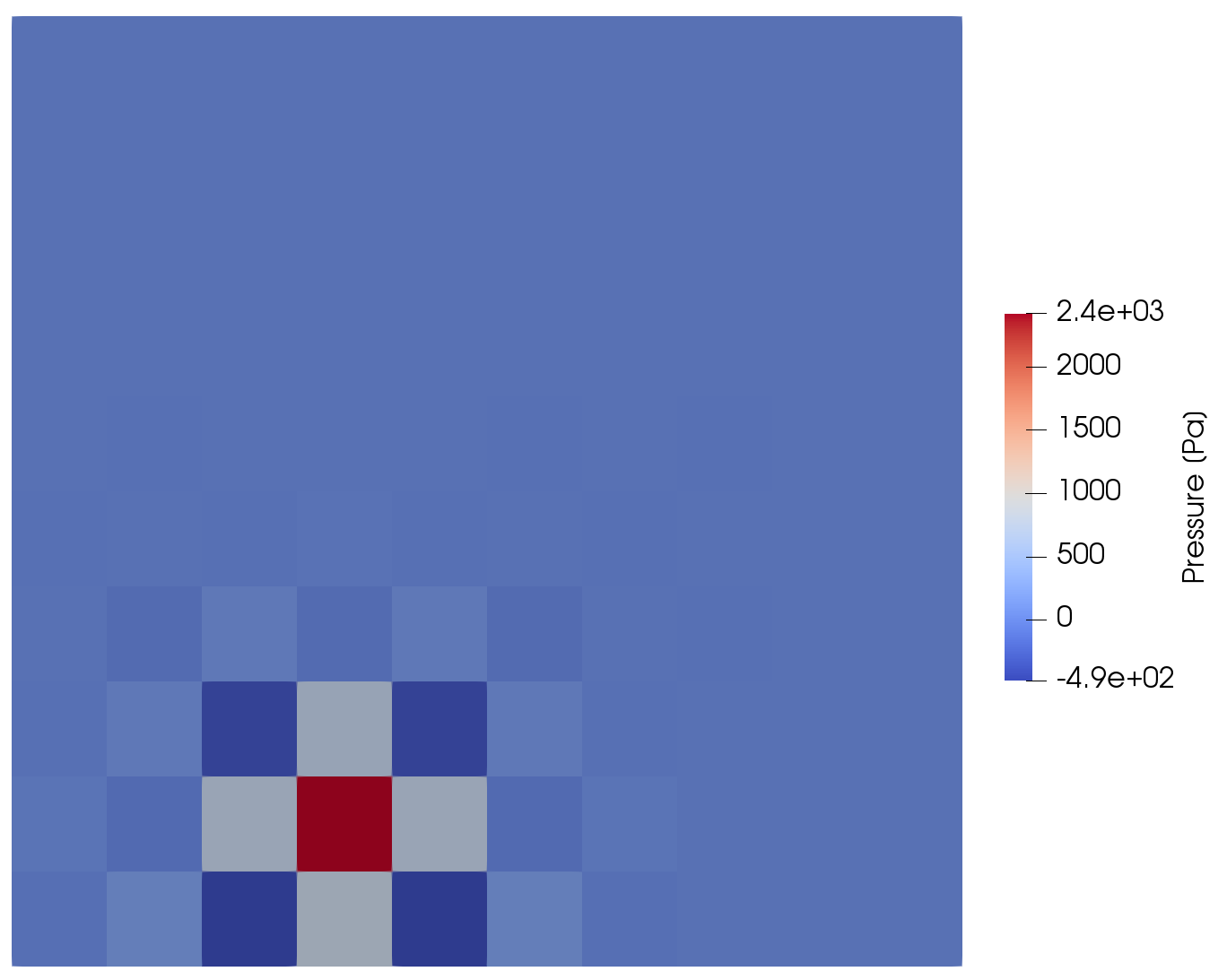}} \\
    \caption{Undrained Barry-Mercer pressure solution at various times with non-iterative fixed-stress scheme}
    \label{fig:undrained_FS_dts}
    \end{figure}
\end{remark}

\begin{remark}
    Some analysis can be done on the non-iterated fixed-stress splitting error term: $-B(u^{n+1} - 2u^n + u^{n-1}) + R(p^{n+1} - 2p^n + p^{n-1})$. The fixed stress assumption of $\sigma^{n+1}_{v} - \sigma^{n}_{v} = \sigma^{n}_{v} - \sigma^{n-1}_{v}$ is a discretized statement of the second time derivative of the total stress being equal to zero. Differentiating Equation \eqref{eq:single_phase_mass} in time, discretizing using the same scheme, applying the fixed-stress assumption, and ignoring forcing terms for simplicity yields 

    \begin{equation}
        p^{n+1} - 2p^n + p^{n-1} + \nabla \cdot \mathbf{\dot{v}}^{n+1} = 0.
    \end{equation}

    \noindent In general, we could use Darcy's Law (Equation \eqref{eq:darcy}) to convert $\mathbf{\dot{v}}$ into the gradient of the pressure change, which would then imply that $p^{n+1} - 2p^n + p^{n-1}$ is related to the Laplacian of the pressure change. This resembles the error terms in both the incremental and non-incremental pressure projection schemes, and would seem promising as a possible source of pressure stabilization. However, in the limit of low permeabilities, $\mathbf{\dot{v}}^{n+1} = 0$, and thus $p^{n+1} - 2p^n + p^{n-1}$ is also zero. This could explain the appearance of oscillations in Figure \ref{fig:undrained_FS_dts}, though the possible effect of the remaining terms in the splitting error is less immediately apparent in general.
\end{remark}

\section{Pressure Jump Stabilization in Fixed-Stress Splitting}

Given the need for pressure stabilization in the case of the iterated fixed-stress scheme (and perhaps for the non-iterated scheme as well), we propose the inclusion of pressure jump stabilization within the flow solver step. This stabilization technique was originally introduced in the context of incompressible flow in \cite{hughes_1987}, and since used successfully in single-phase, immiscible two-phase, and compositional multiphase poromechanics problems \cite{berger2015stabilized, camargo2021macroelement, aronson2023pressure}. Pressure jump stabilization essentially constructs an artificial stabilizing flux which is designed to oppose the element-to-element oscillations characteristic of a lack of inf-sup stability. The stabilizing flux which is added to the standard, TPFA numerical flux for each phase is of the form \cite{aronson2023pressure}

\begin{equation}
    F^{stab, f}_c = \tau V \sum_{\ell=1}^{n_p} \left( x_{c\ell} \rho_\ell k_{r\ell}\right) ^{upw}_{n} [[p^{n+1}-p^n]].
    \label{eq:stab_flux}
\end{equation}

\noindent Here $[[\cdot]]$ denotes the jump operator across a face, and the pressure jump is multiplied by a component mass term which is lagged in time in order to simplify linearization when performing Newton iterations. We denote by $V$ the volume term which is computed in a similar manner to the geometric transmissibility term in the TPFA numerical flux (details on this computation can be found in Appendix A of \cite{white2019two}, for example), and $\tau$ is a stabilization parameter to be defined. Adding this flux to every face in the mesh results in what we will call the global stabilization method \cite{hughes_1987, berger2015stabilized}. The mesh can also be partitioned into macro-elements, and the stabilization can be applied only on the interior of these macro-elements \cite{silvester1990stabilised, silvester1994optimal, camargo2021macroelement}. This has certain theoretical advantages over the global method, including the preservation of mass conservation at the macro-element level. The solution is also less sensitive to over-stabilization by a poor (i.e. too large) choice of stabilization constant $\tau$ due to the fact that, as $\tau \rightarrow \infty$, the scheme reverts to a piecewise constant approximation on every macro-element.

Some theory has also been developed to guide the choice of stabilization constant $\tau$ for local stabilizations. In particular, it was proposed in \cite{silvester1994optimal, camargo2021macroelement} that a good choice of $\tau$ is the one which minimizes the condition number of the pressure-only system over one macro-element. Following \cite{camargo2021macroelement}, where this analysis was performed for immiscible poromechanics on a mesh composed of 3D macro-elements formed by regular hexahedra, the optimal $\tau$ which minimizes the condition number of the undrained pressure Schur complement matrix is computed analytically to be

\begin{equation}
    \tau^* = \left( \frac{9}{32(\lambda + 4G)} \right).
    \label{eq:opt_tau}
\end{equation}

\noindent Here $\lambda$ is the first Lamé parameter of the solid skeleton and $G$ is the second (also known as the shear modulus).

In \cite{aronson2023pressure} it was shown that this value should be scaled up by a constant (3 was selected in that work based on the conditioning studies with a macro-element composed of tetrahedra) in the case of tetrahedral elements, and that these values also perform well for the global stabilization scheme. Inspired by this, we choose in this work to use the global stabilization scheme in conjunction with the fixed stress split, with stabilization constant selected as 

\begin{equation}
    \tau = c \tau^*,
\end{equation}

\noindent with $c = 1$ for hexahedral meshes and $c = 3$ for tetrahedral meshes, unless otherwise specified. Using a global stabilization yields a simple method which allows us to avoid the more complex problem of agglomerating fully unstructured meshes into macro-elements (though this is certainly possible, as seen in \cite{borio2021hybrid, aronson2023pressure}). We expect that a local approach would also sufficiently stabilize the examples shown below. 

As a last note on pressure jump stabilization, we have found that this approach can produce smeared pressure interfaces at the locations in the domain where the high permeability reservoir and low permeability burden regions meet. This is not unexpected given that the stabilization operator penalizes jumps in pressure between cells, and was also seen in \cite{berger2017stabilized, aronson2023pressure}. In \cite{berger2017stabilized} it was proposed that these can be removed by more careful tuning of the stabilization strength $\tau$. Alternatively, one could also only apply the stabilization within the actual undrained region, namely the burden regions surrounding the reservoir, as in \cite{aronson2023pressure}. In this work we elect to use the latter, as it avoids more cumbersome parameter tuning processes, and in our applications it is usually known where the interface between the high and low permeability lies. 

Of course, pressure jump stabilization is far from the only possible method that could be applied here. However, one of the main advantages of sequential coupling schemes is the possibility of using existing, single-physics solvers which are already mature and well optimized, and pressure jump stabilization is minimally invasive when compared to other potential stabilization schemes. For example, using residual-based \cite{wan2003stabilized} or projection \cite{white2008stabilized} stabilization schemes would require components that are likely not readily available in a standard, finite volume flow solver. The former would need to evaluate the mechanical residual in each element during the flow solve, while the latter would require projection onto the piecewise linear, dual mesh basis in the case of a piecewise constant finite volume solver \cite{dohrmann2004stabilized}. 

Pressure jump stabilization, on the other hand, is easily added to a pure finite volume flow solver, as a loop over faces to compute the physical flux already exists. The stabilization flux also takes the same form as the TPFA flux, being composed of an upwind mass term and a potential difference. One simply needs to augment this physical flux based on the flow variables, and the only mechanics information that is needed are the material parameters used to compute the stabilization constant.

\section{Results}

With the stabilized scheme defined, we now consider examples involving compositional multiphase poromechanics. We start with a smaller example on a hexahedral grid, then consider the High Island 24L test case (HI24L) from the Gulf of Mexico, which has been studied in the context of carbon sequestration in, for example, \cite{t2022deformation}. The examples are simulated using GEOS \cite{GEOS}, an open-source, high performance simulator being developed for the simulation of the multiphysics processes associated with coupled subsurface flows and geomechanics. The use of GEOS also allows us to make comparisons to fully implicit solutions, when appropriate. For each example we will show that the number of sequential iterations required to reach convergence within a time step increases as undrained conditions are approached, and pressure jump stabilization accelerates this sequential convergence. We say that the sequential iterations have converged when the coupled residual norm is a factor of $1.0 \times 10^{8}$ smaller than its initial value, and we solve both the solid and fluid subproblems with a Newton-Raphson method with the same tolerance. The linear problems produced are then solved with GMRES with algebraic multigrid and multigrid reduction preconditioners for the solid and fluid problems, respectively \cite{bui2021multigrid}. By appealing to the sequential convergence rates we can also show that the previously defined optimal stabilization parameters are effective at improving efficiency without the risk of over-stabilization. We will also show that the unstabilized fixed-stress results exhibit the same spurious pressure oscillations as the unstabilized fully implicit scheme, and once again pressure jump stabilization can smooth these without adversely affecting the rest of the solution quality.

\subsection{Staircase}

\begin{figure}
\centering
\includegraphics[width=0.5\textwidth]{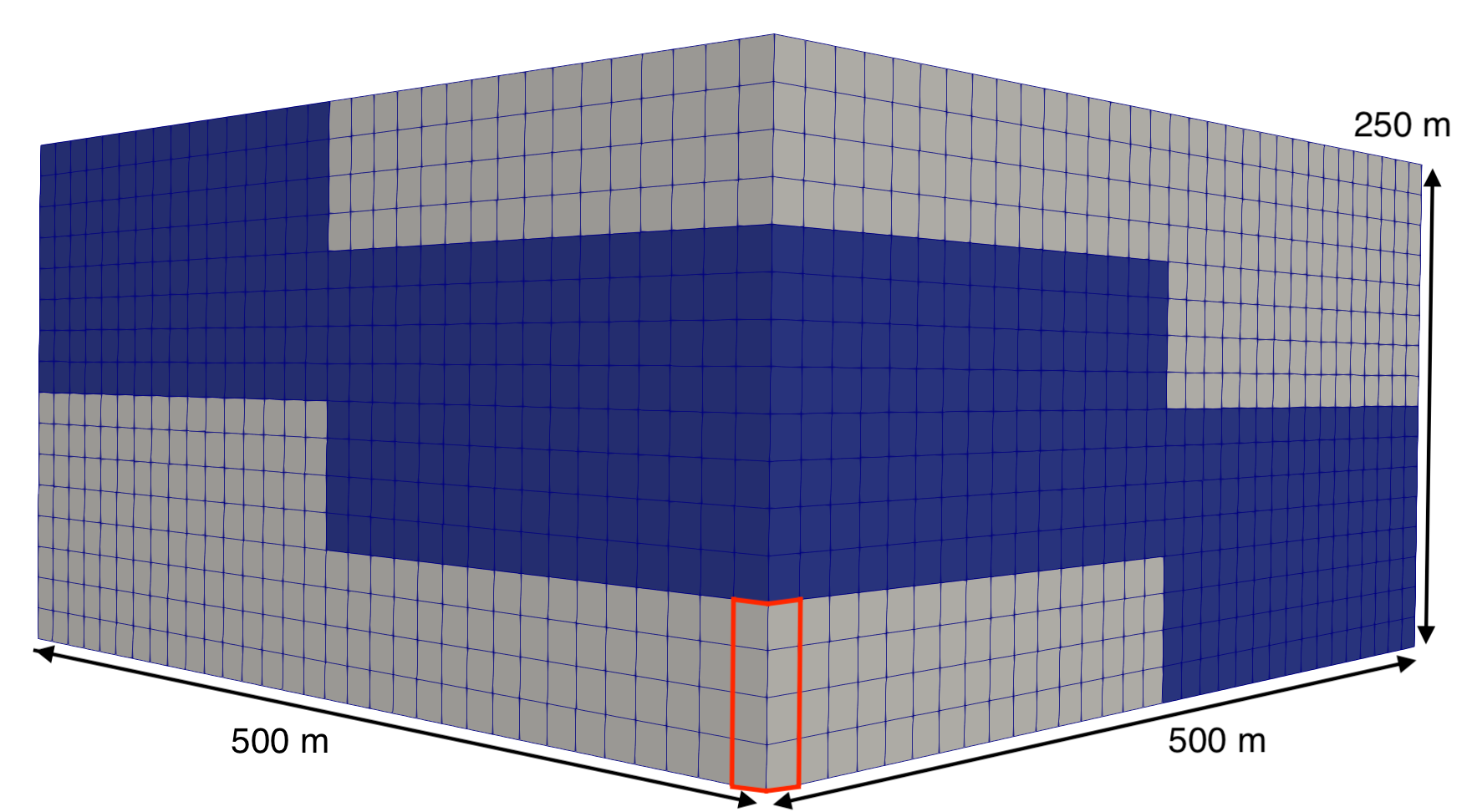}
\caption{Staircase mesh. Gray regions represent the channel while blue represents the barrier region. Outlined in red are cells where CO$_2$ is injected. }
\label{fig:stair_setup}
\end{figure}

We start by considering the staircase example studied in the context of immiscible two-phase flows in \cite{white2019two, camargo2021macroelement}. The simulation mesh is shown in Figure \ref{fig:stair_setup} and consists of a high permeability channel spiraling around a low permeability burden region. The solid parameters of the channel and burden regions are given in Table \ref{tab:params_solid}, and we note that we use burden permeability values of $9.8\times 10^{-14}$ m$^2$, $9.8\times 10^{-17}$ m$^2$, and $9.8\times 10^{-20}$ m$^2$ in some of our tests to show the change in performance of the fixed-stress scheme as the region changes from drained to undrained. The fluid parameters are as in Table \ref{tab:params_fluid}. CO$_2$ is injected at a rate of 10 kg/s into each of the four cells outlined in red in Figure \ref{fig:stair_setup}. 

\begin{table}[ht]
\centering
\caption{Solid material parameters used in synthetic aquifer example}
\begin{tabular}{ccc}
\hline
{ \textbf{Parameter}} & { \textbf{Channel Value}} & { \textbf{Burden Value}}  \\ \hline
Porosity                 & 0.2                            & 0.05                         \\
Permeability             & $9.8\times 10^{-13}$ m$^2$                        & $9.8\times 10^{-14}$ - $9.8\times 10^{-20}$ m$^2$                    \\
Skeleton bulk modulus       & 5 GPa                          & 5 GPa                        \\
Skeleton Poisson ratio       & 0.25                          & 0.25                        \\
Skeleton Biot coefficient   & 1                              & 1                            \\
Skeleton density            & 2650 kg/m$^3$   & 2650 kg/m$^3$                             
\end{tabular}
\label{tab:params_solid}
\end{table}

\begin{table}[ht]
\centering
\caption{Fluid material parameters used in synthetic aquifer example}
\begin{tabular}{cc}
\hline
{ \textbf{Parameter}} & { \textbf{Value}}  \\ \hline
CO$_2$ solubility           & \cite{duan2003improved}        \\
CO$_2$ phase density        & \cite{span1996new}          \\
CO$_2$ phase viscosity      & \cite{fenghour1998viscosity}      \\
Brine phase density      & \cite{phillips1981technical}        \\
Brine phase viscosity    & \cite{phillips1981technical}                           
\end{tabular}
\label{tab:params_fluid}
\end{table}

\begin{figure}
\centering
\subfloat[$\delta t$ = 1 hour, no stabilization]{\label{sfig:stair_1hour_nostab}\includegraphics[width=.5\textwidth]{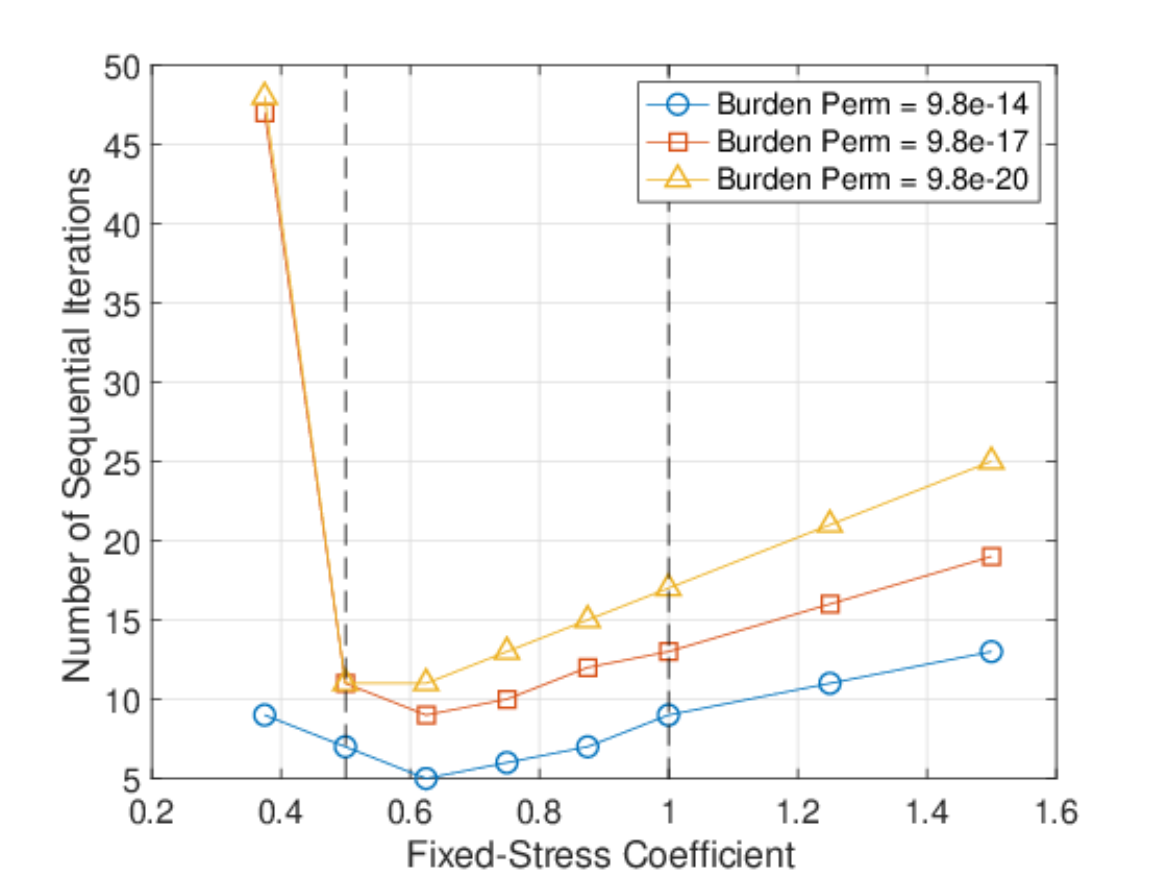}}\hfill
\subfloat[$\delta t$ = 1 hour, burden only stabilization]{\label{sfig:stair_1hour_stab}\includegraphics[width=.5\textwidth]{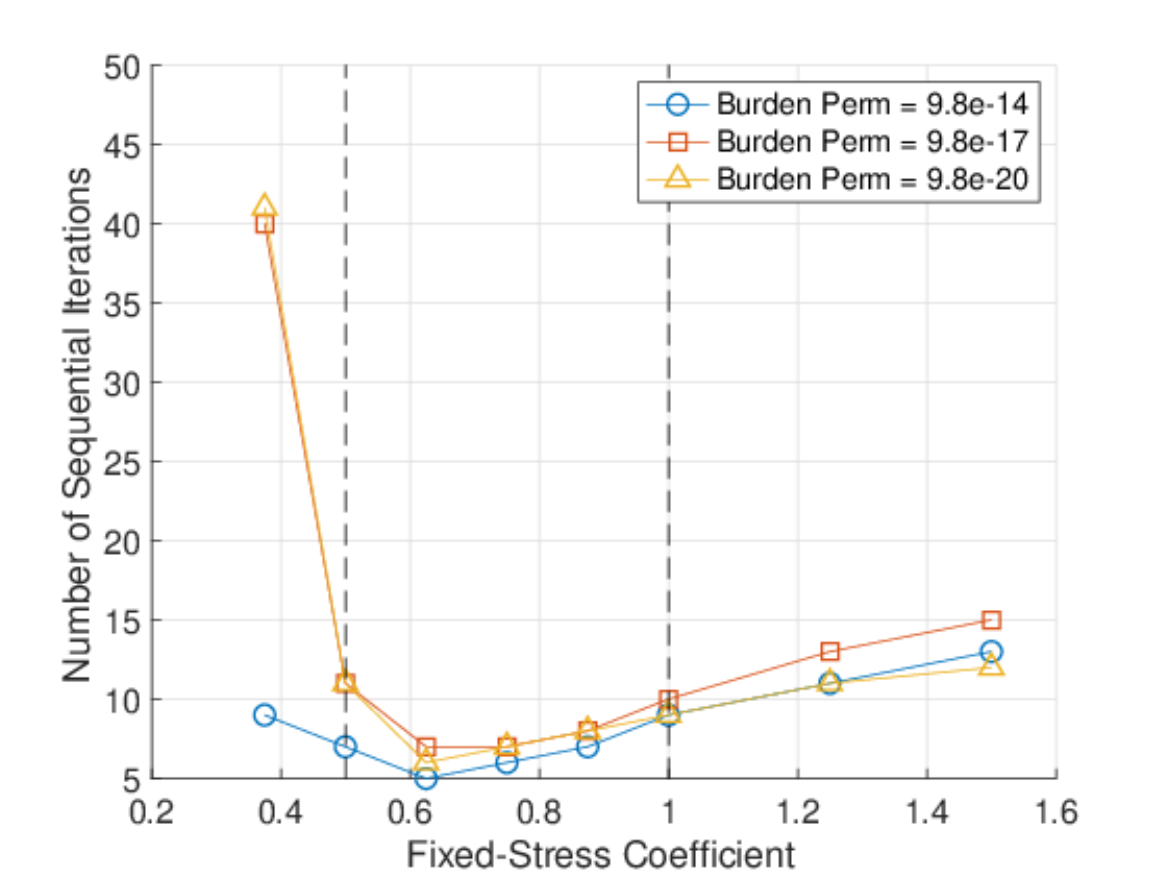}}\\
\subfloat[$\delta t$ = 1 day, no stabilization]{\label{sfig:stair_1day_nostab}\includegraphics[width=.5\textwidth]{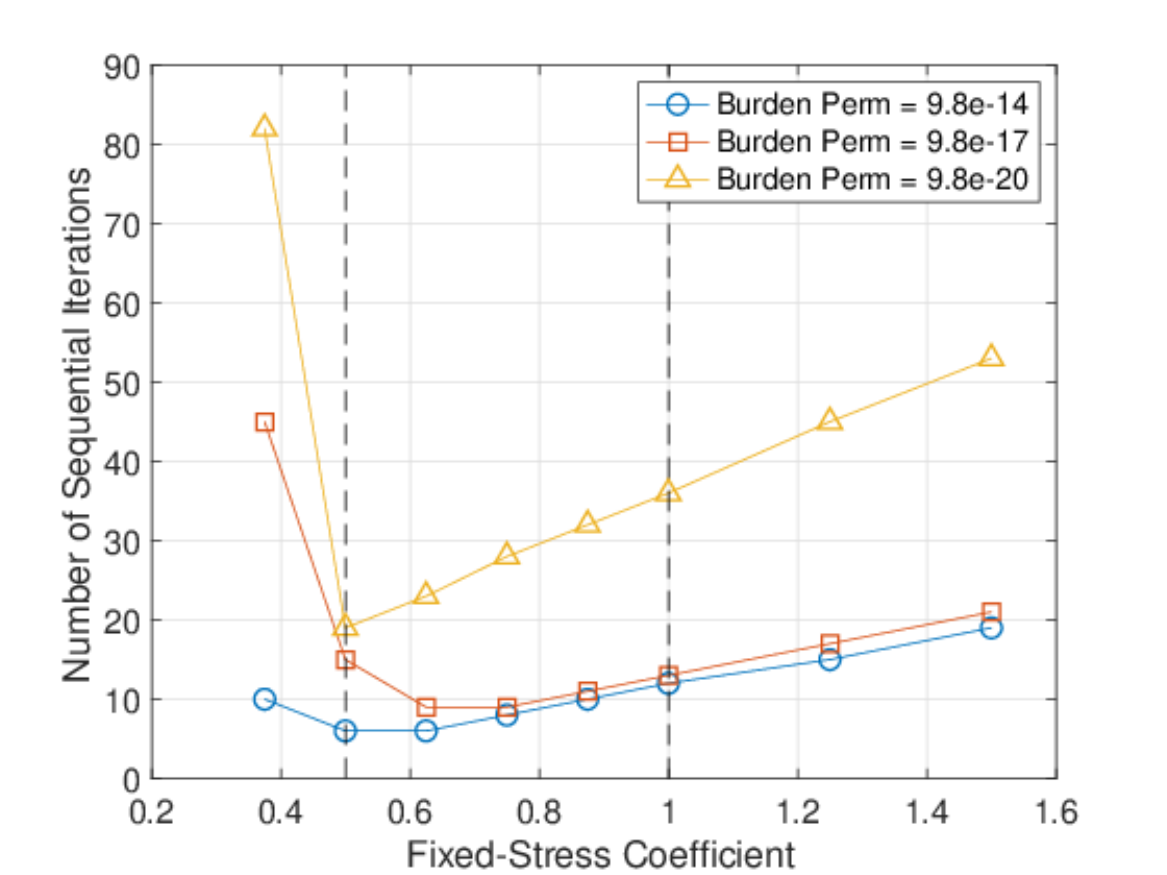}}\hfill
\subfloat[$\delta t$ = 1 day, burden only stabilization]{\label{sfig:stair_1day_stab}\includegraphics[width=.5\textwidth]{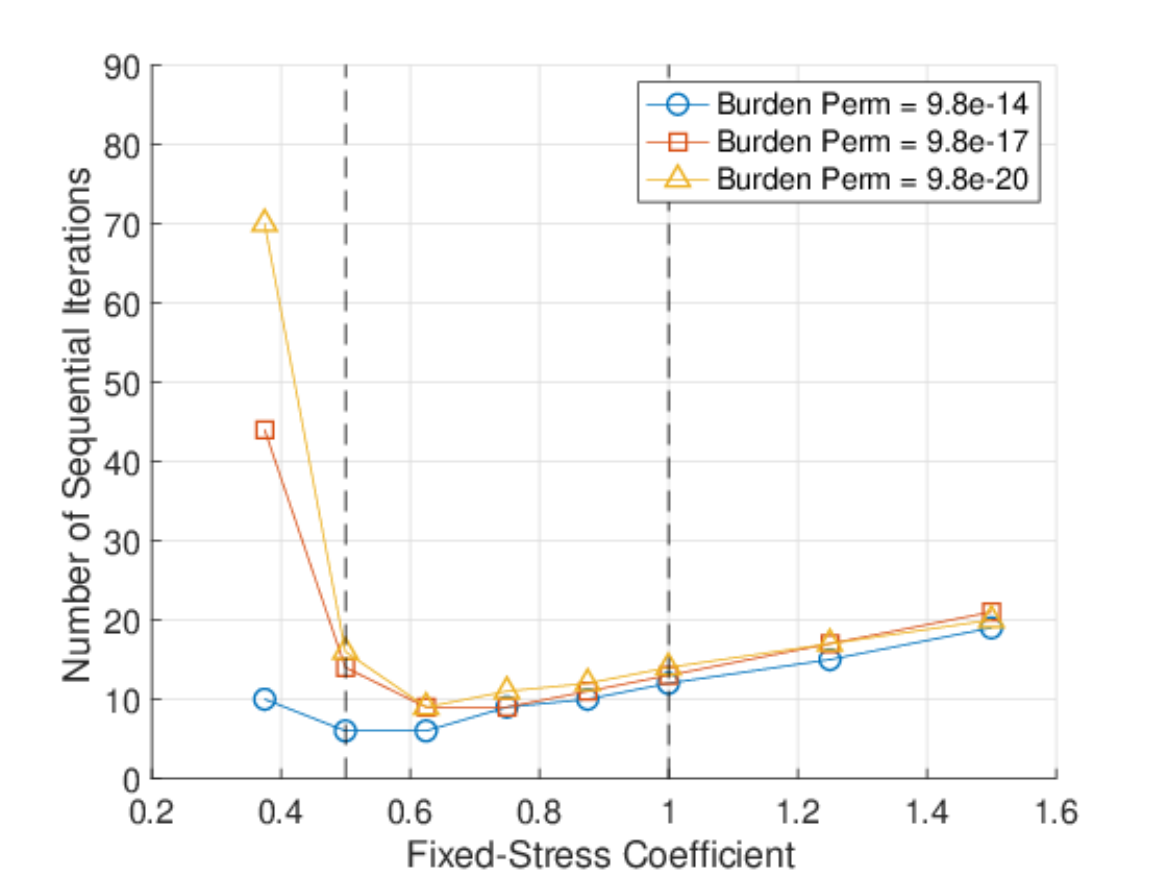}}\\
\subfloat[$\delta t$ = 1 month, no stabilization]{\label{sfig:stair_1month_nostab}\includegraphics[width=.5\textwidth]{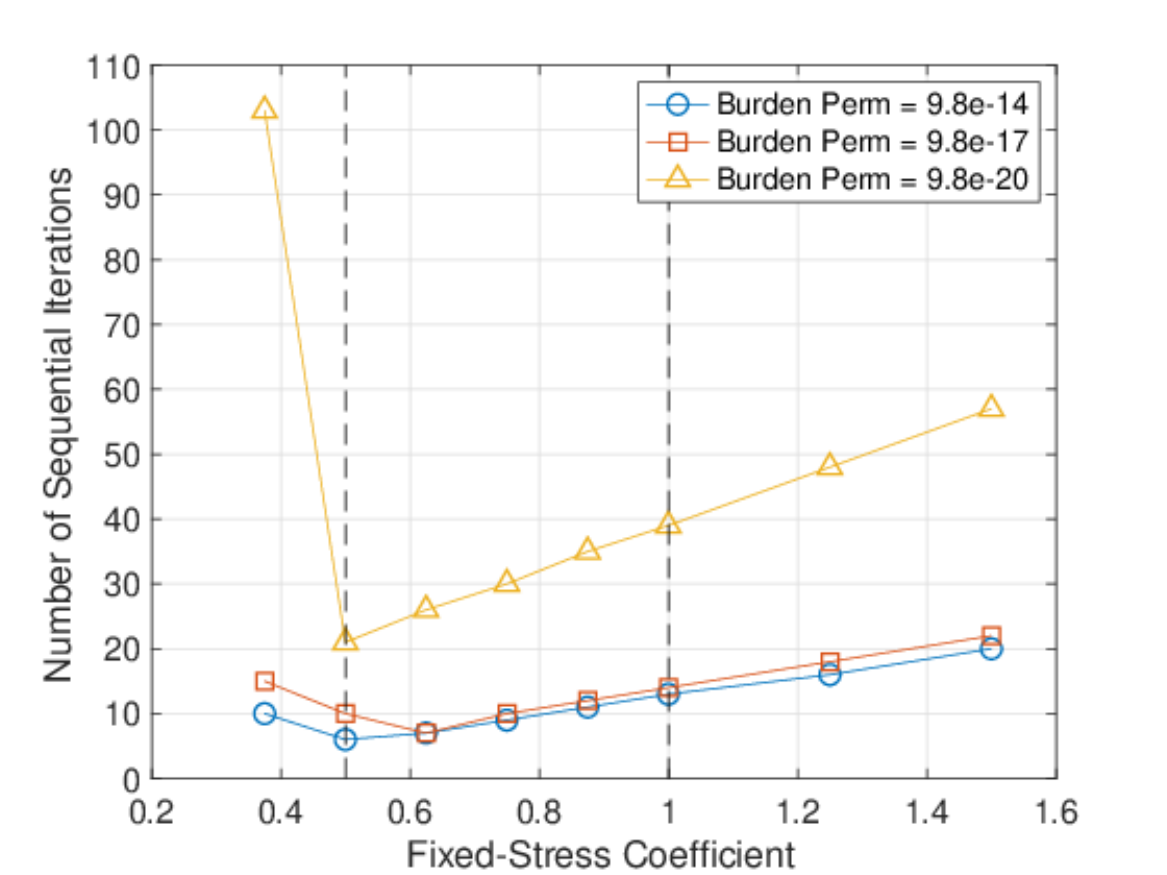}}\hfill
\subfloat[$\delta t$ = 1 month, burden only stabilization]{\label{sfig:stair_1month_stab}\includegraphics[width=.5\textwidth]{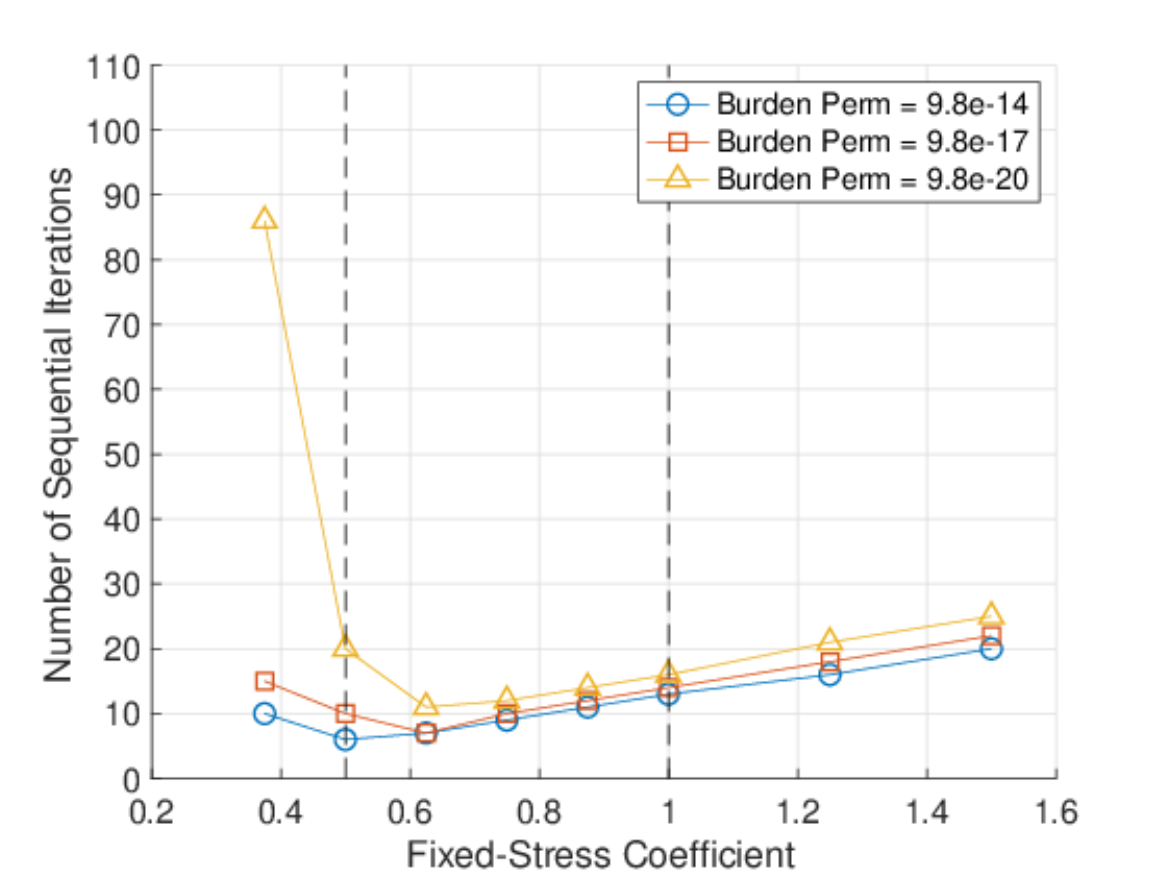}}\\
\caption{Staircase: Number of fixed-stress iterations with varying fixed-stress coefficients $\alpha$}
\label{fig:stair_iter_FScoeff}
\end{figure}

To begin, we consider the convergence behavior of the iterative scheme on the first time step of the simulation. We record the number of sequential iterations taken to reach a convergence tolerance of $1\times 10^{-8}$ for a variety of choices of time step sizes, barrier permeability values, and coefficients $\alpha$ (as defined in Section 7.2), and these results are summarized in Figure \ref{fig:stair_iter_FScoeff}. Each row of plots corresponds to a different time step size, while the left and right columns include results without and with stabilization, respectively. The horizontal axis of each figure represents the coefficient $\alpha$, and the colored lines correspond to the iterations required for convergence at different permeability values. The two dashed vertical lines represent $\alpha = 0.5$, which is the smallest value for which convergence is proven, and $\alpha = 1.0$, which corresponds to the conservation of the physical total stress in the definition of the fixed-stress scheme. It is clear from all of these results that, without stabilization, the convergence rate degrades as undrained conditions are approached in the barrier region. The inclusion of stabilization essentially makes the convergence rate agnostic to the barrier region permeability in the relevant range of $\alpha \geq 0.5$, accelerating the convergence in undrained problems while having little effect when the barrier is drained. These results are similar to those generated when comparing an inherently stable discretization with an unstable one \cite{storvik2018optimization}, and we note that when stabilization is included we also recover the behavior expected in which there is some $\alpha \in [0.5, 1]$ where the convergence is fastest \cite{storvik2019optimization}, even in the undrained setting. Finally, we note that these stabilized results are generated using stabilization in the barrier region only, and the results using stabilization throughout the domain are very similar. 

\begin{figure}
\centering
\includegraphics[width=0.5\textwidth]{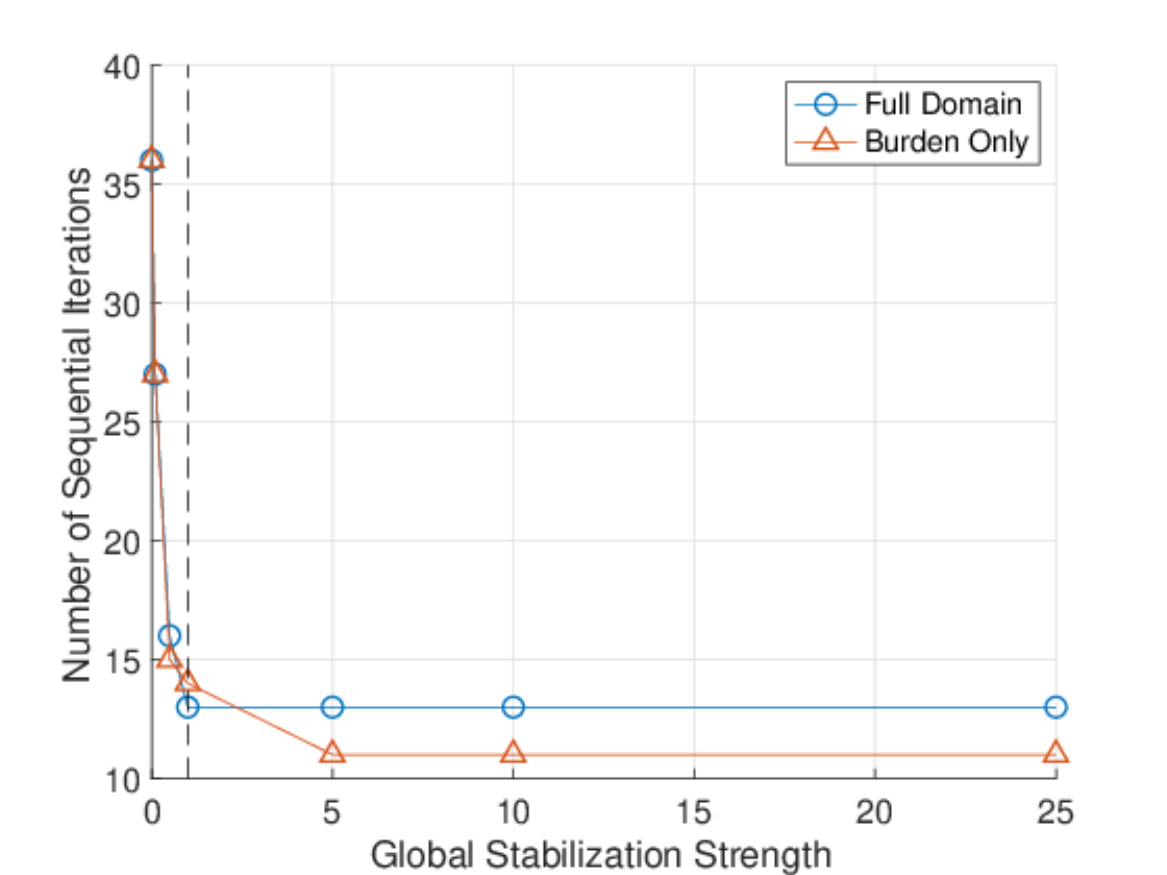}
\caption{Staircase: Number of fixed-stress iterations with varying pressure stabilization coefficient $c$ and undrained burden region}
\label{fig:stab_strength_stair}
\end{figure}

We can also study the performance of the iterative scheme as a function of the stabilization strength defined by the constant $c$ in the previous section. Figure \ref{fig:stab_strength_stair} shows the number of iterations required for convergence at the initial time step of one day with a barrier permeability of $9.8\times 10^{-20}$ m$^2$ and $\alpha = 1$. Here we directly compare the results obtained with stabilization everywhere in the domain and stabilization in the barrier region only. Finally, the vertical dashed line represents $c = 1$, which is the optimal value based on condition number estimates for hexahedral meshes. For values of $c$ less than one the number of required iterations increases dramatically, and for both the full domain and burden only stabilization the previously derived optimum of $c = 1$ is close to the knee of the plot. This is desirable, as this likely means we are adding enough stabilization to smooth spurious oscillations without overdamping the solution and lends credence to our previous analyses.  

\begin{figure}
\centering
\subfloat[No stabilization]{\label{sfig:stair_nostab_30days}\includegraphics[width=.5\textwidth]{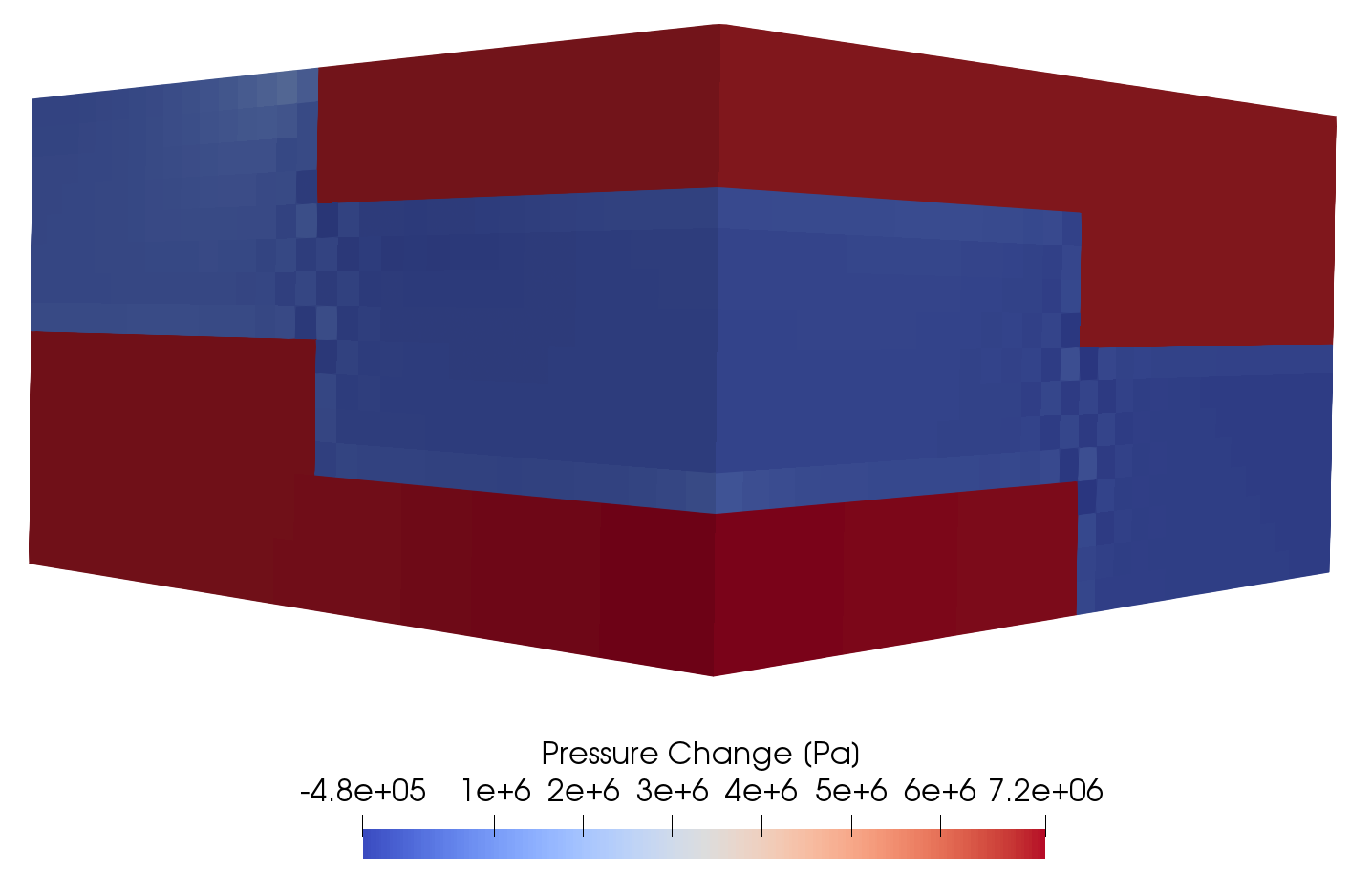}}\hfill
\subfloat[No stabilization (fully implicit)]{\label{sfig:stair_nostab_30days_fim}\includegraphics[width=.5\textwidth]{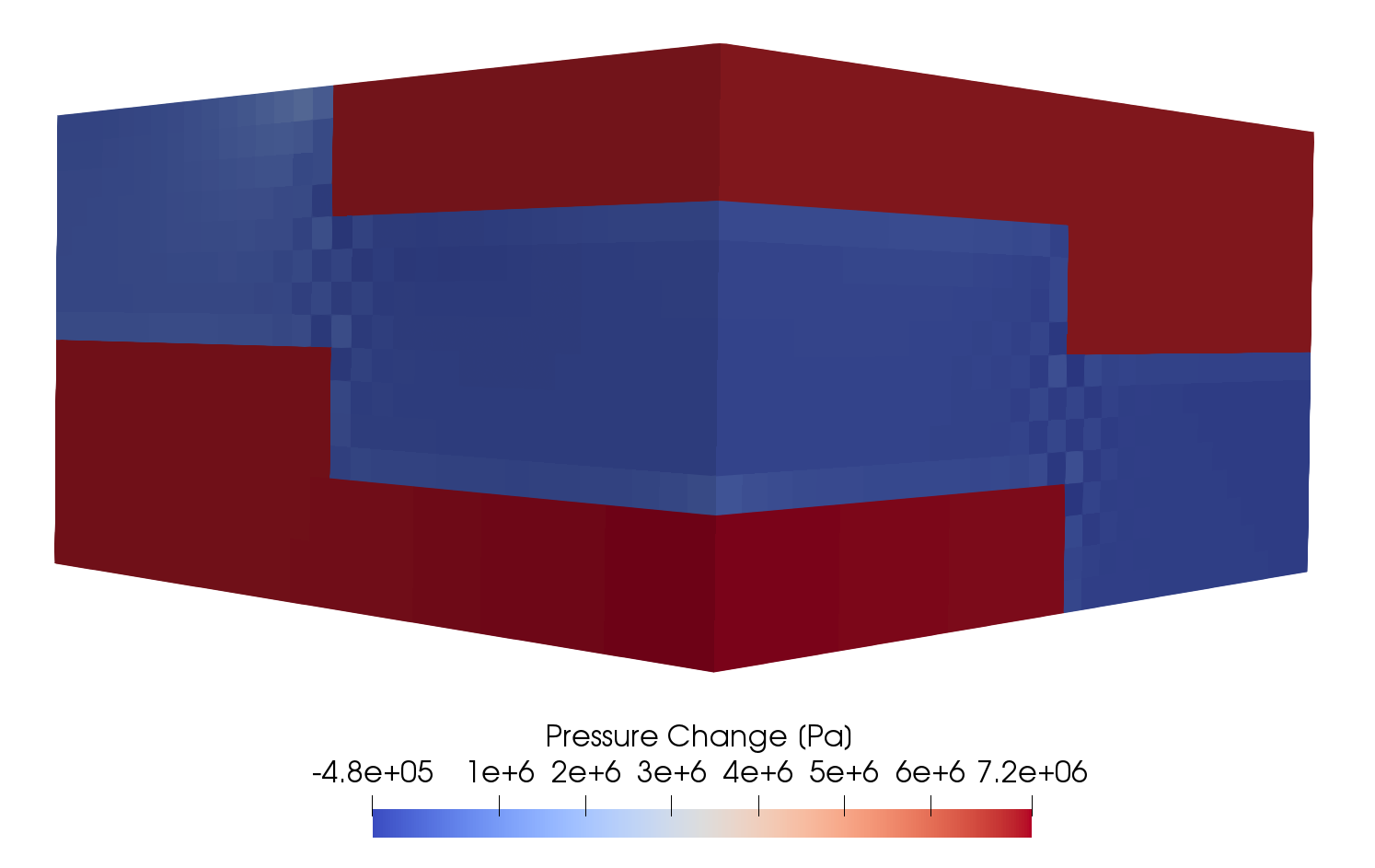}}\\
\subfloat[Stabilization in full domain]{\label{sfig:stair_full_30days}\includegraphics[width=.5\textwidth]{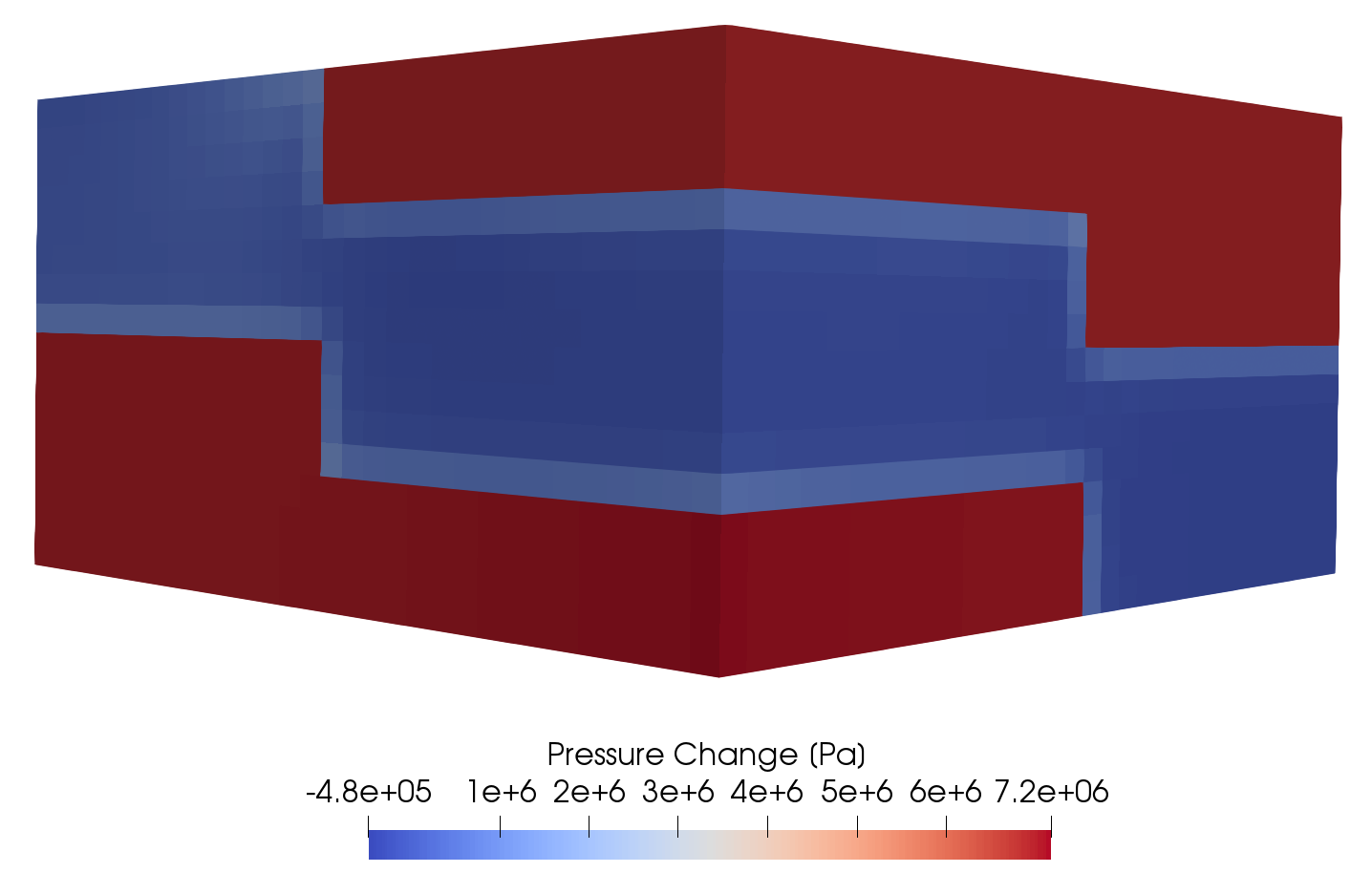}}\hfill
\subfloat[Stabilization in burden only]{\label{sfig:stair_burden_30days}\includegraphics[width=.5\textwidth]{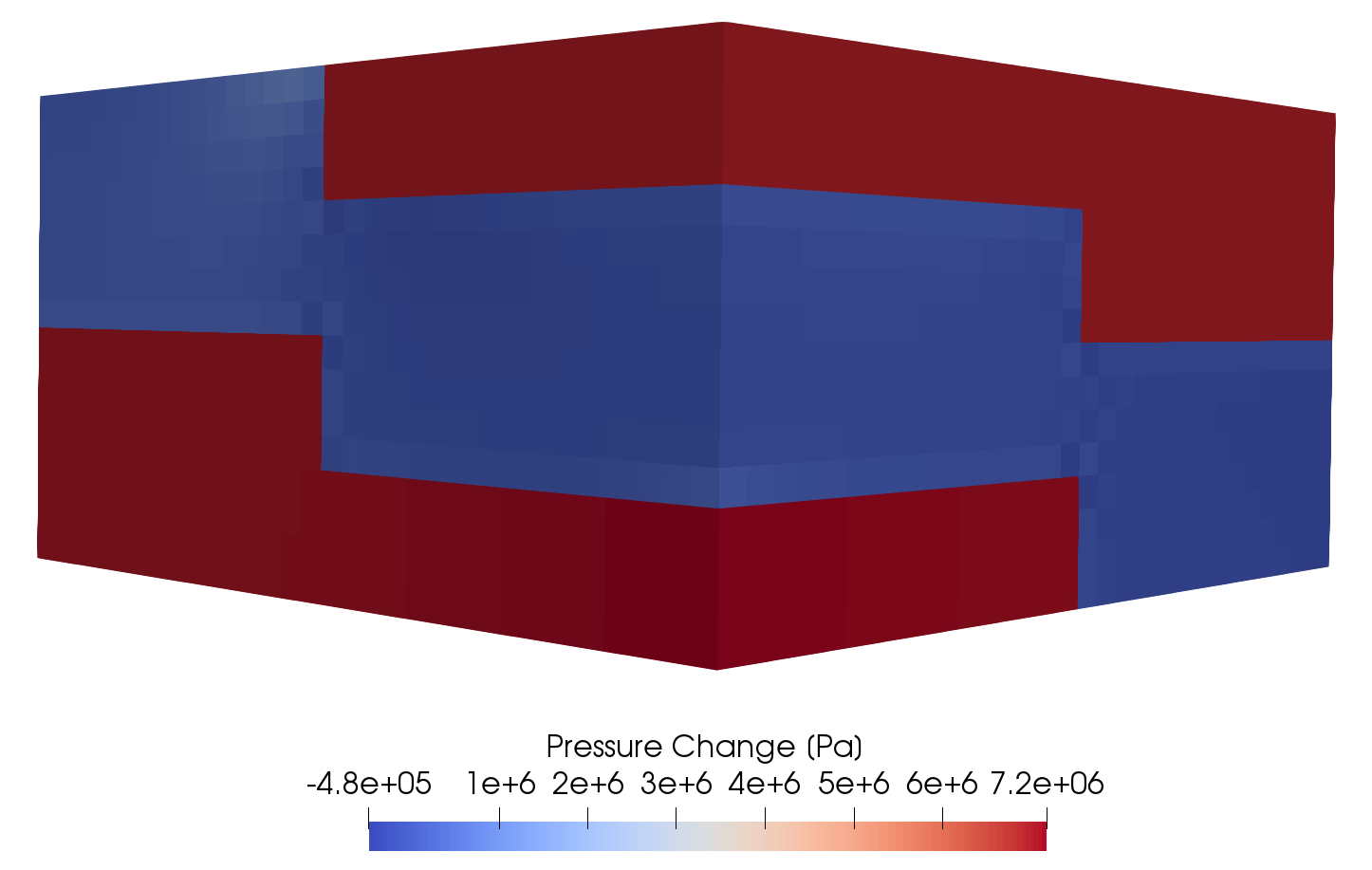}}\\
\caption{Staircase: Pressure field after 30 days of injection with undrained burden region}
\label{fig:stair_30days}
\end{figure}

Next we consider the results of the simulation after taking 30 time steps of size equal to one day, with a barrier permeability of $9.8\times 10^{-20}$ m$^2$ and $\alpha = 1$. Figure \ref{sfig:stair_nostab_30days} shows the resulting pressure field at the end of the simulation, and it is clear that spurious pressure oscillations are present in the barrier region. For reference, Figure \ref{sfig:stair_nostab_30days_fim} shows pressure resulting from an equivalent fully implicit simulation, which demonstrates that the fixed-stress scheme does recover the exact same spurious pressure modes in the undrained case. Figures \ref{sfig:stair_full_30days} and \ref{sfig:stair_burden_30days} show the pressure fields resulting from the fixed-stress scheme with stabilization in the full domain and in the barrier region only, respectively. These results match well with the fully implicit results seen in \cite{aronson2023pressure}, where both are effective at removing the spurious pressure mode in the undrained region. Applying stabilization across the interface between the low and high permeability regions results in a smearing of the interface, while just applying the stabilization flux in the burden retains the sharpness of the interface while still smoothing the spurious oscillations.  

\begin{figure}
\centering
\subfloat[No stabilization]{\label{sfig:plume_nostab_30days}\includegraphics[width=.5\textwidth]{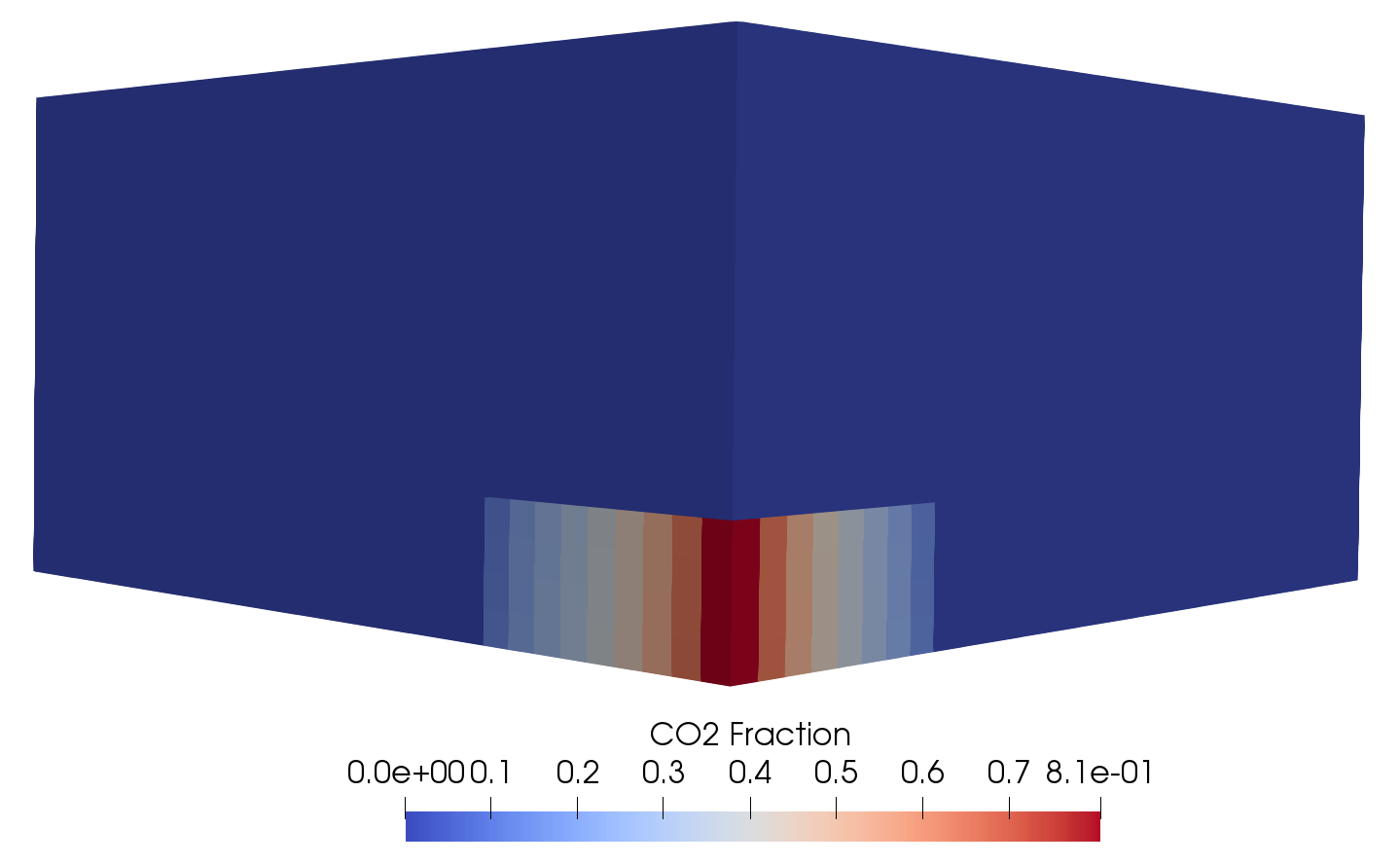}}\hfill
\subfloat[No stabilization (fully implicit)]{\label{sfig:plume_nostab_30days_fim}\includegraphics[width=.5\textwidth]{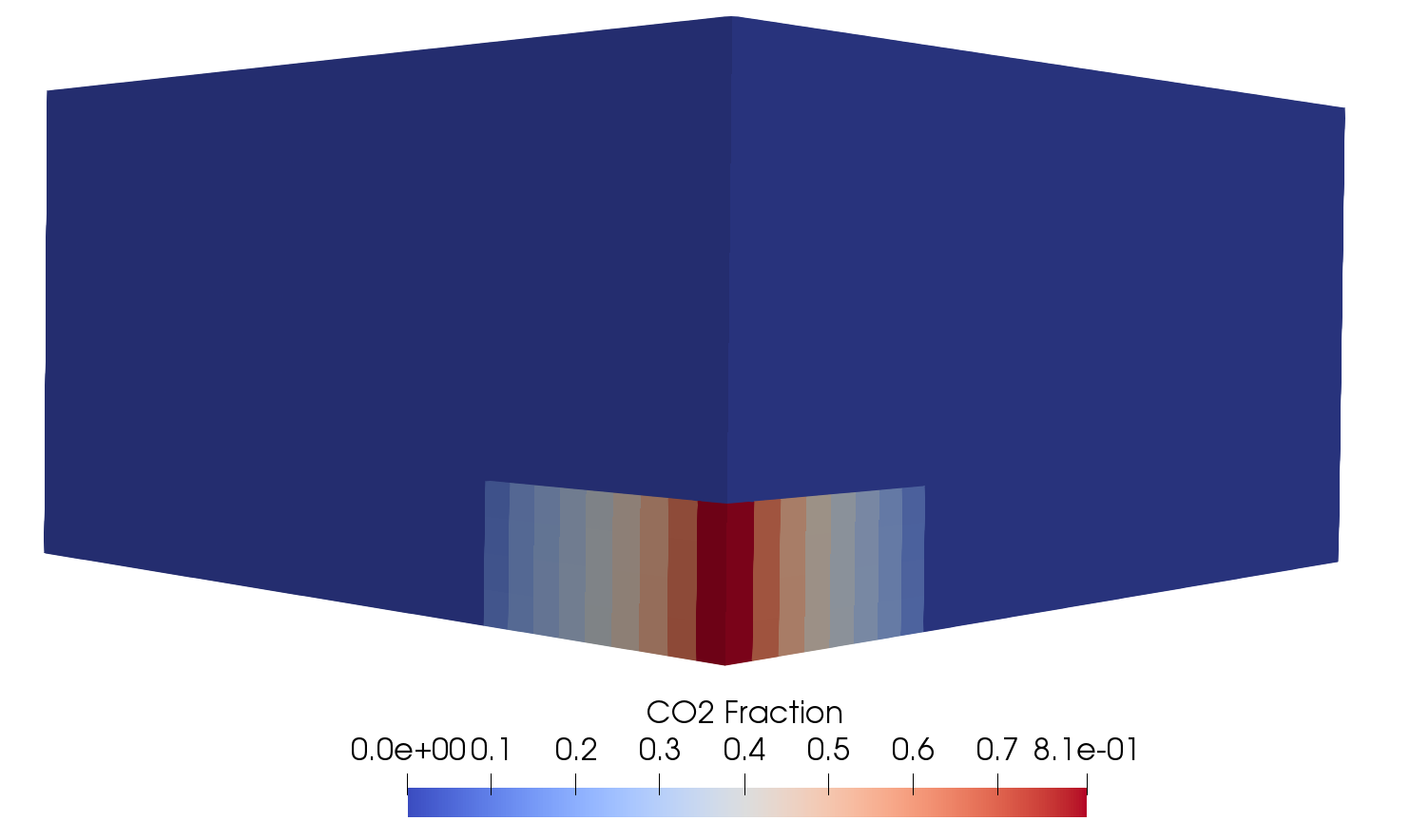}}\\
\subfloat[Stabilization in full domain, $c=1$]{\label{sfig:plume_full_30days}\includegraphics[width=.5\textwidth]{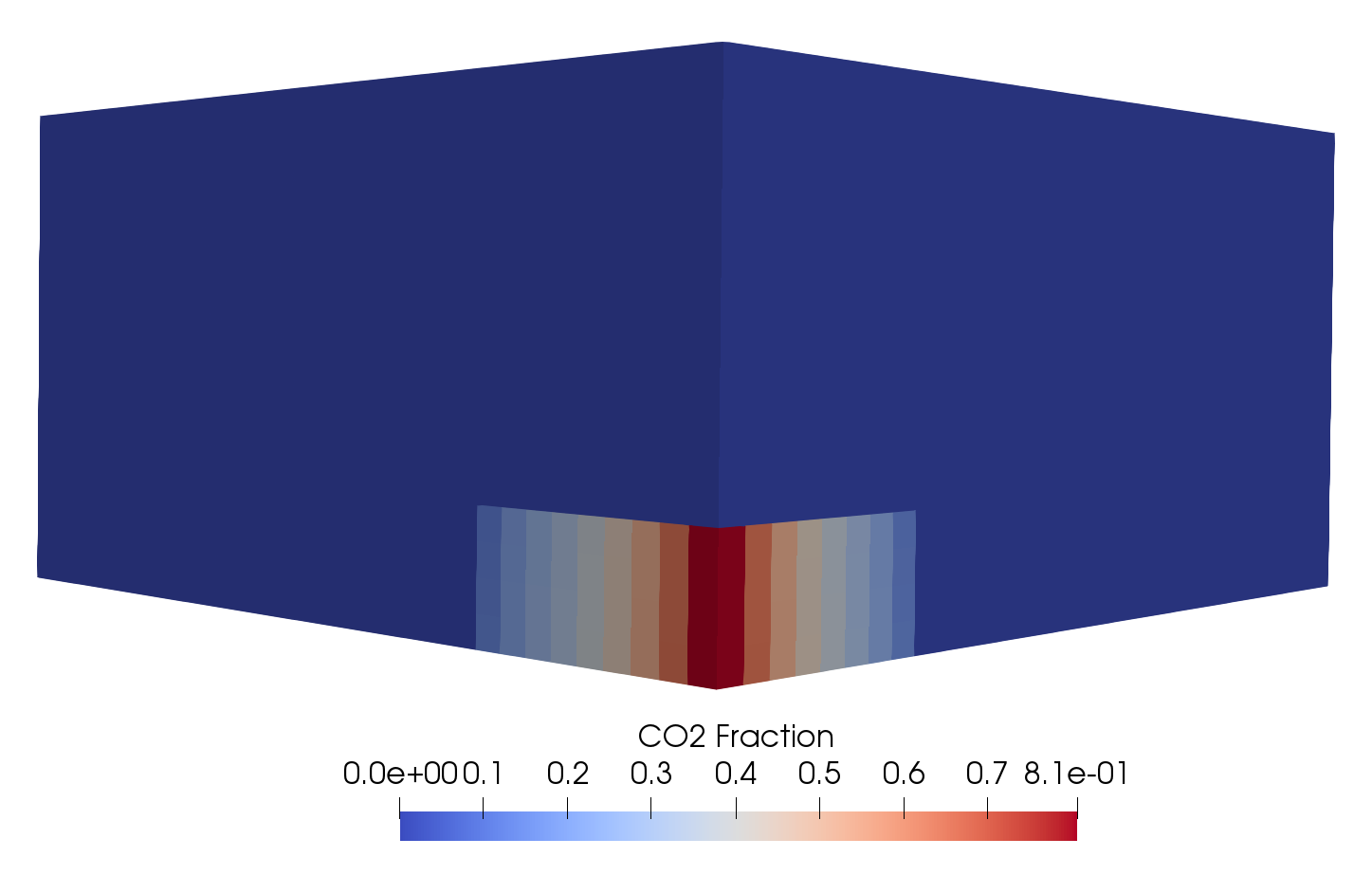}}\hfill
\subfloat[Stabilization in full domain, $c=25$]{\label{sfig:plume_fullC25_30days}\includegraphics[width=.5\textwidth]{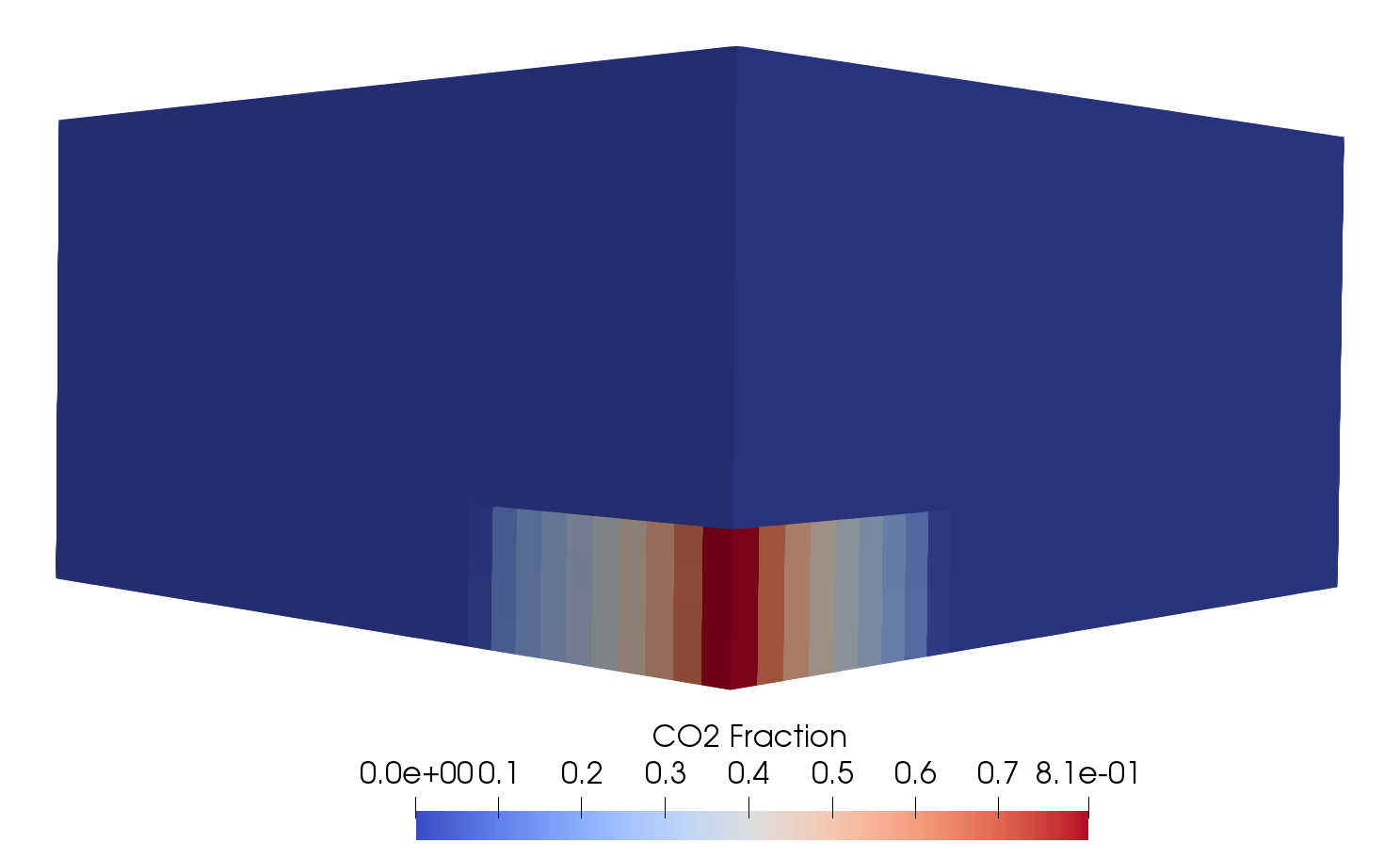}}\\
\caption{Staircase: CO$_2$-rich phase volume fraction after 30 days of injection with undrained burden region}
\label{fig:plume_30days}
\end{figure}

Finally, we can also consider the saturation fields at the end of the simulation, to demonstrate that the stabilization flux has a negligible effect on the quality of the results in the drained region. Figure \ref{fig:plume_30days} shows the volume fraction of CO$_2$-rich phase in the domain using various methodologies. In the first row we see that the results generated with the unstabilized fixed-stress method match those generated by a fully implicit simulation. Figure \ref{sfig:plume_full_30days} is obtained using stabilization everywhere in the domain with constant $c = 1$. We see that this also matches the unstabilized results, and in particular even when the interface between the two regions is smoothed by the stabilization this does not result in a artificial leakage of CO$_2$ into the barrier. In fact, Figure \ref{sfig:plume_fullC25_30days} shows that, even when this constant is increased to $c = 25$, the only result in that the plume travels slightly farther in the channel region (about one cell farther), and there is still no noticeable artificial leakages. If stabilization is applied only in the burden regions, there is certainly no effect on the CO$_2$ plume in the reservoir, by construction.   

\subsection{HI24L}

\begin{figure}
\centering
\subfloat[Full domain mesh]{\label{sfig:HI_full_setup}\includegraphics[width=.5\textwidth]{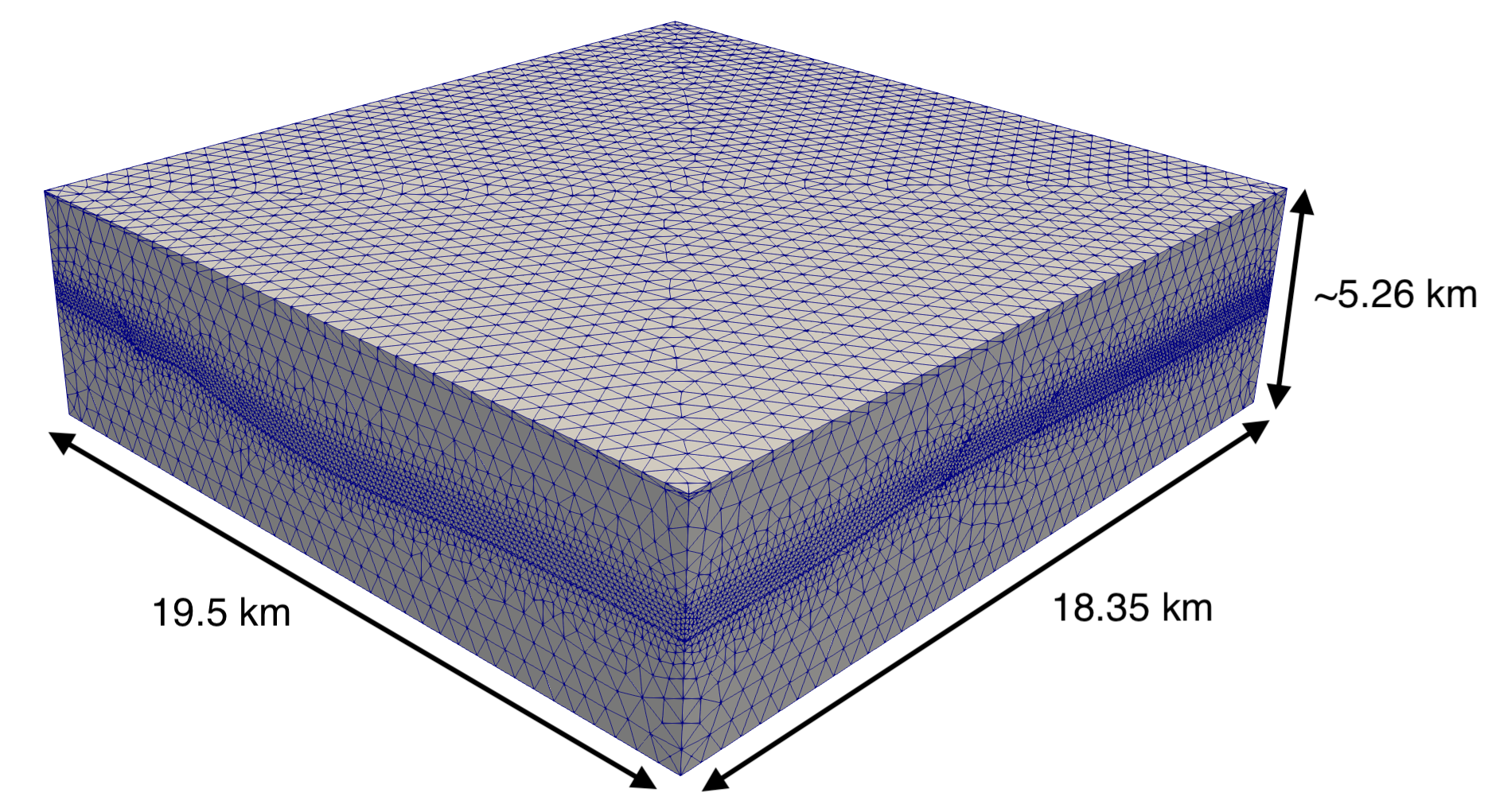}}\hfill
\subfloat[Reservoir mesh colored by permeability]{\label{sfig:HI_res_setup}\includegraphics[width=.5\textwidth]{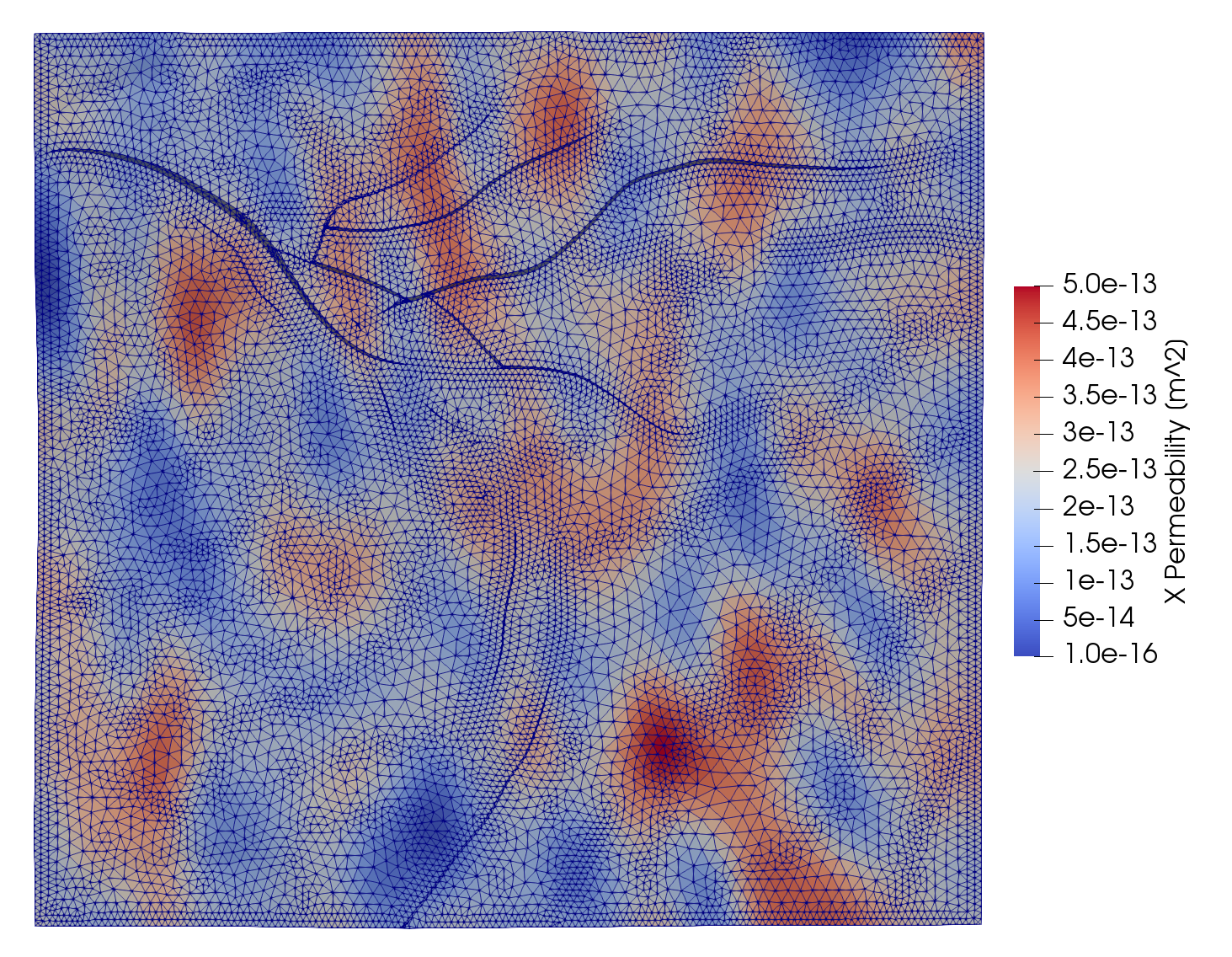}}\\
\caption{HI24L problem setup}
\label{fig:HI_setup}
\end{figure}

To consider the effect of inf-sup instabilities on a more realistic scenario, we now consider the simulation of CO$_2$ injection into a geologic model developed for an offshore field in the Gulf of Mexico, HI24L. The region around HI24L has been extensively studied for its large CO$_2$ storage potential, as summarized in \cite{ruiz2019characterization}. The full simulation mesh is shown in Figure \ref{sfig:HI_full_setup}, and contains a reservoir (shown in Figure \ref{sfig:HI_res_setup}) between two burden regions. Table \ref{tab:params_solid_HI} summarizes the material properties of the solid skeleton, and we note that we will use burden permeabilities of $1.0\times 10^{-14}$ m$^2$, $1.0\times 10^{-17}$ m$^2$, and $1.0\times 10^{-20}$ m$^2$ in our numerical tests to determine the behavior of the fixed-stress scheme from drained to undrained conditions. The fluid parameters are again summarized in Table \ref{tab:params_fluid}. We inject CO$_2$ at a rate of 30 kg/s into the reservoir near the center in the horizontal directions, and at the bottom of the reservoir in the vertical direction.

\begin{table}[ht]
\centering
\caption{Solid material parameters used in HI24L}
\begin{tabular}{ccc}
\hline
{ \textbf{Parameter}} & { \textbf{Reservoir Value}} & { \textbf{Burden Value}}  \\ \hline
Porosity                 & Varies (0.05 - 0.25, correlated with permeability)                            & 0.05                         \\
Permeability             & Varies (Figure \ref{sfig:HI_res_setup})                       & $1.0\times 10^{-14}$ - $1.0\times 10^{-20}$ m$^2$                    \\
Skeleton bulk modulus       & 9.4 GPa                          & 11.5 GPa                        \\
Skeleton Poisson ratio  & 0.25  & 0.3 \\
Skeleton Biot coefficient   & 1                              & 1                            \\
Skeleton density            & 2700 kg/m$^3$   & 2700 kg/m$^3$                             
\end{tabular}
\label{tab:params_solid_HI}
\end{table}

\begin{figure}
\centering
\subfloat[No stabilization]{\label{sfig:HI_nostab_FS_iter}\includegraphics[width=.5\textwidth]{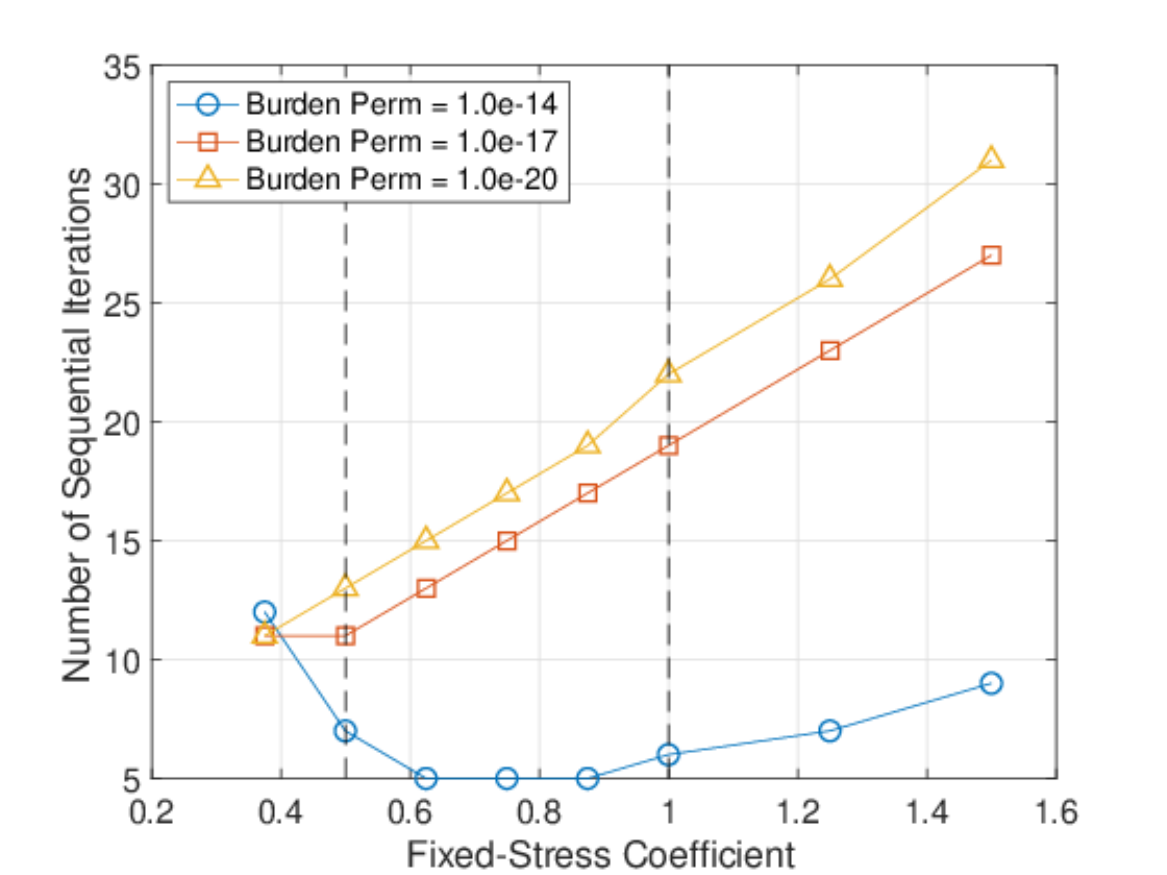}}\hfill
\subfloat[Burden only stabilization]{\label{sfig:HI_stab_FS_iter}\includegraphics[width=.5\textwidth]{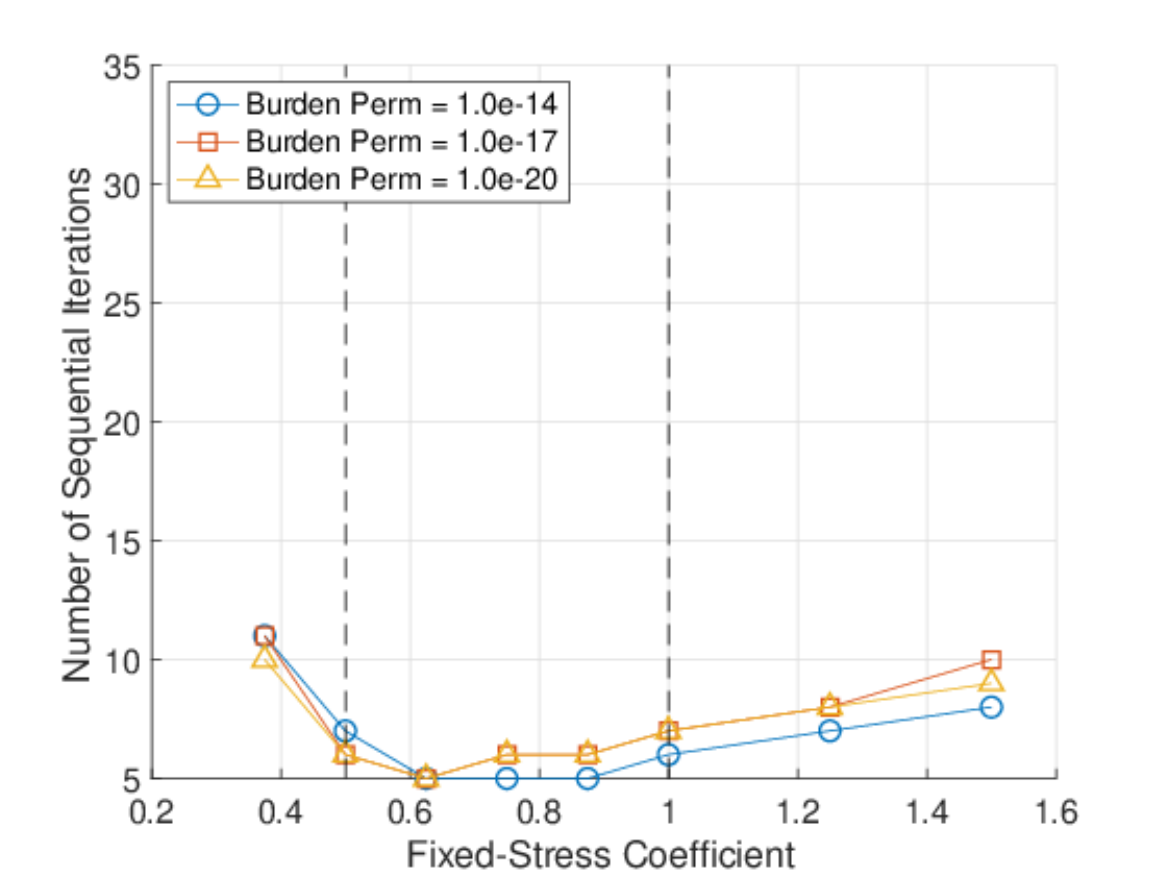}}\\
\caption{HI24L: Number of fixed-stress iterations with varying fixed-stress coefficients $\alpha$}
\label{fig:HI_iter_FScoeff}
\end{figure}

As in the previous example, we start by considering the convergence of the iterative fixed-stress scheme in the first time step. For brevity we only consider an initial time step size of one day, and Figure \ref{fig:HI_iter_FScoeff} details the results for various burden permeability values and values of $\alpha$. Without stabilization, the number of iterations required to reach a tolerance of $1.0\times 10^{-8}$ grows as the burden region becomes undrained, but the inclusion of stabilization almost completely removes this effect. Note once again that the stabilized results in Figure \ref{fig:HI_iter_FScoeff} are generated using stabilization only in the burden region, but the results when stabilization is included throughout the domain are similar.

\begin{figure}
\centering
\includegraphics[width=.5\textwidth]{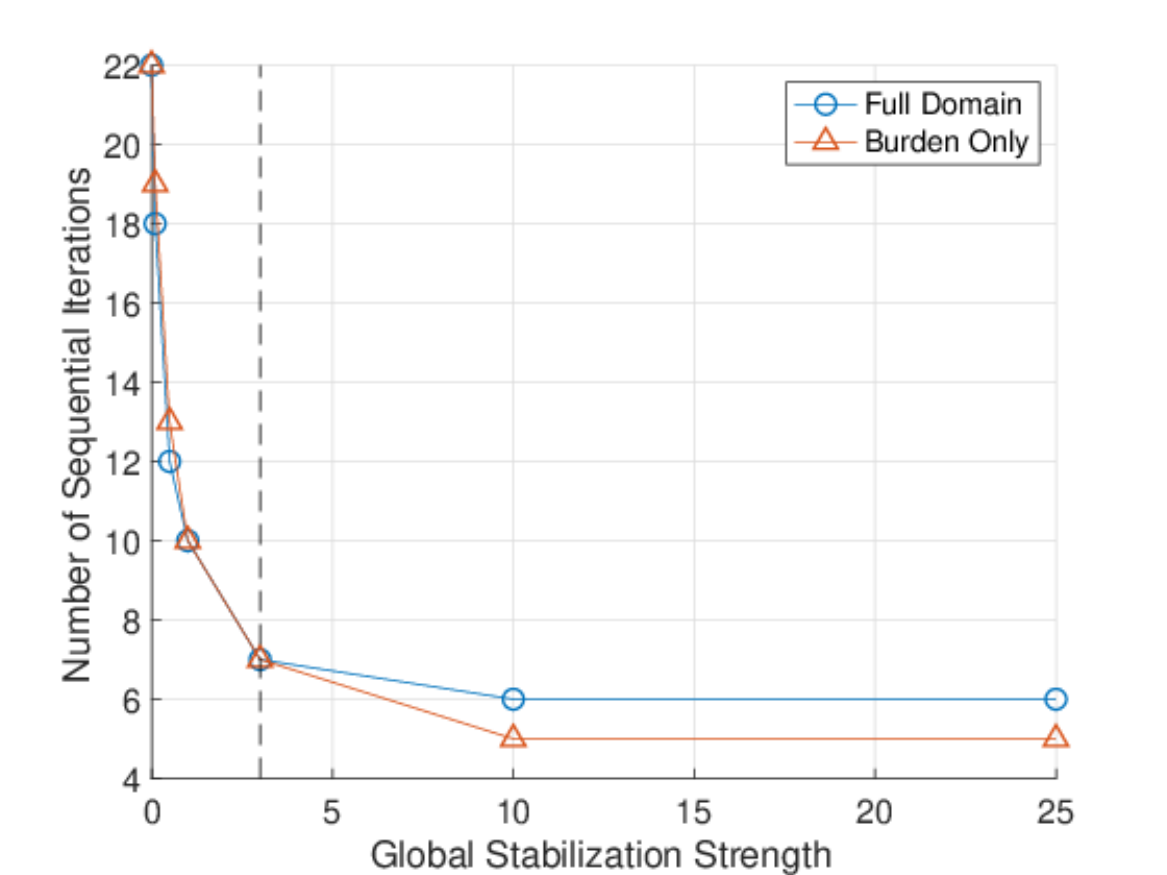}
\caption{HI24L: Number of fixed-stress iterations with varying pressure stabilization coefficient $c$ and undrained burden regions}
\label{fig:stab_strength_hi}
\end{figure}

Fixing $\alpha = 1$ and burden permeability equal to $1.0\times 10^{-20}$ m$^2$ during the initial time step of one day, we can again study the effect of the stabilization strength on the convergence rate. Figure \ref{fig:stab_strength_hi} shows the number of iterations required for convergence with varying values of the constant $c$. The dashed vertical line represents the value selected according to the previous section, $c = 3$. Much like in the previous example, this constant seems to strike a good balance between stability and accuracy, as it is once again in the knee of the plot and thus seems to provide the best efficiency without risking over-smoothing the physical solution. 

\begin{figure}
\centering
\includegraphics[width=0.5\textwidth]{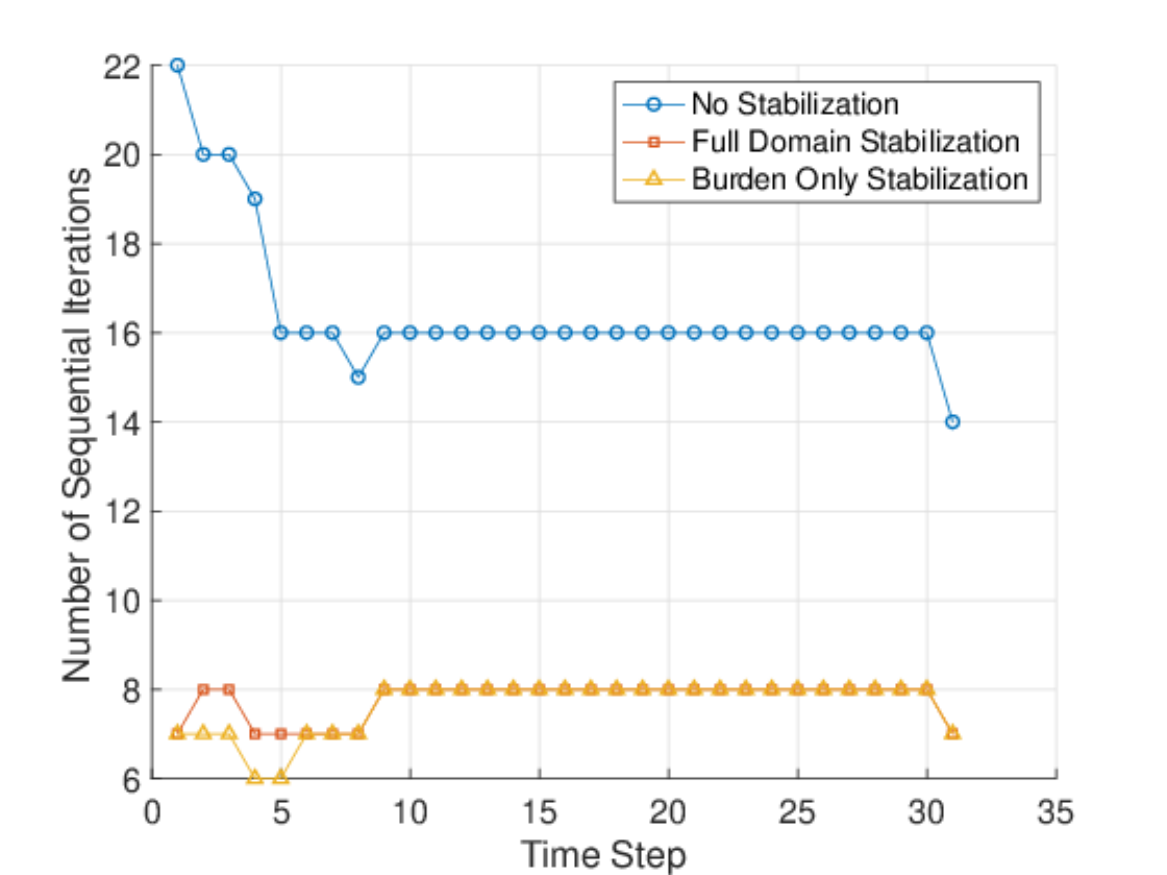}
\caption{HI24L: Number of fixed-stress iterations per time step during 2 year simulation with undrained burden regions}
\label{fig:timestep_hi}
\end{figure}

Next we consider simulating injection for two years, starting with a time step of one day and allowing this to double until a maximum time step size of 30 days. Figure \ref{fig:timestep_hi} shows the number of fixed-stress iterations taken at each time step in the simulations with and without stabilization. Even as the time step size grows as time advances, the inclusion of stabilization reduces the number of iterations by approximately half. 

\begin{figure}
\centering
\subfloat[No stabilization]{\label{sfig:HI_nostab_FS}\includegraphics[width=.5\textwidth]{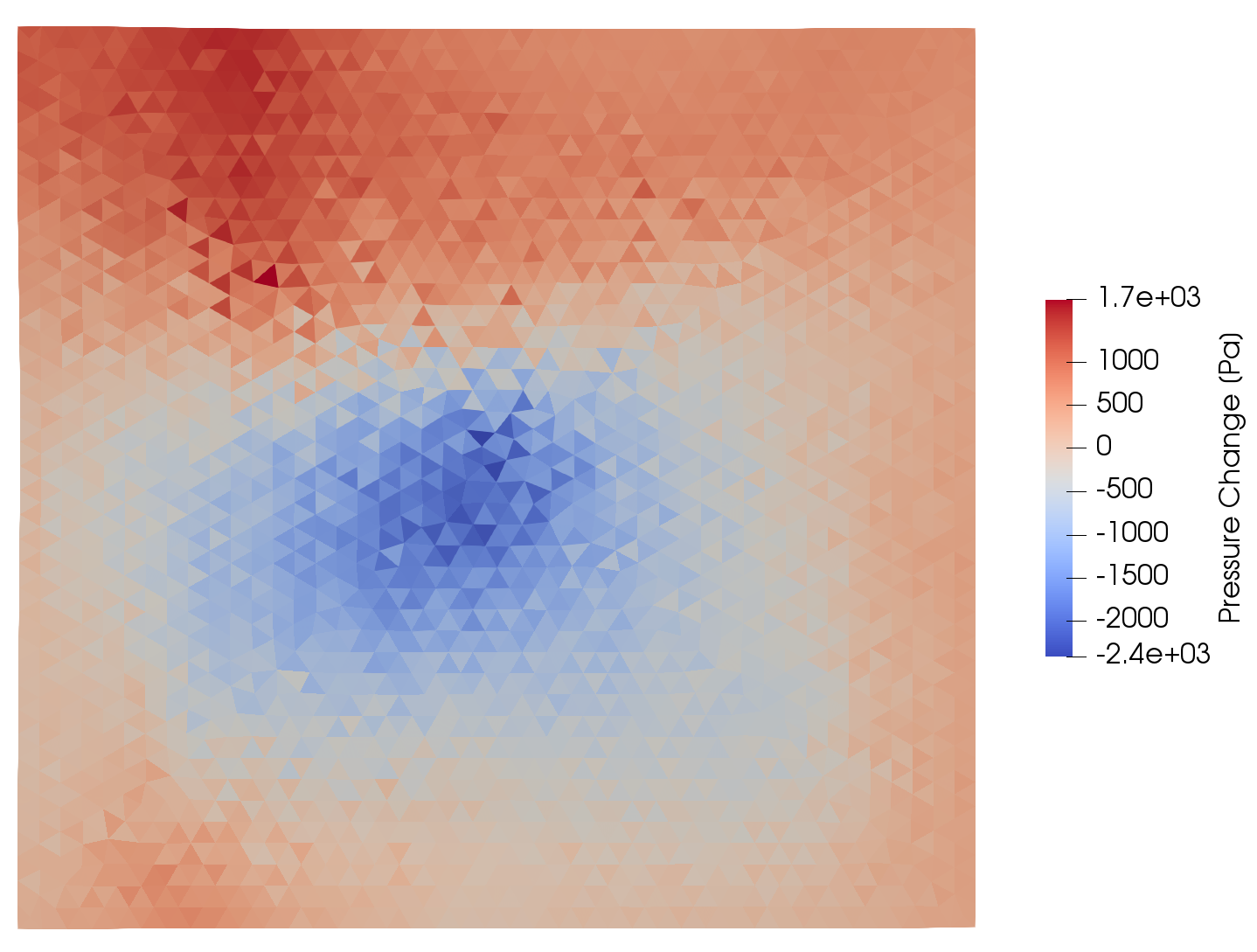}}\hfill
\subfloat[No stabilization (fully implicit)]{\label{sfig:HI_nostab_FS_fim}\includegraphics[width=.5\textwidth]{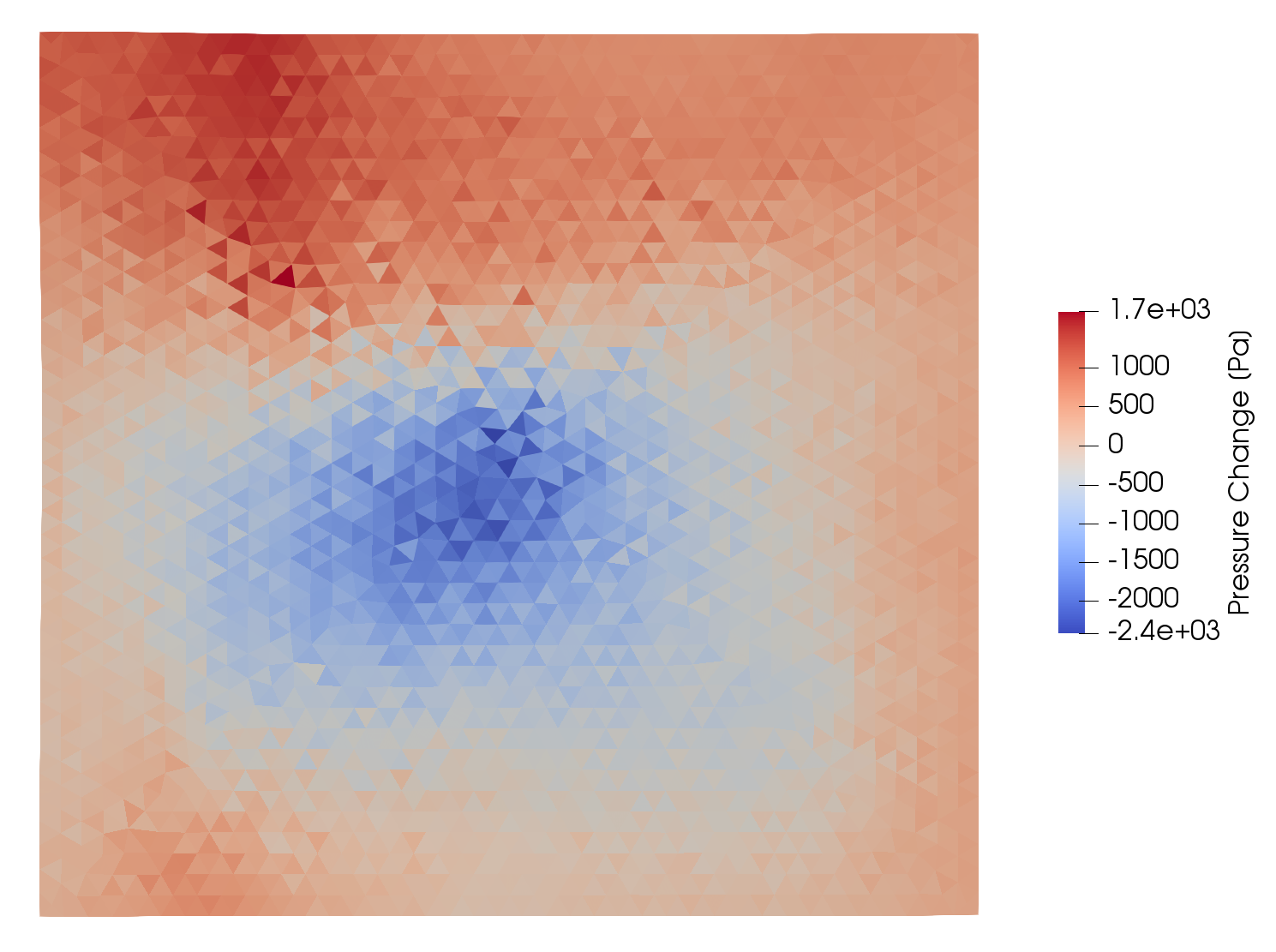}}\\
\subfloat[Full domain stabilization]{\label{sfig:HI_full_FS}\includegraphics[width=.5\textwidth]{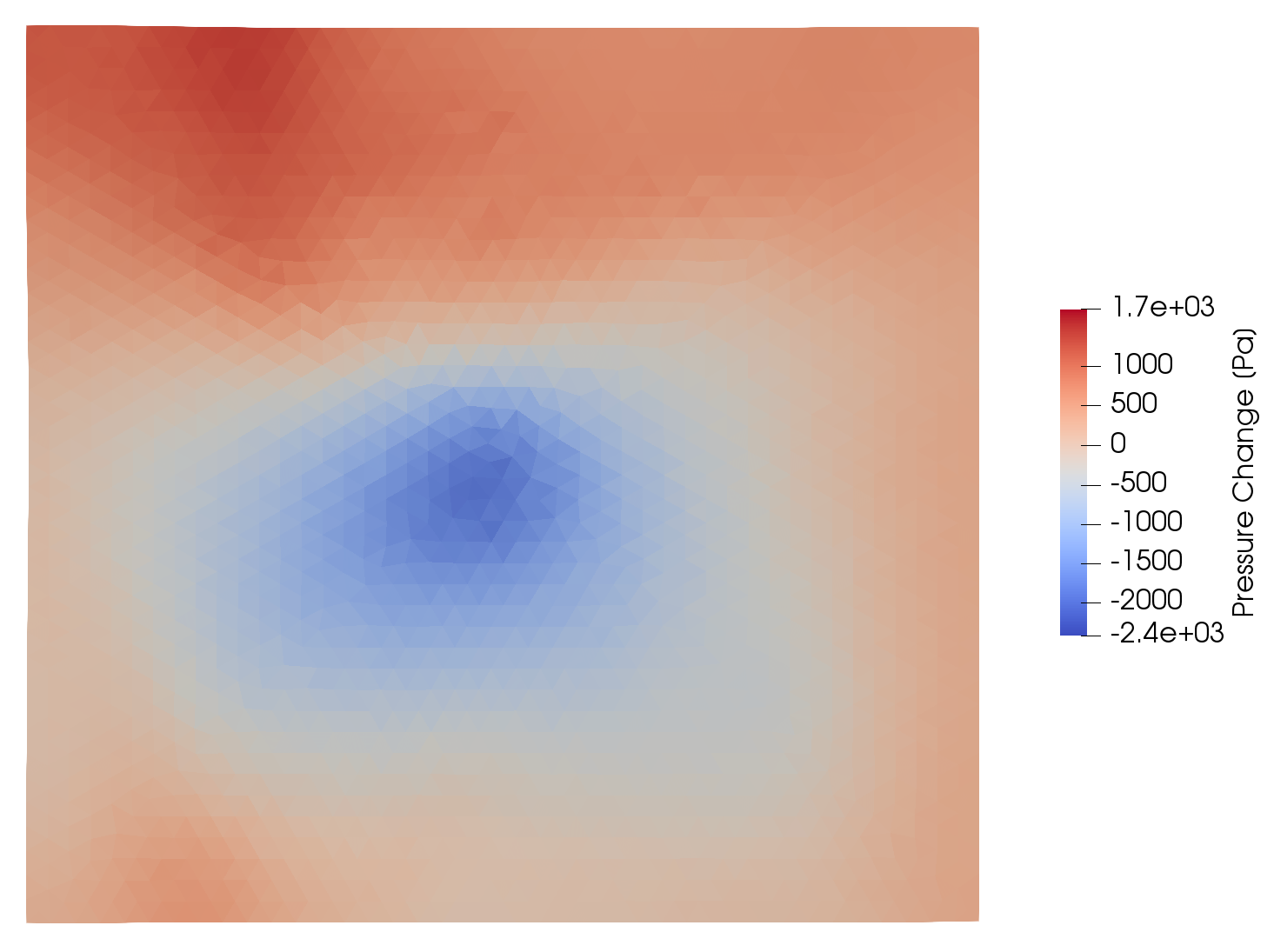}}\hfill
\subfloat[Burden only stabilization]{\label{sfig:HI_burden_FS}\includegraphics[width=.5\textwidth]{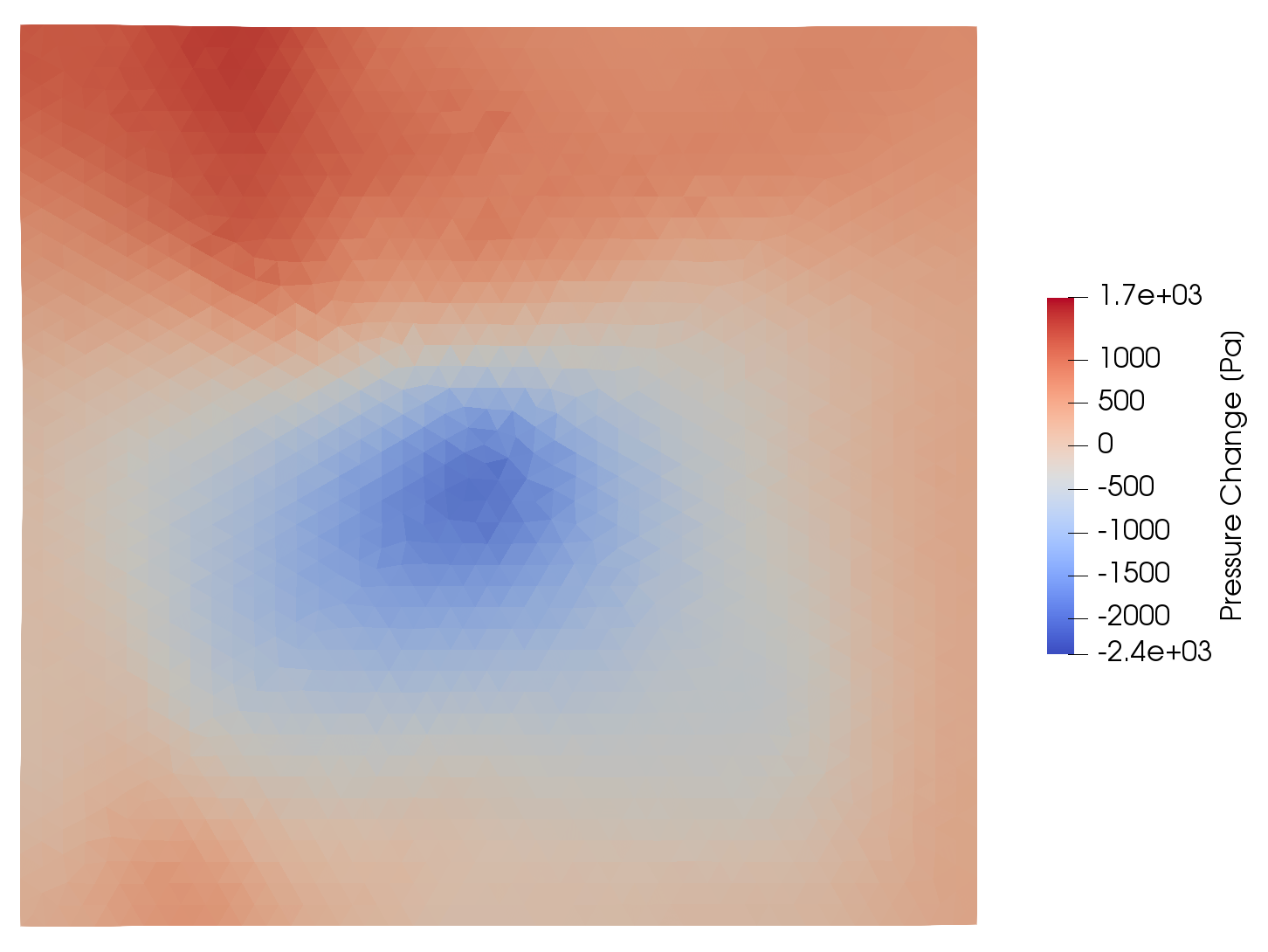}}\\
\caption{HI24L: Seafloor pressure field after 2 years of injection with undrained burden regions}
\label{fig:HI_pressure}
\end{figure}

\begin{figure}
\centering
\subfloat[No stabilization]{\label{sfig:HI_slice_nostab}\includegraphics[width=.7\textwidth]{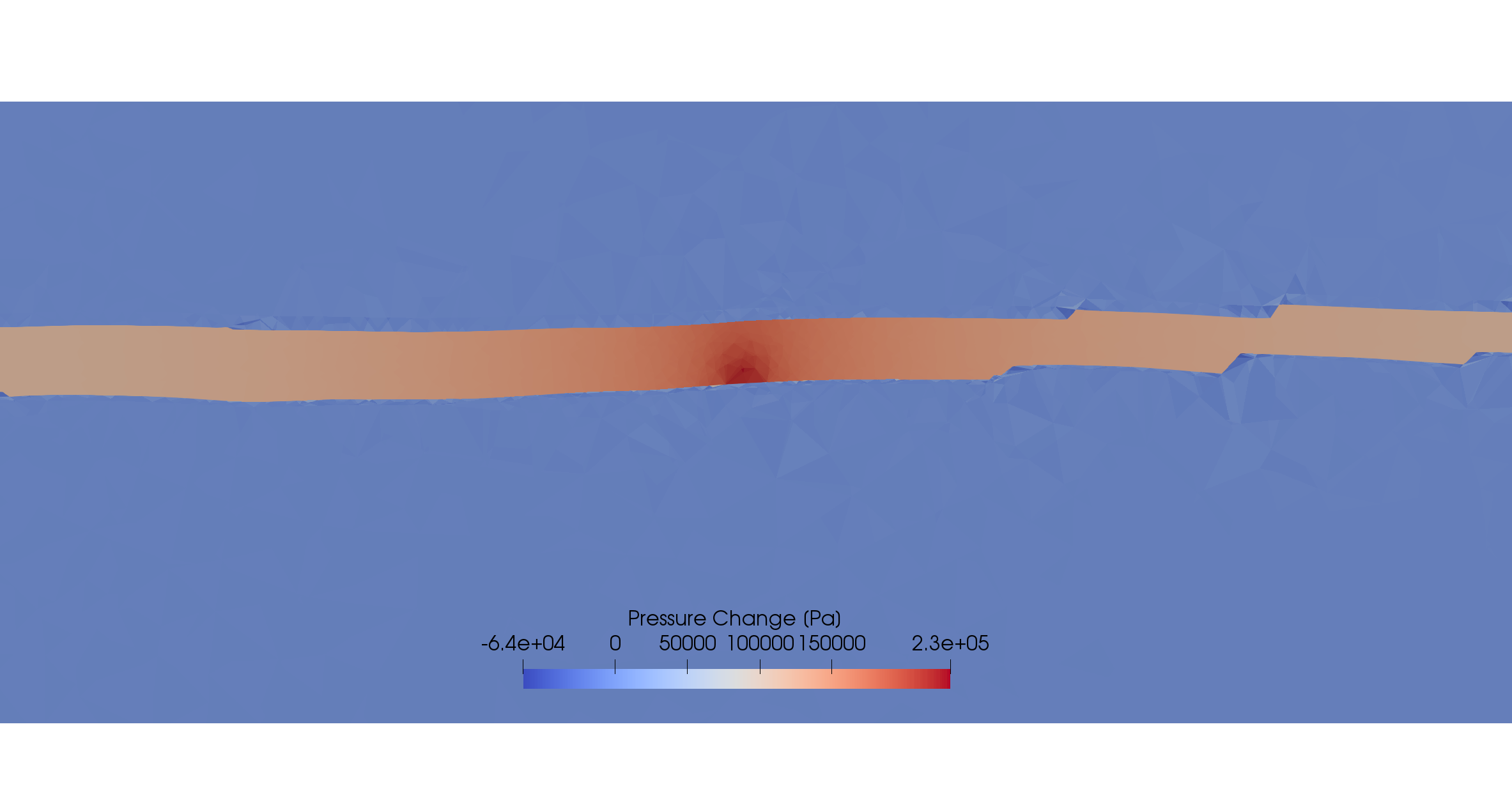}}\\
\subfloat[Stabilization in full domain]{\label{sfig:HI_slice_full}\includegraphics[width=.7\textwidth]{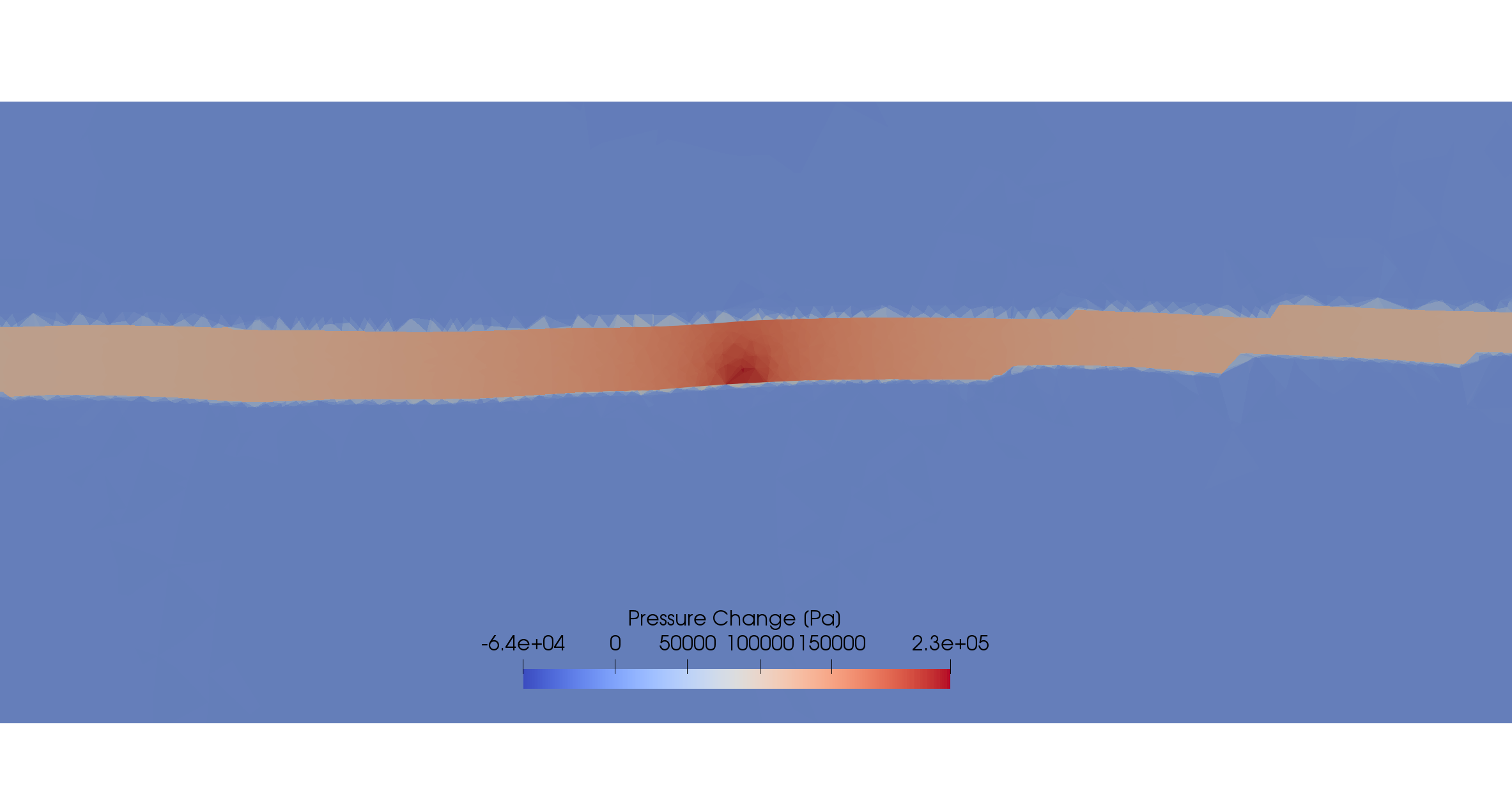}}\\
\subfloat[Stabilization in burden only]{\label{sfig:HI_slice_burden}\includegraphics[width=.7\textwidth]{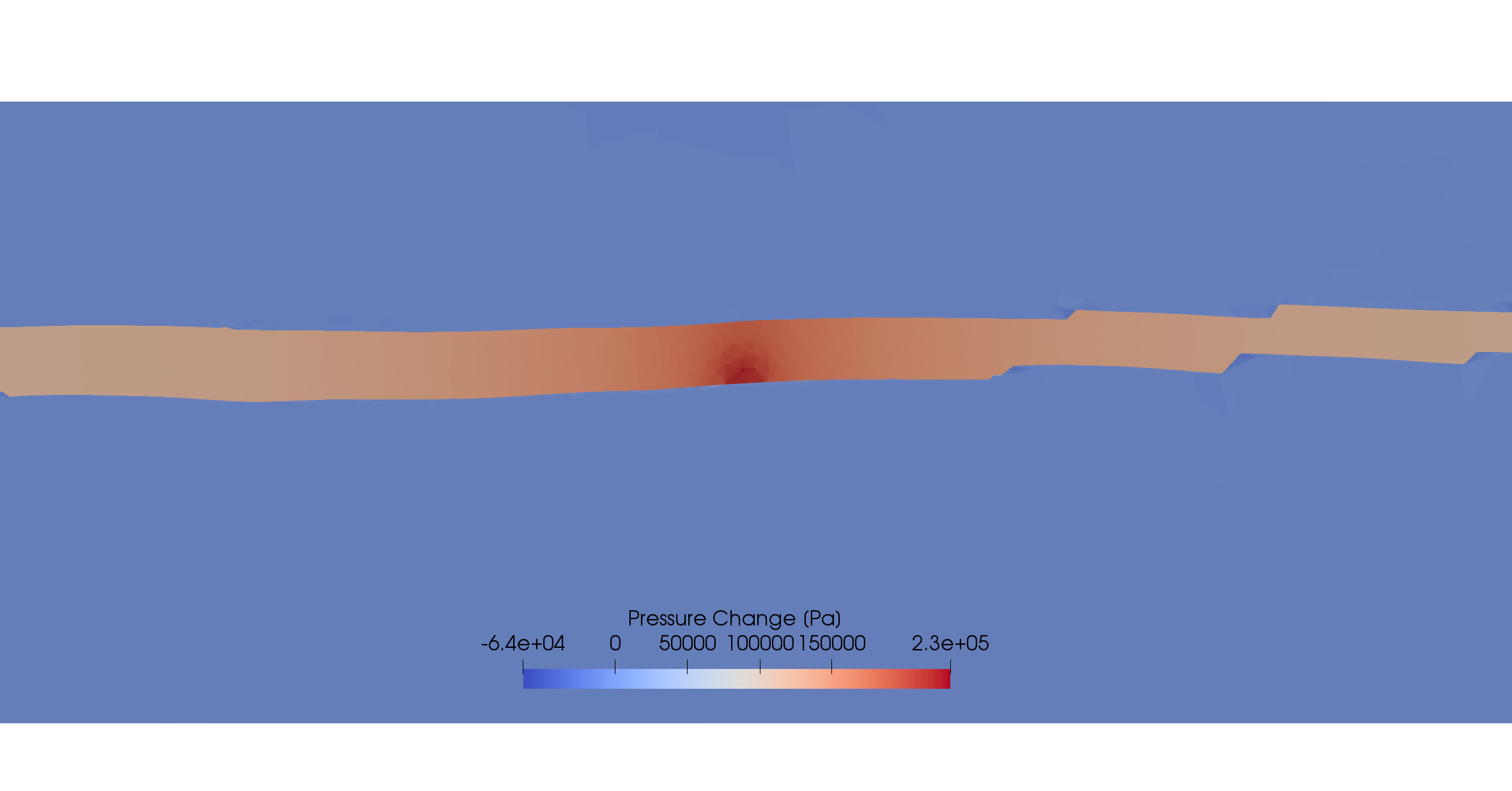}}\\
\caption{HI24L: Vertical slice of pressure field after 2 years of injection with undrained burden regions}
\label{fig:slice_HI}
\end{figure}

Figure \ref{fig:HI_pressure} shows the pressure field at the seafloor, or the top surface of the overburden region, after two years of injection. Clearly there are spurious pressure oscillations superimposed on the physical solution when no stabilization is applied, and the oscillations are the same in the fixed-stress and fully implicit results. Including stabilization effectively removes the spurious oscillations from the solution, and at the seafloor there does not seem to be a significant difference between the solutions obtained by applying stabilization throughout the domain and applying stabilization only in the burden regions. To see the difference between these two stabilization strategies, we consider the pressure field in a vertical slice of the domain near the injection location, shown in Figure \ref{fig:slice_HI}. Without stabilization there are oscillations throughout the burden region, and these are smoothed with stabilization. However, we see that the stabilization flux again smooths the pressure at the interface between the reservoir and the burden regions, while this is remedied by the burden only strategy. 

\begin{figure}
\centering
\subfloat[No stabilization]{\label{sfig:HI_slice_nostab_plume}\includegraphics[width=.65\textwidth]{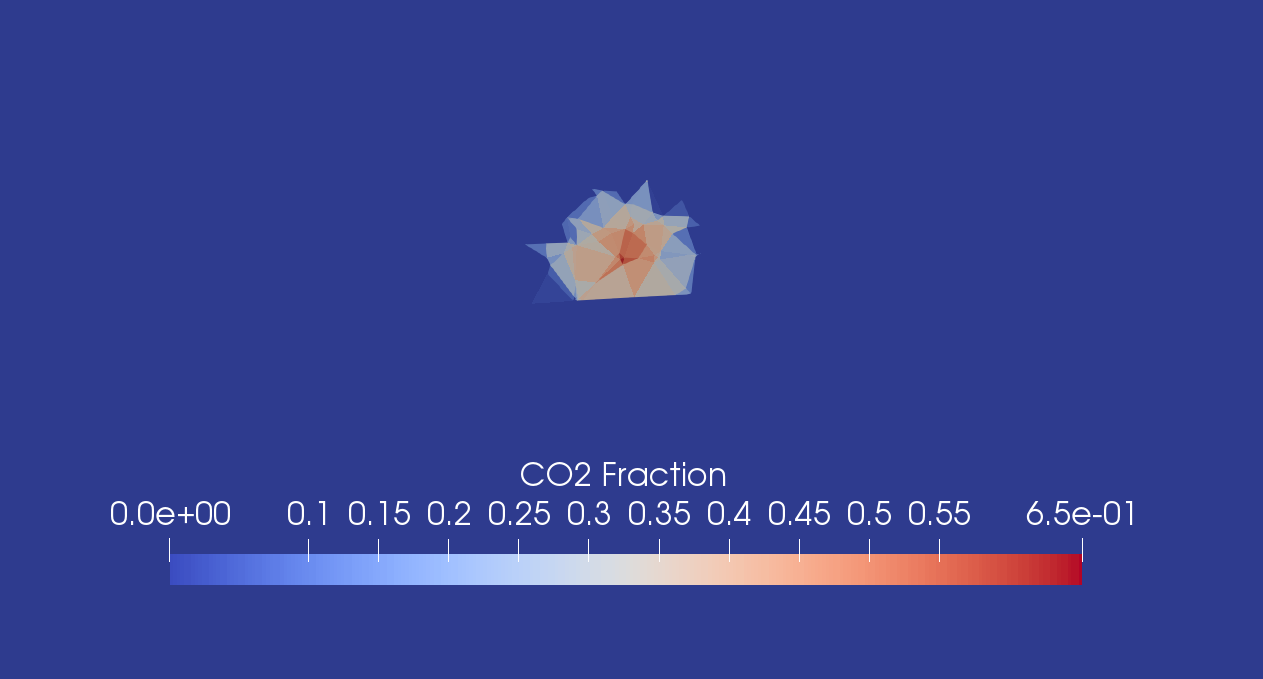}}\\
\subfloat[Stabilization in full domain, $c = 3$]{\label{sfig:HI_slice_full_plume}\includegraphics[width=.65\textwidth]{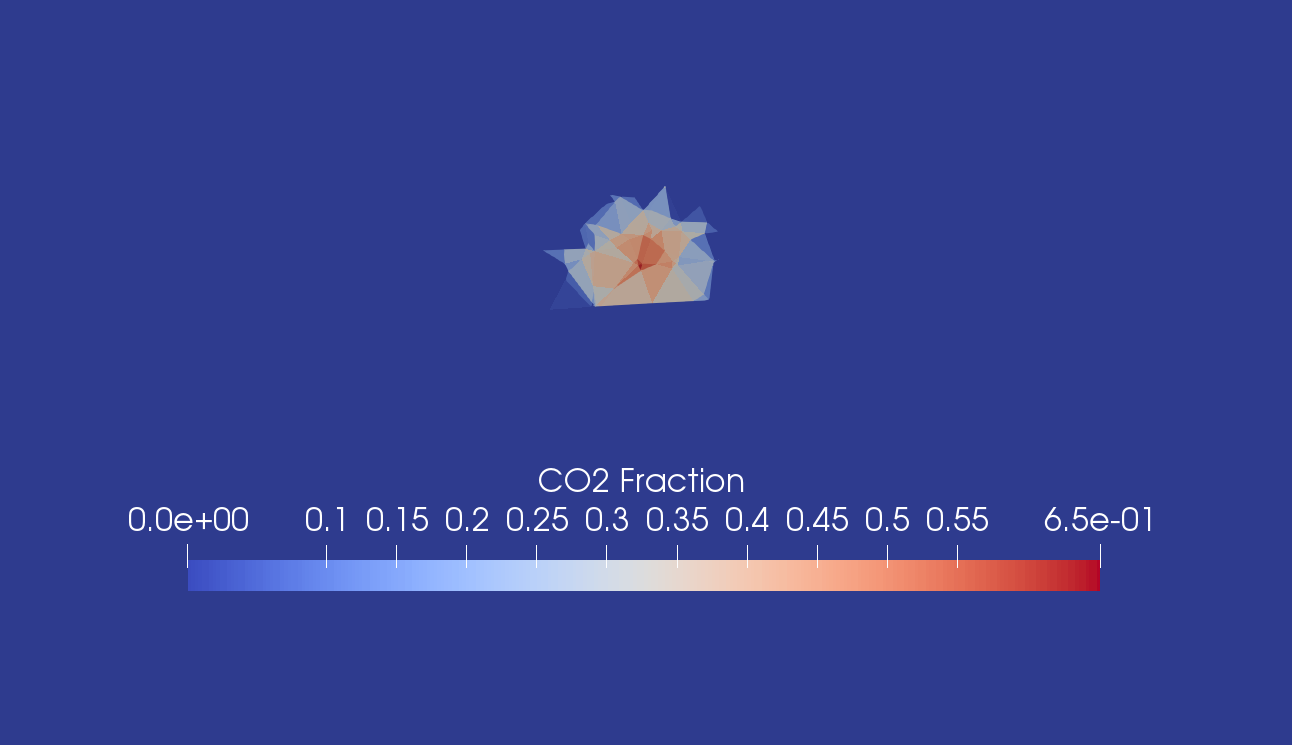}}\\
\subfloat[Stabilization in full domain, $c = 25$]{\label{sfig:HI_slice_fullc25_plume}\includegraphics[width=.65\textwidth]{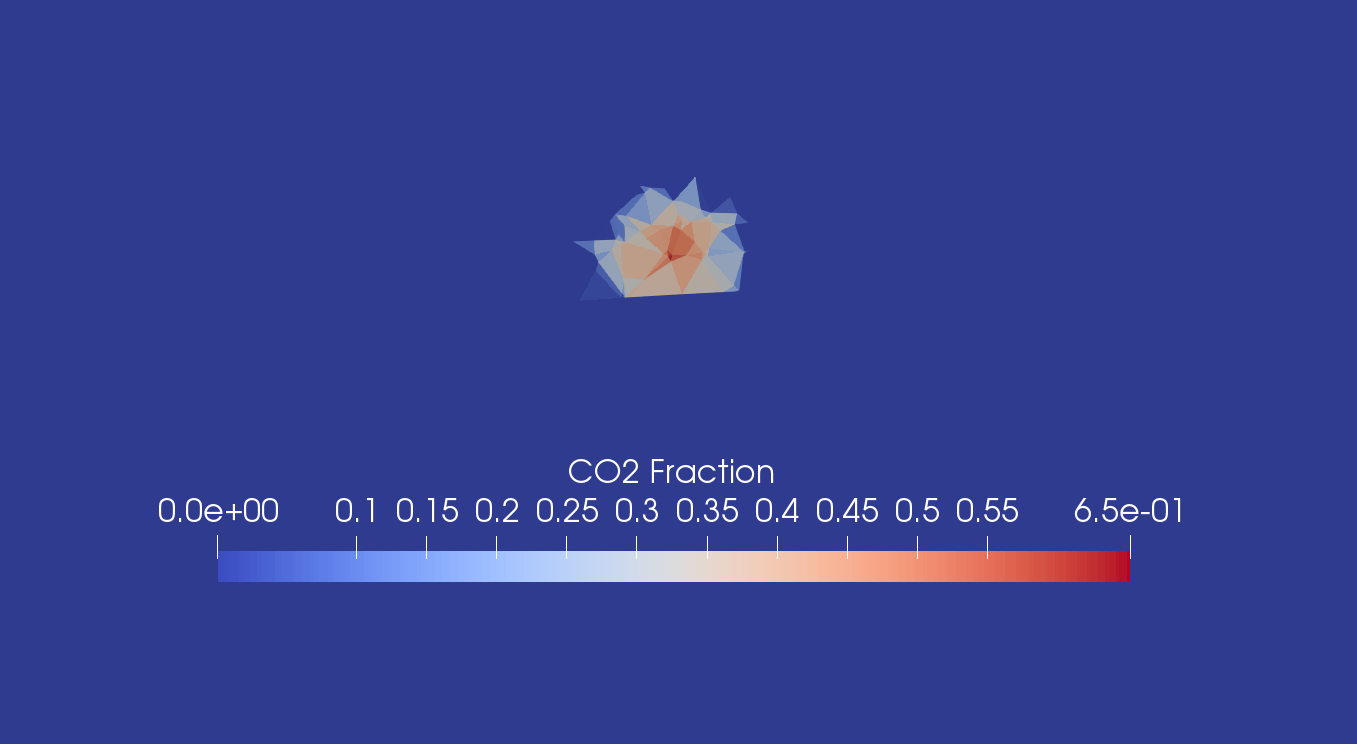}}\\
\caption{HI24L: Vertical slice of CO$_2$-rich phase volume fraction after 2 years of injection with undrained burden regions}
\label{fig:slice_HI_plume}
\end{figure}

\begin{figure}
\centering
\subfloat[Stabilization in full domain]{\label{sfig:HI_uplift_full}\includegraphics[width=.5\textwidth]{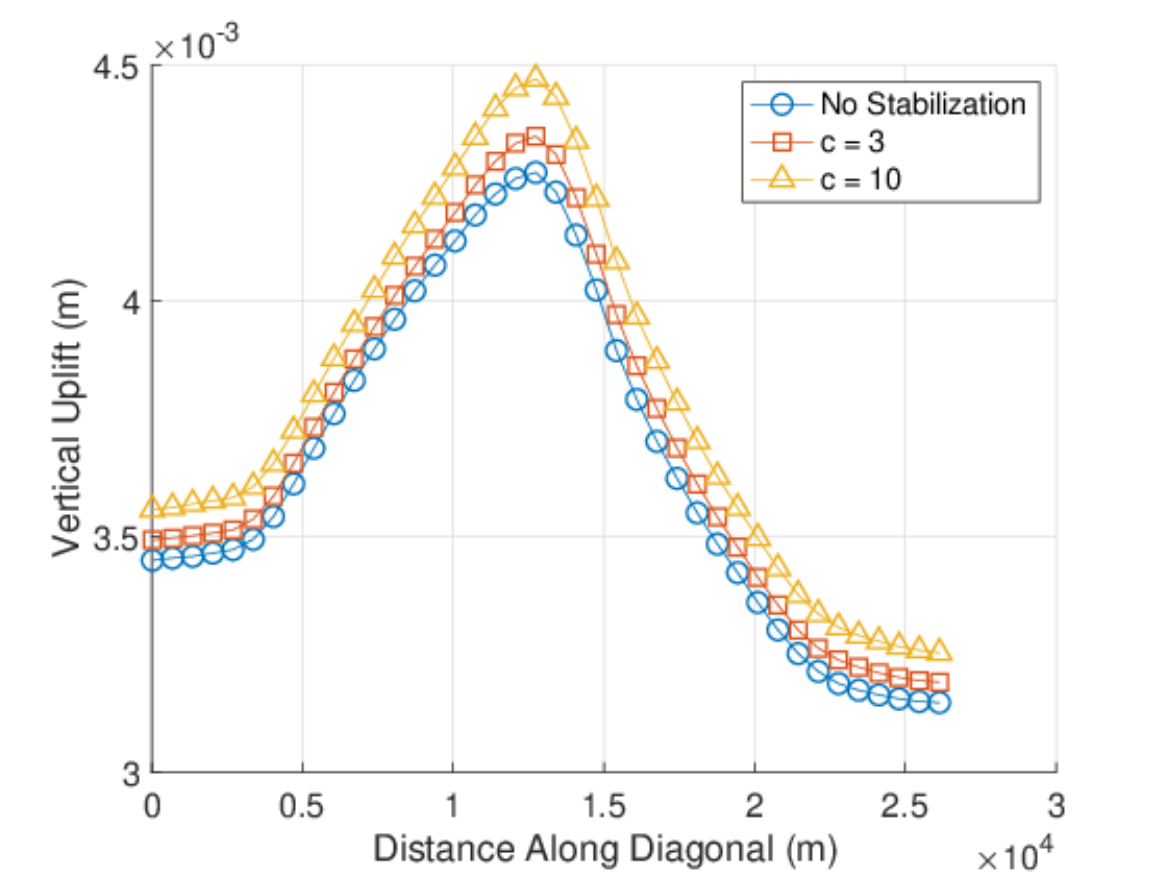}}\hfill
\subfloat[Stabilization in burden only]{\label{sfig:HI_uplift_burden}\includegraphics[width=.5\textwidth]{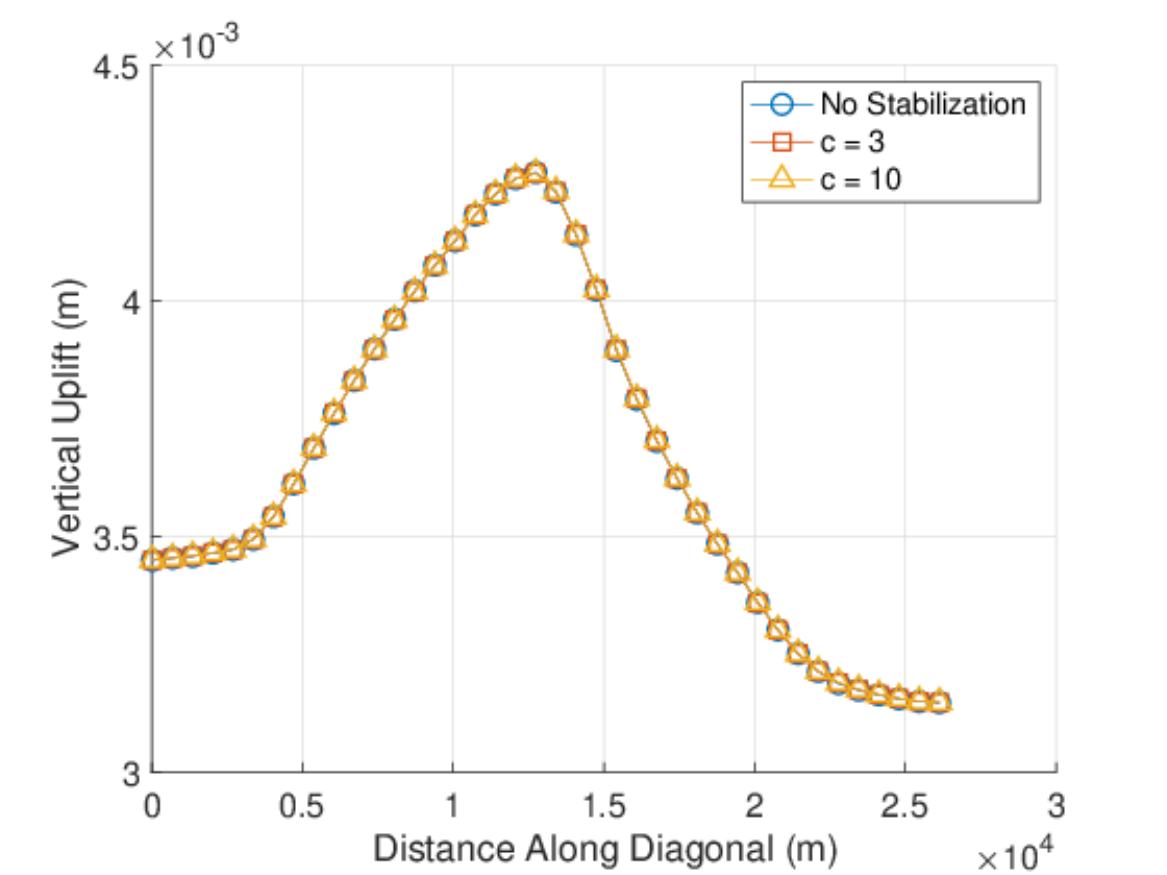}}\\
\caption{HI24L: Uplift of seafloor along horizontal diagonal after 2 years of injection with undrained burden regions}
\label{fig:uplift_HI}
\end{figure}

Note that, like the previous example, the smoothed pressure interface present when applying stabilization between the reservoir and burden cells does not result in erroneous CO$_2$ plume results. Indeed, Figure \ref{fig:slice_HI_plume} shows a zoomed view of the CO$_2$-rich phase volume fraction in the same slice, and we see that even when the stabilization strength is increased the final CO$_2$ plume remains unchanged. With this example, however, we can show that the smeared pressure interface can affect predictions of the displacement of the burden. Figure \ref{fig:uplift_HI} shows the vertical displacement, or uplift, of the skeleton along the diagonal of the seafloor with both stabilization strategies. Stabilization across the interface results in increased uplift of the seafloor, and this effect worsens with increasing stabilization strength. Applying stabilization only in the burden remedies this and makes the displacement estimate agnostic to the stabilization strength. 

\section{Conclusions}

In this work we have considered the role of pressure stability within the fixed-stress sequential implicit method for coupled poromechanics. By appealing to examples from the incompressible flow literature, we clarified that it is splitting error which can result in splitting methods seemingly possessing pressure stability not attained by the equivalent fully implicit method, not simply the replacement of the discrete saddle point system with sequences of positive definite systems. This allowed us to reconcile results in the literature which seem to be in conflict regarding the pressure stability of the fixed-stress scheme. We also showed that, when applied on realistic simulations of CO$_2$ injection into formations with undrained burden regions, the iterative fixed-stress scheme produced the same oscillatory pressure solutions as the fully implicit method. The inclusion of global pressure jump stabilization effectively removed these oscillations and simultaneously greatly improved the convergence rate of the iterative method. 

Given the possibility of obtaining accurate overburden pressure fields using stabilization, the most attractive future direction for this work is the integration of contact modeling of faults and fractures within our poromechanical simulations. While the HI24L mesh used here did include faults, these were purely geometric factors in the results, meaning that there was no modeling of the actual mechanics of the faults. However, being able to accurately predict any fault reactivation or seismic activity resulting from stress increases as CO$_2$ is injected into the subsurface is critical for safety evaluations and monitoring. The interplay between stabilization schemes and the contact problem could be a very interesting avenue of research and is relatively unexplored.

\section*{Acknowledgements}
\label{sec::acknow}
Funding was provided by TotalEnergies and Chevron through the FC-MAELSTROM project.
Portions of this work were performed under the auspices of the U.S. Department of
Energy by Lawrence Livermore National Laboratory under Contract DE-AC52-07-NA27344 (LLNL-JRNL-851214).

\FloatBarrier


\begin{thebibliography}{10}
\expandafter\ifx\csname url\endcsname\relax
  \def\url#1{\texttt{#1}}\fi
\expandafter\ifx\csname urlprefix\endcsname\relax\def\urlprefix{URL }\fi
\expandafter\ifx\csname href\endcsname\relax
  \def\href#1#2{#2} \def\path#1{#1}\fi

\bibitem{ipcc}
S.~Benson, P.~Cook, IPCC Special Report on Carbon dioxide Capture and Storage,
  Cambridge University Press, 2005, Ch. Underground geological storage.

\bibitem{rutqvist2012geomechanics}
J.~Rutqvist, The geomechanics of co 2 storage in deep sedimentary formations,
  Geotechnical and Geological Engineering 30 (2012) 525--551.

\bibitem{zienkiewicz1990static}
O.~C. Zienkiewicz, A.~Chan, M.~Pastor, D.~Paul, T.~Shiomi, Static and dynamic
  behaviour of soils: a rational approach to quantitative solutions. i. fully
  saturated problems, Proceedings of the Royal Society of London. A.
  Mathematical and Physical Sciences 429~(1877) (1990) 285--309.

\bibitem{truty2006stabilized}
A.~Truty, T.~Zimmermann, Stabilized mixed finite element formulations for
  materially nonlinear partially saturated two-phase media, Computer methods in
  applied mechanics and engineering 195~(13-16) (2006) 1517--1546.

\bibitem{aguilar2008numerical}
G.~Aguilar, F.~Gaspar, F.~Lisbona, C.~Rodrigo, Numerical stabilization of
  biot's consolidation model by a perturbation on the flow equation,
  International journal for numerical methods in engineering 75~(11) (2008)
  1282--1300.

\bibitem{preisig2011stabilization}
M.~Preisig, J.~H. Pr{\'e}vost, Stabilization procedures in coupled
  poromechanics problems: A critical assessment, International Journal for
  Numerical and Analytical Methods in Geomechanics 35~(11) (2011) 1207--1225.

\bibitem{wan2003stabilized}
J.~Wan, Stabilized finite element methods for coupled geomechanics and
  multiphase flow, Stanford university, 2003.

\bibitem{white2008stabilized}
J.~A. White, R.~I. Borja, Stabilized low-order finite elements for coupled
  solid-deformation/fluid-diffusion and their application to fault zone
  transients, Computer Methods in Applied Mechanics and Engineering 197~(49-50)
  (2008) 4353--4366.

\bibitem{berger2015stabilized}
L.~Berger, R.~Bordas, D.~Kay, S.~Tavener, Stabilized lowest-order finite
  element approximation for linear three-field poroelasticity, SIAM Journal on
  Scientific Computing 37~(5) (2015) A2222--A2245.

\bibitem{camargo2021macroelement}
J.~T. Camargo, J.~A. White, R.~I. Borja, A macroelement stabilization for mixed
  finite element/finite volume discretizations of multiphase poromechanics,
  Computational Geosciences 25 (2021) 775--792.

\bibitem{frigo2021efficient}
M.~Frigo, N.~Castelletto, M.~Ferronato, J.~A. White, Efficient solvers for
  hybridized three-field mixed finite element coupled poromechanics, Computers
  \& Mathematics with Applications 91 (2021) 36--52.

\bibitem{garipov2018unified}
T.~T. Garipov, P.~Tomin, R.~Rin, D.~V. Voskov, H.~A. Tchelepi, Unified
  thermo-compositional-mechanical framework for reservoir simulation,
  Computational Geosciences 22 (2018) 1039--1057.

\bibitem{bui2020scalable}
Q.~M. Bui, D.~Osei-Kuffuor, N.~Castelletto, J.~A. White, A scalable multigrid
  reduction framework for multiphase poromechanics of heterogeneous media, SIAM
  Journal on Scientific Computing 42~(2) (2020) B379--B396.

\bibitem{bui2021multigrid}
Q.~M. Bui, F.~P. Hamon, N.~Castelletto, D.~Osei-Kuffuor, R.~R. Settgast, J.~A.
  White, Multigrid reduction preconditioning framework for coupled processes in
  porous and fractured media, Computer Methods in Applied Mechanics and
  Engineering 387 (2021) 114111.

\bibitem{dean2006comparison}
R.~H. Dean, X.~Gai, C.~M. Stone, S.~E. Minkoff, A comparison of techniques for
  coupling porous flow and geomechanics, Spe Journal 11~(01) (2006) 132--140.

\bibitem{park1983stabilization}
K.~Park, Stabilization of partitioned solution procedure for pore fluid-soil
  interaction analysis, International Journal for Numerical Methods in
  Engineering 19~(11) (1983) 1669--1673.

\bibitem{settari2001advances}
A.~Settari, D.~A. Walters, Advances in coupled geomechanical and reservoir
  modeling with applications to reservoir compaction, Spe Journal 6~(03) (2001)
  334--342.

\bibitem{mikelic2013convergence}
A.~Mikeli{\'c}, M.~F. Wheeler, Convergence of iterative coupling for coupled
  flow and geomechanics, Computational Geosciences 17 (2013) 455--461.

\bibitem{settari1998coupled}
A.~Settari, F.~Mounts, A coupled reservoir and geomechanical simulation system,
  Spe Journal 3~(03) (1998) 219--226.

\bibitem{kim2011stability}
J.~Kim, H.~A. Tchelepi, R.~Juanes, Stability and convergence of sequential
  methods for coupled flow and geomechanics: Fixed-stress and fixed-strain
  splits, Computer Methods in Applied Mechanics and Engineering 200~(13-16)
  (2011) 1591--1606.

\bibitem{kim2011stability_spe}
J.~Kim, H.~A. Tchelepi, R.~Juanes, Stability, accuracy, and efficiency of
  sequential methods for coupled flow and geomechanics, SPE Journal 16~(02)
  (2011) 249--262.

\bibitem{castelletto2015accuracy}
N.~Castelletto, J.~A. White, H.~Tchelepi, Accuracy and convergence properties
  of the fixed-stress iterative solution of two-way coupled poromechanics,
  International Journal for Numerical and Analytical Methods in Geomechanics
  39~(14) (2015) 1593--1618.

\bibitem{kim2011stability_drained}
J.~Kim, H.~A. Tchelepi, R.~Juanes, Stability and convergence of sequential
  methods for coupled flow and geomechanics: Drained and undrained splits,
  Computer Methods in Applied Mechanics and Engineering 200~(23-24) (2011)
  2094--2116.

\bibitem{janenko1971method}
N.~N. Janenko, The method of fractional steps, Vol. 160, Springer, 1971.

\bibitem{yoon2018spatial}
H.~C. Yoon, J.~Kim, Spatial stability for the monolithic and sequential methods
  with various space discretizations in poroelasticity, International Journal
  for Numerical Methods in Engineering 114~(7) (2018) 694--718.

\bibitem{storvik2019optimization}
E.~Storvik, J.~W. Both, K.~Kumar, J.~M. Nordbotten, F.~A. Radu, On the
  optimization of the fixed-stress splitting for biot's equations,
  International Journal for Numerical Methods in Engineering 120~(2) (2019)
  179--194.

\bibitem{aronson2023pressure}
R.~Aronson, F.~Hamon, N.~Castelletto, J.~White, H.~Tchelepi, Pressure jump
  stabilization for compositional poromechanics on unstructured meshes, in: SPE
  Reservoir Simulation Conference, SPE, 2023, p. D011S003R003.

\bibitem{storvik2018optimization}
E.~Storvik, J.~W. Both, K.~Kumar, J.~M. Nordbotten, F.~A. Radu, On the
  optimization of the fixed-stress splitting for biot's equations, arXiv
  preprint arXiv:1811.06242 (2018).

\bibitem{biot1941general}
M.~A. Biot, General theory of three-dimensional consolidation, Journal of
  applied physics 12~(2) (1941) 155--164.

\bibitem{coussy2004poromechanics}
O.~Coussy, Poromechanics, John Wiley \& Sons, 2004.

\bibitem{wang2000theory}
H.~Wang, Theory of linear poroelasticity with applications to geomechanics and
  hydrogeology, Vol.~2, Princeton university press, 2000.

\bibitem{murad1992improved}
M.~A. Murad, A.~F. Loula, Improved accuracy in finite element analysis of
  biot's consolidation problem, Computer Methods in Applied Mechanics and
  Engineering 95~(3) (1992) 359--382.

\bibitem{murad1994stability}
M.~A. Murad, A.~F. Loula, On stability and convergence of finite element
  approximations of biot's consolidation problem, International Journal for
  Numerical Methods in Engineering 37~(4) (1994) 645--667.

\bibitem{ferronato2010fully}
M.~Ferronato, N.~Castelletto, G.~Gambolati, A fully coupled 3-d mixed finite
  element model of biot consolidation, Journal of Computational Physics
  229~(12) (2010) 4813--4830.

\bibitem{nordbotten2014cell}
J.~M. Nordbotten, Cell-centered finite volume discretizations for deformable
  porous media, International journal for numerical methods in engineering
  100~(6) (2014) 399--418.

\bibitem{liu2004discontinuous}
R.~Liu, Discontinuous Galerkin finite element solution for poromechanics, The
  University of Texas at Austin, 2004.

\bibitem{choo2018enriched}
J.~Choo, S.~Lee, Enriched galerkin finite elements for coupled poromechanics
  with local mass conservation, Computer Methods in Applied Mechanics and
  Engineering 341 (2018) 311--332.

\bibitem{voskov2012comparison}
D.~V. Voskov, H.~A. Tchelepi, Comparison of nonlinear formulations for
  two-phase multi-component eos based simulation, Journal of Petroleum Science
  and Engineering 82 (2012) 101--111.

\bibitem{aziz1979petroleum}
K.~Aziz, Petroleum reservoir simulation, Applied Science Publishers 476 (1979).

\bibitem{white2019two}
J.~A. White, N.~Castelletto, S.~Klevtsov, Q.~M. Bui, D.~Osei-Kuffuor, H.~A.
  Tchelepi, A two-stage preconditioner for multiphase poromechanics in
  reservoir simulation, Computer Methods in Applied Mechanics and Engineering
  357 (2019) 112575.

\bibitem{terekhov2017cell}
K.~M. Terekhov, B.~T. Mallison, H.~A. Tchelepi, Cell-centered nonlinear
  finite-volume methods for the heterogeneous anisotropic diffusion problem,
  Journal of Computational Physics 330 (2017) 245--267.

\bibitem{both2019numerical}
J.~W. Both, U.~K{\"o}cher, Numerical investigation on the fixed-stress
  splitting scheme for biot’s equations: Optimality of the tuning parameter,
  in: Numerical Mathematics and Advanced Applications ENUMATH 2017, Springer,
  2019, pp. 789--797.

\bibitem{guermond1998stability}
J.-L. Guermond, L.~Quartapelle, On stability and convergence of projection
  methods based on pressure poisson equation, International Journal for
  Numerical Methods in Fluids 26~(9) (1998) 1039--1053.

\bibitem{guermond2006overview}
J.-L. Guermond, P.~Minev, J.~Shen, An overview of projection methods for
  incompressible flows, Computer methods in applied mechanics and engineering
  195~(44-47) (2006) 6011--6045.

\bibitem{chorin1968numerical}
A.~J. Chorin, Numerical solution of the navier-stokes equations, Mathematics of
  computation 22~(104) (1968) 745--762.

\bibitem{chorin1969convergence}
A.~J. Chorin, On the convergence of discrete approximations to the
  navier-stokes equations, Mathematics of computation 23~(106) (1969) 341--353.

\bibitem{badia2007convergence}
S.~Badia, R.~Codina, Convergence analysis of the fem approximation of the first
  order projection method for incompressible flows with and without the inf-sup
  condition, Numerische Mathematik 107 (2007) 533--557.

\bibitem{goda1979multistep}
K.~Goda, A multistep technique with implicit difference schemes for calculating
  two-or three-dimensional cavity flows, Journal of computational physics
  30~(1) (1979) 76--95.

\bibitem{barry1999exact}
S.~Barry, G.~Mercer, Exact solutions for two-dimensional time-dependent flow
  and deformation within a poroelastic medium, Journal of applied mechanics
  66~(2) (1999) 536--540.

\bibitem{hughes_1987}
T.~J. Hughes, L.~P. Franca,
  \href{https://www.sciencedirect.com/science/article/pii/0045782587901848}{A
  new finite element formulation for computational fluid dynamics: Vii. the
  stokes problem with various well-posed boundary conditions: Symmetric
  formulations that converge for all velocity/pressure spaces}, Computer
  Methods in Applied Mechanics and Engineering 65~(1) (1987) 85--96.
\newblock \href {https://doi.org/https://doi.org/10.1016/0045-7825(87)90184-8}
  {\path{doi:https://doi.org/10.1016/0045-7825(87)90184-8}}.
\newline\urlprefix\url{https://www.sciencedirect.com/science/article/pii/0045782587901848}

\bibitem{silvester1990stabilised}
D.~J. Silvester, N.~Kechkar, Stabilised bilinear-constant velocity-pressure
  finite elements for the conjugate gradient solution of the stokes problem,
  Computer Methods in Applied Mechanics and Engineering 79~(1) (1990) 71--86.

\bibitem{silvester1994optimal}
D.~Silvester, Optimal low order finite element methods for incompressible flow,
  Computer methods in applied mechanics and engineering 111~(3-4) (1994)
  357--368.

\bibitem{borio2021hybrid}
A.~Borio, F.~P. Hamon, N.~Castelletto, J.~A. White, R.~R. Settgast, Hybrid
  mimetic finite-difference and virtual element formulation for coupled
  poromechanics, Computer Methods in Applied Mechanics and Engineering 383
  (2021) 113917.

\bibitem{berger2017stabilized}
L.~Berger, R.~Bordas, D.~Kay, S.~Tavener, A stabilized finite element method
  for finite-strain three-field poroelasticity, Computational Mechanics 60
  (2017) 51--68.

\bibitem{dohrmann2004stabilized}
C.~R. Dohrmann, P.~B. Bochev, A stabilized finite element method for the stokes
  problem based on polynomial pressure projections, International Journal for
  Numerical Methods in Fluids 46~(2) (2004) 183--201.

\bibitem{t2022deformation}
J.~T~Camargo, F.~Hamon, A.~Mazuyer, T.~Meckel, N.~Castelletto, J.~A. White,
  Deformation monitoring feasibility for offshore carbon storage in the
  gulf-of-mexico, Fran{\c{c}}ois and Mazuyer, Antoine and Meckel, Tip and
  Castelletto, Nicola and White, Joshua A., Deformation Monitoring Feasibility
  for Offshore Carbon Storage in the Gulf-of-Mexico (October 23, 2022) (2022).

\bibitem{GEOS}
GEOS, {Next-gen simulation for geologic carbon storage},
  \url{http://www.geos.dev} (2023).

\bibitem{duan2003improved}
Z.~Duan, R.~Sun, An improved model calculating co2 solubility in pure water and
  aqueous nacl solutions from 273 to 533 k and from 0 to 2000 bar, Chemical
  geology 193~(3-4) (2003) 257--271.

\bibitem{span1996new}
R.~Span, W.~Wagner, A new equation of state for carbon dioxide covering the
  fluid region from the triple-point temperature to 1100 k at pressures up to
  800 mpa, Journal of physical and chemical reference data 25~(6) (1996)
  1509--1596.

\bibitem{fenghour1998viscosity}
A.~Fenghour, W.~A. Wakeham, V.~Vesovic, The viscosity of carbon dioxide,
  Journal of physical and chemical reference data 27~(1) (1998) 31--44.

\bibitem{phillips1981technical}
S.~Phillips, A technical databook for geothermal energy utilization (1981).

\bibitem{ruiz2019characterization}
I.~Ruiz, Characterization of the high island 24l field for modeling and
  estimating co2 storage capacity in the offshore texas state waters, gulf of
  mexico, Ph.D. thesis (2019).

\end{thebibliography}

\end{document}